\documentclass{amsart}
\usepackage{amssymb}

\usepackage{amscd}
\usepackage{thmdefs}

\input{tcilatex}
\theoremstyle{definition}
\theoremstyle{remark}
\numberwithin{equation}{section}

\input tcilatex

\begin{document}
\title[R-transform Theory On a Tower of Amalgamated NCPSpaces]{R-transform Theory On a Tower of Amalgamated Noncommutative Probability
Spaces }
\author{Ilwoo Cho}
\address{Dep. of Math, Univ of Iowa, Iowa Sity, IA, U.S.A}
\email{icho@math.uiowa.edu}
\thanks{}
\maketitle

\begin{abstract}
In [12], we observed Amalgamated R-transform Theory. Different from the
original definition of Voiculescu and Speicher (in [10] and [1]), we define
R-transforms of operator-valued random variable(s) by operator-valued formal
series. By doing that we can establish Amalgamated R-transform calculus. In
this paper, we will concentrate on observing R-transforms on a tower of
amalgamated noncommutative probability spaces. We will consider the
R-transform calculus and relation between R-transforms on the tower and
distributions with respect to the tower. As application, we will observe the
compatibility of the given tower and commuting ladders of amalgamated
noncommutative probability spaces.
\end{abstract}

\strut

Voiculescu developed Free Probability Theory. Here, the classical concept of
Independence in Probability theory is replaced by a noncommutative analogue
called Freeness (See [8]). There are two approaches to study Free
Probability Theory. One of them is the original analytic approach of
Voiculescu (See [8] and [7]) and the other one is the combinatorial approach
of Speicher and Nica (See [1], [6] and [7]).

To observe the free additive convolution and free multiplicative convolution
of two distributions induced by free random variables in a noncommutative
probability space (over $B=\Bbb{C}$), Voiculescu defined R-transforms and
S-transforms, respectively. These show that to study distributions is to
study certain ($B$-) formal series for arbitrary noncommutative
indeterminants.

Speicher defined the free cumulants which are the main objects in the
combinatorial approach of Free Probability Theory. And he developed free
probability theory by using the combinatorics and lattice theory on
collections of noncrossing partitions (See [7]). Also, Speicher considered
the operator-valued free probability theory, which is also defined and
observed analytically by Voiculescu, when $\Bbb{C}$ is replaced to an
arbitrary algebra $B$ (See [1] and [7]). Nica defined R-transforms of
several random variables (See [6]). He defined these R-transforms as
multivariable formal series in noncommutative several indeterminants. To
observe the R-transform, the M\"{o}bius Inversion under the embedding of
lattices plays a key role (See [1],[7],[5],[9],[13] and [17]).\strut 

In [9], [19] and [20], we observed the amalgamated R-transform calculus.
Actually, amalgamated R-transforms are defined originally by Voiculescu (See
[7]) and are characterized combinatorially by Speicher (See [1]). In [9], we
defined amalgamated R-transforms slightly different, compared with those in
[1] and [7], because of the bimodule map property. We defined them as $B$%
-formal series, in $\Theta _{B}^{s},$ in $s$-noncommutative indeterminents
and tried to characterize, like in [6] and [7]. The main tool which is
considered, for studying amalgamated R-transform calculus, is the
amalgamated boxed convolution, \frame{*}$_{B}$ defined in [9]. After
defining boxed convolution over an arbitrary algebra $B,$ we could get that

\strut

\begin{center}
$R_{x_{1},...,x_{s}}\,\,\frame{*}_{B}\,%
\,R_{y_{1},...,y_{s}}^{symm(1_{B})}=R_{x_{1}y_{1},...,x_{s}y_{s}},$ for any $%
s\in \Bbb{N},$
\end{center}

\strut

where $x_{j}$'s and $y_{j}$'s are free $B$-valued random variables. However,
different from scalar-valued case (in [6] and [7]), in general,

\strut

\begin{center}
$R_{x_{1},...,x_{s}}\,\,\,\frame{*}_{B}\,\,\,R_{y_{1},...,y_{s}}\neq
R_{x_{1}y_{1},...,x_{s}y_{s}},$
\end{center}

\strut

even if $x_{j}$'s and $y_{j}$'s are free over $B.$

\strut

In this paper, we will consider a tower of algebras. i.e

\strut

\begin{center}
$\Bbb{C}\subset B\subset A_{1}\subset A_{2}\subset ....\subset A_{\infty },$
\end{center}

\strut

where $A_{1}$ is an algebra over $B$ and $A_{j+1}$ is an algebra over $%
A_{j}, $ for all $j=1,2,\cdot \cdot \cdot $, and $A_{\infty }$ is an
enveloping algebra of this tower of algebra. Notice that $A_{n}$ is an
algebra over $B,$ again, for all $n\in \Bbb{N}.$

\strut

In Chapter 1, we will observe Amalgamated R-transform Theory established in
[9]. In Chapter 2, we will consider the R-transform theory on a tower of
algebras. In Chapter 3, we will study the relation between R-transforms on a
tower of algebras and operator-valued (or scalar-valued) distributions on
that tower. In Chapter 4, as application, we will define a operator-valued
semicircularity on a tower of algebras and observe their R-transform theory.
Also, after assuming the compatibility on the tower, we will observe the
R-transform calculus and operator-valued distributions.

\strut

\strut

\strut

\section{Preliminaries}

\strut

\strut

\subsection{Amalgamated Free Probability Theory}

\strut

\strut

In this section, we will summarize and introduced the basic results from [1]
and [9]. Throughout this section, let $B$ be a unital algebra. The algebraic
pair $(A,\varphi )$ is said to be a noncommutative probability space over $B$
(shortly, NCPSpace over $B$) if $A$ is an algebra over $B$ (i.e $%
1_{B}=1_{A}\in B\subset A$) and $\varphi :A\rightarrow B$ is a $B$%
-functional (or a conditional expectation) ; $\varphi $ satisfies

\begin{center}
$\varphi (b)=b,$ for all $b\in B$
\end{center}

and

\begin{center}
$\varphi (bxb^{\prime })=b\varphi (x)b^{\prime },$ for all $b,b^{\prime }\in
B$ and $x\in A.$
\end{center}

\strut

\strut Let $(A,\varphi )$ be a NCPSpace over $B.$ Then, for the given $B$%
-functional, we can determine a moment multiplicative function $\widehat{%
\varphi }=(\varphi ^{(n)})_{n=1}^{\infty }\in I(A,B),$ where

\strut

\begin{center}
$\varphi ^{(n)}(a_{1}\otimes ...\otimes a_{n})=\varphi (a_{1}....a_{n}),$
\end{center}

\strut \strut

for all $a_{1}\otimes ...\otimes a_{n}\in A^{\otimes _{B}n},$ $\forall n\in 
\Bbb{N}.$

\strut

\strut We will denote noncrossing partitions over $\{1,...,n\}$ ($n\in \Bbb{N%
}$) by $NC(n).$ Define an ordering on $NC(n)$ ;

\strut

$\theta =\{V_{1},...,V_{k}\}\leq \pi =\{W_{1},...,W_{l}\}$\strut $\overset{%
def}{\Leftrightarrow }$ For each block $V_{j}\in \theta $, there exists only
one block $W_{p}\in \pi $ such that $V_{j}\subset W_{p},$ for $j=1,...,k$
and $p=1,...,l.$

\strut

Then $(NC(n),\leq )$ is a complete lattice with its minimal element $%
0_{n}=\{(1),...,(n)\}$ and its maximal element $1_{n}=\{(1,...,n)\}$. We
define the incidence algebra $I_{2}$ by a set of all complex-valued\
functions $\eta $ on $\cup _{n=1}^{\infty }\left( NC(n)\times NC(n)\right) $
satisfying $\eta (\theta ,\pi )=0,$ whenever $\theta \nleq \pi .$ Then,
under the convolution

\begin{center}
$*:I_{2}\times I_{2}\rightarrow \Bbb{C}$
\end{center}

defined by

\begin{center}
$\eta _{1}*\eta _{2}(\theta ,\pi )=\underset{\theta \leq \sigma \leq \pi }{%
\sum }\eta _{1}(\theta ,\sigma )\cdot \eta _{2}(\sigma ,\pi ),$
\end{center}

$I_{2}$ is indeed an algebra of complex-valued functions. Denote zeta,
M\"{o}bius and delta functions in the incidence algebra $I_{2}$ by $\zeta ,$ 
$\mu $ and $\delta ,$ respectively. i.e

\strut

\begin{center}
$\zeta (\theta ,\pi )=\left\{ 
\begin{array}{lll}
1 &  & \theta \leq \pi \\ 
0 &  & otherwise,
\end{array}
\right. $
\end{center}

\strut

\begin{center}
$\delta (\theta ,\pi )=\left\{ 
\begin{array}{lll}
1 &  & \theta =\pi \\ 
0 &  & otherwise,
\end{array}
\right. $
\end{center}

\strut

and $\mu $ is the ($*$)-inverse of $\zeta .$ Notice that $\delta $ is the ($%
* $)-identity of $I_{2}.$ By using the same notation ($*$), we can define a
convolution between $I(A,B)$ and $I_{2}$ by

\strut

\begin{center}
$\widehat{f}\,*\,\eta \left( a_{1},...,a_{n}\,;\,\pi \right) =\underset{\pi
\in NC(n)}{\sum }\widehat{f}(\pi )(a_{1}\otimes ...\otimes a_{n})\eta (\pi
,1_{n}),$
\end{center}

\strut

where $\widehat{f}\in I(A,B)$, $\eta \in I_{1},$ $\pi \in NC(n)$ and $%
a_{j}\in A$ ($j=1,...,n$), for all $n\in \Bbb{N}.$ Notice that $\widehat{f}%
*\eta \in I(A,B),$ too. Let $\widehat{\varphi }$ be a moment multiplicative
function in $I(A,B)$ which we determined before. Then we can naturally
define a cumulant multiplicative function $\widehat{c}=(c^{(n)})_{n=1}^{%
\infty }\in I(A,B)$ by

\begin{center}
$\widehat{c}=\widehat{\varphi }*\mu $ \ \ \ or \ \ $\widehat{\varphi }=%
\widehat{c}*\zeta .$
\end{center}

This says that if we have a moment\ multiplicative function, then we always
get a cumulant multiplicative function and vice versa, by $(*).$ This
relation is so-called ''M\"{o}bius Inversion''. More precisely, we have

\strut

\begin{center}
$
\begin{array}{ll}
\varphi (a_{1}...a_{n}) & =\varphi ^{(n)}(a_{1}\otimes ...\otimes a_{n}) \\ 
& =\underset{\pi \in NC(n)}{\sum }\widehat{c}(\pi )(a_{1}\otimes ...\otimes
a_{n})\zeta (\pi ,1_{n}) \\ 
& =\underset{\pi \in NC(n)}{\sum }\widehat{c}(\pi )(a_{1}\otimes ...\otimes
a_{n}),
\end{array}
$
\end{center}

\strut

for all $a_{j}\in A$ and $n\in \Bbb{N}.$ Or equivalently,

\strut

\begin{center}
$
\begin{array}{ll}
c^{(n)}(a_{1}\otimes ...\otimes a_{n}) & =\underset{\pi \in NC(n)}{\sum }%
\widehat{\varphi }(\pi )(a_{1}\otimes ...\otimes a_{n})\mu (\pi ,1_{n}).
\end{array}
$
\end{center}

\strut \strut

Now, let $(A_{i},\varphi _{i})$ be NCPSpaces over $B,$ for all $i\in I.$
Then we can define a amalgamated free product of $A_{i}$ 's and amalgamated
free product of $\varphi _{i}$'s by

\begin{center}
$A\equiv *_{B}A_{i}$ \ \ and \ $\varphi \equiv *_{i}\varphi _{i},$
\end{center}

respectively. Then, by Voiculescu, $(A,\varphi )$ is again a NCPSpace over $%
B $ and, as a vector space, $A$ can be represented by

\begin{center}
\strut $A=B\oplus \left( \oplus _{n=1}^{\infty }\left( \underset{i_{1}\neq
...\neq i_{n}}{\oplus }(A_{i_{1}}\ominus B)\otimes ...\otimes
(A_{i_{n}}\ominus B)\right) \right) ,$
\end{center}

where $A_{i_{j}}\ominus B=\ker \varphi _{i_{j}}.$ We will use Speicher's
combinatorial definition of amalgamated free product of $B$-functionals ;

\strut

\begin{definition}
Let $(A_{i},\varphi _{i})$ be NCPSpaces over $B,$ for all $i\in I.$ Then $%
\varphi =*_{i}\varphi _{i}$ is the amalgamated free product of $B$%
-functionals $\varphi _{i}$'s on $A=*_{B}A_{i}$ if the cumulant
multiplicative function $\widehat{c}=\widehat{\varphi }*\mu \in I(A,B)$ has
its restriction to $\underset{i\in I}{\cup }A_{i},$ $\underset{i\in I}{%
\oplus }\widehat{c_{i}},$ where $\widehat{c_{i}}$ is the cumulant
multiplicative function induced by $\varphi _{i},$ for all $i\in I$ and, for
each $n\in \Bbb{N},$

\strut 

\begin{center}
$c^{(n)}(a_{1}\otimes ...\otimes a_{n})=\left\{ 
\begin{array}{lll}
c_{i}^{(n)}(a_{1}\otimes ...\otimes a_{n}) &  & \text{if }\forall a_{j}\in
A_{i} \\ 
0_{B} &  & otherwise.
\end{array}
\right. $
\end{center}
\end{definition}

\strut

Now, we will observe the freeness over $B$ ;

\strut

\begin{definition}
Let $(A,\varphi )$\strut be a NCPSpace over $B.$

\strut 

(1) Subalgebras containing $B,$ $A_{i}\subset A$ ($i\in I$) are free (over $B
$) if we let $\varphi _{i}=\varphi \mid _{A_{i}},$ for all $i\in I,$ then $%
*_{i}\varphi _{i}$ has its cumulant multiplicative function $\widehat{c}$
such that its restriction to $\underset{i\in I}{\cup }A_{i}$ is $\underset{%
i\in I}{\oplus }\widehat{c_{i}},$ where $\widehat{c_{i}}$ is the cumulant
multiplicative function induced by each $\varphi _{i},$ for all $i\in I.$

\strut 

(2) Sebsets $X_{i}$ ($i\in I$) are free (over $B$) if subalgebras $A_{i}$'s
generated by $B$ and $X_{i}$'s are free in the sense of (1). i.e If we let $%
A_{i}=A\lg \left( X_{i},B\right) ,$ for all $i\in I,$ then $A_{i}$'s are
free over $B.$
\end{definition}

\strut

In [1], Speicher showed that the above combinatorial freeness with
amalgamation can be used alternatively with respect to Voiculescu's original
freeness with amalgamation.

\strut

Let $(A,\varphi )$ be a NCPSpace over $B$ and let $x_{1},...,x_{s}$ be $B$%
-valued random variables ($s\in \Bbb{N}$). Define $(i_{1},...,i_{n})$-th
moment of $x_{1},...,x_{s}$ by

\strut

\begin{center}
$\varphi (x_{i_{1}}b_{i_{2}}x_{i_{2}}...b_{i_{n}}x_{i_{n}}),$
\end{center}

\strut

for arbitrary $b_{i_{2}},...,b_{i_{n}}\in B,$ where $(i_{1},...,i_{n})\in
\{1,...,s\}^{n},$ $\forall n\in \Bbb{N}.$ Similarly, define a symmetric $%
(i_{1},...,i_{n})$-th moment by the fixed $b_{0}\in B$ by

\strut

\begin{center}
$\varphi (x_{i_{1}}b_{0}x_{i_{2}}...b_{0}x_{i_{n}}).$
\end{center}

\strut

If $b_{0}=1_{B},$ then we call this symmetric moments, trivial moments.

\strut

Cumulants defined below are main tool of combinatorial free probability
theory ; in [9], we defined the $(i_{1},...,i_{n})$-th cumulant of $%
x_{1},...,x_{s}$ by

\strut

\begin{center}
$k_{n}(x_{i_{1}},...,x_{i_{n}})=c^{(n)}(x_{i_{1}}\otimes
b_{i_{2}}x_{i_{2}}\otimes ...\otimes b_{i_{n}}x_{i_{n}}),$
\end{center}

\strut

for $b_{i_{2}},...,b_{i_{n}}\in B,$ arbitrary, and $(i_{1},...,i_{n})\in
\{1,...,s\}^{n},$ $\forall n\in \Bbb{N},$ where $\widehat{c}%
=(c^{(n)})_{n=1}^{\infty }$ is the cumulant multiplicative function induced
by $\varphi .$ Notice that, by M\"{o}bius inversion, we can always take such 
$B$-value whenever we have $(i_{1},...,i_{n})$-th moment of $%
x_{1},...,x_{s}. $ And, vice versa, if we have cumulants, then we can always
take moments. Hence we can define a symmetric $(i_{1},...,i_{n})$-th\
cumulant by $b_{0}\in B$ of $x_{1},...,x_{s}$ by

\strut

\begin{center}
$k_{n}^{symm(b_{0})}(x_{i_{1}},...,x_{i_{n}})=c^{(n)}(x_{i_{1}}\otimes
b_{0}x_{i_{2}}\otimes ...\otimes b_{0}x_{i_{n}}).$
\end{center}

\strut

If $b_{0}=1_{B},$ then it is said to be trivial cumulants of $%
x_{1},...,x_{s} $.

\strut

By Speicher, it is shown that subalgebras $A_{i}$ ($i\in I$) are free over $%
B $ if and only if all mixed cumulants vanish.

\strut

\begin{proposition}
(See [1] and [9]) Let $(A,\varphi )$ be a NCPSpace over $B$ and let $%
x_{1},...,x_{s}\in (A,\varphi )$ be $B$-valued random variables ($s\in \Bbb{N%
}$). Then $x_{1},...,x_{s}$ are free if and only if all their mixed
cumulants vanish. $\square $
\end{proposition}

\strut

\strut

\strut

\strut

\subsection{Amalgamated R-transform Theory}

\strut

\strut

In this section, we will define an R-transform of several $B$-valued random
variables. Note that to study R-transforms is to study operator-valued
distributions. R-transforms with single variable is defined by Voiculescu
(over $B,$ in particular, $B=\Bbb{C}$. See [8] and [7]). Over $\Bbb{C},$
Nica defined multi-variable R-transforms in [6]. In [9], we extended his
concepts, over $B.$ R-transforms of $B$-valued random variables can be
defined as $B$-formal series with its $(i_{1},...,i_{n})$-th coefficients, $%
(i_{1},...,i_{n})$-th cumulants of $B$-valued random variables, where $%
(i_{1},...,i_{n})\in \{1,...,s\}^{n},$ $\forall n\in \Bbb{N}.$

\strut

\begin{definition}
Let $(A,\varphi )$ be a NCPSpace over $B$ and let $x_{1},...,x_{s}\in
(A,\varphi )$ be $B$-valued random variables ($s\in \Bbb{N}$). Let $%
z_{1},...,z_{s}$ be noncommutative indeterminants. Define a moment series of 
$x_{1},...,x_{s}$, as a $B$-formal series, by

\strut 

\begin{center}
$M_{x_{1},...,x_{s}}(z_{1},...,z_{s})=\sum_{n=1}^{\infty }\underset{%
i_{1},..,i_{n}\in \{1,...,s\}}{\sum }\varphi
(x_{i_{1}}b_{i_{2}}x_{i_{2}}...b_{i_{n}}x_{i_{n}})\,z_{i_{1}}...z_{i_{n}},$
\end{center}

\strut 

where $b_{i_{2}},...,b_{i_{n}}\in B$ are arbitrary for all $%
(i_{2},...,i_{n})\in \{1,...,s\}^{n-1},$ $\forall n\in \Bbb{N}.$

\strut 

Define an R-transform of $x_{1},...,x_{s}$, as a $B$-formal series, by

\strut 

\begin{center}
$R_{x_{1},...,x_{s}}(z_{1},...,z_{s})=\sum_{n=1}^{\infty }\underset{%
i_{1},...,i_{n}\in \{1,...,s\}}{\sum }k_{n}(x_{i_{1}},...,x_{i_{n}})%
\,z_{i_{1}}...z_{i_{n}},$
\end{center}

\strut with

\begin{center}
$k_{n}(x_{i_{1}},...,x_{i_{n}})=c^{(n)}(x_{i_{1}}\otimes
b_{i_{2}}x_{i_{2}}\otimes ...\otimes b_{i_{n}}x_{i_{n}}),$
\end{center}

\strut 

where $b_{i_{2}},...,b_{i_{n}}\in B$ are arbitrary for all $%
(i_{2},...,i_{n})\in \{1,...,s\}^{n-1},$ $\forall n\in \Bbb{N}.$ Here, $%
\widehat{c}=(c^{(n)})_{n=1}^{\infty }$ is a cumulant multiplicative function
induced by $\varphi $ in $I(A,B).$
\end{definition}

\strut

Denote a set of all $B$-formal series with $s$-noncommutative indeterminants
($s\in \Bbb{N}$), by $\Theta _{B}^{s}$. i.e if $g\in \Theta _{B}^{s},$ then

\begin{center}
$g(z_{1},...,z_{s})=\sum_{n=1}^{\infty }\underset{i_{1},...,i_{n}\in
\{1,...,s\}}{\sum }b_{i_{1},...,i_{n}}\,z_{i_{1}}...z_{i_{n}},$
\end{center}

where $b_{i_{1},...,i_{n}}\in B,$ for all $(i_{1},...,i_{n})\in
\{1,...,s\}^{n},$ $\forall n\in \Bbb{N}.$ Trivially, by definition, $%
M_{x_{1},...,x_{s}},$ $R_{x_{1},...,x_{s}}\in \Theta _{B}^{s}.$ By $\mathcal{%
R}_{B}^{s},$\ we denote a set of all R-transforms of $s$-$B$-valued random
variables. Recall that, set-theoratically,

\begin{center}
$\Theta _{B}^{s}=\mathcal{R}_{B}^{s},$ sor all $s\in \Bbb{N}.$
\end{center}

\strut

We can also define symmetric moment series and symmetric R-transform by $%
b_{0}\in B,$ by

\strut

\begin{center}
$M_{x_{1},...,x_{s}}^{symm(b_{0})}(z_{1},...,z_{s})=\sum_{n=1}^{\infty }%
\underset{i_{1},...,i_{n}\in \{1,...,s\}}{\sum }\varphi
(x_{i_{1}}b_{0}x_{i_{2}}...b_{0}x_{i_{n}})\,z_{i_{1}}...z_{i_{n}}$
\end{center}

and

\begin{center}
$R_{x_{1},...,x_{s}}^{symm(b_{0})}(z_{1},...,z_{s})=\sum_{n=1}^{\infty }%
\underset{i_{1},..,i_{n}\in \{1,...,s\}}{\sum }%
k_{n}^{symm(b_{0})}(x_{i_{1}},...,x_{i_{n}})\,z_{i_{1}}...z_{i_{n}},$
\end{center}

with

\begin{center}
$k_{n}^{symm(b_{0})}(x_{i_{1}},...,x_{i_{n}})=c^{(n)}(x_{i_{1}}\otimes
b_{0}x_{i_{2}}\otimes ...\otimes b_{0}x_{i_{n}}),$
\end{center}

for all $(i_{1},...,i_{n})\in \{1,...,s\}^{n},$ $\forall n\in \Bbb{N}.$

\strut

If $b_{0}=1_{B},$ then we have trivial moment series and trivial R-transform
of $x_{1},...,x_{s}$ denoted by $M_{x_{1},...,x_{s}}^{t}$ and $%
R_{x_{1},...,x_{s}}^{t},$ respectively.

\strut

The followings are known in [1] and [9] ;

\strut

\begin{proposition}
Let $(A,\varphi )$ be a NCPSpace over $B$ and let $%
x_{1},...,x_{s},y_{1},...,y_{p}\in (A,\varphi )$ be $B$-valued random
variables, where $s,p\in \Bbb{N}.$ Suppose that $\{x_{1},...,x_{s}\}$ and $%
\{y_{1},...,y_{p}\}$ are free in $(A,\varphi ).$ Then

\strut 

(1) $
R_{x_{1},...,x_{s},y_{1},...,y_{p}}(z_{1},...,z_{s+p})=R_{x_{1},...,x_{s}}(z_{1},...,z_{s})+R_{y_{1},...,y_{p}}(z_{1},...,z_{p}).
$

\strut 

(2) If $s=p,$ then $R_{x_{1}+y_{1},...,x_{s}+y_{s}}(z_{1},...,z_{s})=\left(
R_{x_{1},...,x_{s}}+R_{y_{1},...,y_{s}}\right) (z_{1},...,z_{s}).$

$\square $
\end{proposition}

\strut

The above proposition is proved by the characterization of freeness with
respect to cumulants. i.e $\{x_{1},...,x_{s}\}$ and $\{y_{1},...,y_{p}\}$
are free in $(A,\varphi )$ if and only if their mixed cumulants vanish. Thus
we have

\strut

$k_{n}(p_{i_{1}},...,p_{i_{n}})=c^{(n)}(p_{i_{1}}\otimes
b_{i_{2}}p_{i_{2}}\otimes ...\otimes b_{i_{n}}p_{i_{n}})$

$\ \ \ \ =\left( \widehat{c_{x}}\oplus \widehat{c_{y}}\right)
^{(n)}(p_{i_{1}}\otimes b_{i_{2}}p_{i_{2}}\otimes ...\otimes
b_{i_{n}}p_{i_{n}})$

$\ \ \ \ =\left\{ 
\begin{array}{lll}
k_{n}(x_{i_{1}},...,x_{i_{n}}) &  & or \\ 
k_{n}(y_{i_{1}},...,y_{i_{n}}) &  & 
\end{array}
\right. $

\strut

and if $s=p,$ then

$k_{n}(x_{i_{1}}+y_{i_{1}},...,x_{i_{n}}+y_{i_{n}})$

$\ =c^{(n)}\left( (x_{i_{1}}+y_{i_{1}})\otimes
b_{i_{2}}(x_{i_{2}}+y_{i_{2}})\otimes ...\otimes
b_{i_{n}}(x_{i_{n}}+y_{i_{n}})\right) $

$\ \ =c^{(n)}(x_{i_{1}}\otimes b_{i_{2}}x_{i_{2}}\otimes ...\otimes
b_{i_{n}}x_{i_{n}})+c^{(n)}(y_{i_{1}}\otimes b_{i_{2}}y_{i_{2}}\otimes
...\otimes b_{i_{n}}y_{i_{n}})+[Mixed]$

\strut

where $[Mixed]$ is the sum of mixed cumulants of $x_{j}$'s and $y_{i}$'s, by
the bimodule map property of $c^{(n)}$

\strut

$\ =k_{n}(x_{i_{1}},...,x_{i_{n}})+k_{n}(y_{i_{1}},...,y_{i_{n}})+0_{B}.$

\strut

Note that if $f,g\in \Theta _{B}^{s},$ then we can always choose free $%
\{x_{1},...,x_{s}\}$ and $\{y_{1},...,y_{s}\}$ in (some) NCPSpace over $B,$ $%
(A,\varphi ),$ such that

\begin{center}
$f=R_{x_{1},...,x_{s}}$ \ \ and \ \ $g=R_{y_{1},...,y_{s}}.$
\end{center}

\strut

\begin{definition}
(1) Let $s\in \Bbb{N}.$ Let $(f,g)\in \Theta _{B}^{s}\times \Theta _{B}^{s}.$
Define \frame{*}\thinspace \thinspace $:\Theta _{B}^{s}\times \Theta
_{B}^{s}\rightarrow \Theta _{B}^{s}$ by

\strut 

\begin{center}
$\left( f,g\right) =\left( R_{x_{1},...,x_{s}},\,R_{y_{1},...,y_{s}}\right)
\longmapsto R_{x_{1},...,x_{s}}\,\,\frame{*}\,\,R_{y_{1},...,y_{s}}.$
\end{center}

\strut 

Here, $\{x_{1},...,x_{s}\}$ and $\{y_{1},...,y_{s}\}$ are free in $%
(A,\varphi )$. Suppose that

\strut 

\begin{center}
$coef_{i_{1},..,i_{n}}\left( R_{x_{1},...,x_{s}}\right)
=c^{(n)}(x_{i_{1}}\otimes b_{i_{2}}x_{i_{2}}\otimes ...\otimes
b_{i_{n}}x_{i_{n}})$
\end{center}

and

\begin{center}
$coef_{i_{1},...,i_{n}}(R_{y_{1},...,y_{s}})=c^{(n)}(y_{i_{1}}\otimes
b_{i_{2}}^{\prime }y_{i_{2}}\otimes ...\otimes b_{i_{n}}^{\prime }y_{i_{n}}),
$
\end{center}

\strut 

for all $(i_{1},...,i_{n})\in \{1,...,s\}^{n},$ $n\in \Bbb{N},$ where $%
b_{i_{j}},b_{i_{n}}^{\prime }\in B$ arbitrary. Then

\strut 

$coef_{i_{1},...,i_{n}}\left( R_{x_{1},...,x_{s}}\,\,\frame{*}%
\,\,R_{y_{1},...,y_{s}}\right) $

\strut 

$=\underset{\pi \in NC(n)}{\sum }\left( \widehat{c_{x}}\oplus \widehat{c_{y}}%
\right) (\pi \cup Kr(\pi ))(x_{i_{1}}\otimes y_{i_{1}}\otimes
b_{i_{2}}x_{i_{2}}\otimes b_{i_{2}}^{\prime }y_{i_{2}}\otimes ...\otimes
b_{i_{n}}x_{i_{n}}\otimes b_{i_{n}}^{\prime }y_{i_{n}})$

$\strut $

$\overset{denote}{=}\underset{\pi \in NC(n)}{\sum }\left( k_{\pi }\oplus
k_{Kr(\pi )}\right) (x_{i_{1}},y_{i_{1}},...,x_{i_{n}}y_{i_{n}}),$

\strut \strut 

where $\widehat{c_{x}}\oplus \widehat{c_{y}}=\widehat{c}\mid
_{A_{x}*_{B}A_{y}},$ $A_{x}=A\lg \left( \{x_{i}\}_{i=1}^{s},B\right) $ and $%
A_{y}=A\lg \left( \{y_{i}\}_{i=1}^{s},B\right) $ and where $\pi \cup Kr(\pi )
$ is an alternating union of partitions in $NC(2n)$
\end{definition}

\strut

\begin{proposition}
(See [9])\strut Let $(A,\varphi )$ be a NCPSpace over $B$ and let $%
x_{1},...,x_{s},y_{1},...,y_{s}\in (A,\varphi )$ be $B$-valued random
variables ($s\in \Bbb{N}$). If $\{x_{1},...,x_{s}\}$ and $\{y_{1},...,y_{s}\}
$ are free in $(A,\varphi ),$ then we have

\strut 

$k_{n}(x_{i_{1}}y_{i_{1}},...,x_{i_{n}}y_{i_{n}})$

\strut 

$=\underset{\pi \in NC(n)}{\sum }\left( \widehat{c_{x}}\oplus \widehat{c_{y}}%
\right) (\pi \cup Kr(\pi ))(x_{i_{1}}\otimes y_{i_{1}}\otimes
b_{i_{2}}x_{i_{2}}\otimes y_{i_{2}}\otimes ...\otimes
b_{i_{n}}x_{i_{n}}\otimes y_{i_{n}})$

\strut 

$\overset{denote}{=}\underset{\pi \in NC(n)}{\sum }\left( k_{\pi }\oplus
k_{Kr(\pi )}^{symm(1_{B})}\right)
(x_{i_{1}},y_{i_{1}},...,x_{i_{n}},y_{i_{n}}),$

\strut 

for all $(i_{1},...,i_{n})\in \{1,...,s\}^{n},$ $\forall n\in \Bbb{N},$ $%
b_{i_{2}},...,b_{i_{n}}\in B,$ arbitrary, where $\widehat{c_{x}}\oplus 
\widehat{c_{y}}=\widehat{c}\mid _{A_{x}*_{B}A_{y}},$ $A_{x}=A\lg \left(
\{x_{i}\}_{i=1}^{s},B\right) $ and $A_{y}=A\lg \left(
\{y_{i}\}_{i=1}^{s},B\right) .$ \ $\square $
\end{proposition}

\strut

This shows that ;

\strut

\begin{corollary}
(See [9]) Under the same condition with the previous proposition,

\strut 

\begin{center}
$R_{x_{1},...,x_{s}}\,\,\frame{*}\,%
\,R_{y_{1},...,y_{s}}^{t}=R_{x_{1}y_{1},...,x_{s}y_{s}}.$
\end{center}

$\square $
\end{corollary}

\strut

Notice that, in general, unless $b_{i_{2}}^{\prime }=...=b_{i_{n}}^{\prime
}=1_{B}$ in $B,$

\strut

\begin{center}
$R_{x_{1},...,x_{s}}\,\,\frame{*}\,\,R_{y_{1},...,y_{s}}\neq
R_{x_{1}y_{1},...,x_{s}y_{s}}.$
\end{center}

\strut

However, as we can see above,

\strut

\begin{center}
$R_{x_{1},...,x_{s}}\,\,\frame{*}\,%
\,R_{y_{1},...,y_{s}}^{t}=R_{x_{1}y_{1},...,x_{s}y_{s}}$
\end{center}

and

\begin{center}
$R_{x_{1},...,x_{s}}^{t}\,\,\frame{*}\,%
\,R_{y_{1},...,y_{s}}^{t}=R_{x_{1}y_{1},...,x_{s}y_{s}}^{t},$
\end{center}

\strut

where $\{x_{1},...,x_{s}\}$ and $\{y_{1},...,y_{s}\}$ are free over $B.$
Over $B=\Bbb{C},$ the last equation is proved by Nica and Speicher in [6]
and [7]. Actually, their R-transforms (over $\Bbb{C}$) is our trivial
R-transforms (over $\Bbb{C}$).

\strut

\strut

\strut

\strut

\section{\strut Amalgamated R-transform Theory on a Tower of Algebras}

\strut

\strut

\strut

\strut

\subsection{R-transform Calculus On a Tower of Amalgamated NCPSpaces}

\strut

\strut

\strut Let $B$ be a unital algebra and let $A_{k}$'s be algebras over $B,$
for all $k=1,2,\cdot \cdot \cdot .$ Moreover, assume that $A_{j+1}$ is an
algebra over $A_{j},$ for each $\ j=1,2,.....$ Then we can get a tower of
algebras,

\strut

\begin{center}
$\Bbb{C}\subset B\subset A_{1}\subset A_{2}\subset \cdot \cdot \cdot .$
\end{center}

\strut

Assume that there exists a sequence of conditional expectations $\left(
\varphi _{n}:A_{n}\rightarrow A_{n-1}\right) _{n=1}^{\infty },$ where $%
A_{0}=B.$ We can express this by

\strut

\begin{center}
$\Bbb{C}\subset ^{\varphi _{0}}B\subset ^{\varphi _{1}}A_{1}\subset
^{\varphi _{2}}A_{2}\subset ^{\varphi _{3}}\cdot \cdot \cdot ,$
\end{center}

\strut \strut

where $\varphi _{0}:B\rightarrow \Bbb{C}$ is a fixed linear functional.

\strut \strut \strut

\begin{definition}
Suppose that we have a tower of algebras,

\strut 

\begin{center}
$\Bbb{C}\subset B\subset A_{1}\subset A_{2}\subset A_{3}\subset \cdot \cdot
\cdot $
\end{center}

\strut 

such that $A_{j}$ is an algebra over $A_{j-1},$ with $A_{0}=B,$ for each $%
j\in \Bbb{N}.$ And assume that there is a sequence of conditional
expectations, $\varphi _{n}:A_{n}\rightarrow A_{n-1},$ for all $n\in \Bbb{N}.
$ Then the tower of algebras together with conditional expectations,

\strut 

\begin{center}
$\Bbb{C}\subset ^{\varphi _{0}}B\subset ^{\varphi _{1}}A_{1}\subset
^{\varphi _{2}}A_{2}\subset ^{\varphi _{3}}A_{3}\subset ^{\varphi _{4}}\cdot
\cdot \cdot $
\end{center}

\strut 

is called a tower of amalgamated noncommutative probability spaces, where $%
\varphi _{0}:B\rightarrow \Bbb{C}$ is the fixed linear functional.
\end{definition}

\strut \strut \strut

If we fix $\ j\in \Bbb{N},$ then the inclusion $[A_{j}\subset ^{\varphi
_{j+1}}A_{j+1}]$ satisfies the R-transform theory observed in Chapter 1. i.e
for each $j,$ we can regard $[A_{j}\subset ^{\varphi _{j+1}}A_{j+1}]$ as $%
[B\subset ^{\varphi }A]$ observed in Chapter 1.

\strut \strut

We are interested in the case $[\Bbb{C}\subset ^{E_{N}}A_{N}],$ where $N\in 
\Bbb{N}$ and

\strut

\begin{center}
$E_{N}:A_{N}\rightarrow \Bbb{C}$
\end{center}

\strut

is the linear functional defined by

\strut

\begin{center}
$E_{N}=\varphi _{0}\circ \varphi _{1}\circ \varphi _{2}\circ \cdot \cdot
\cdot \circ \varphi _{N}\overset{denote}{=}\varphi _{0}\varphi
_{1}...\varphi _{N}.$
\end{center}

\strut

Indeed, $E_{N}:A_{N}\rightarrow \Bbb{C}$ is a linear functional ;

\strut

$\ \ \ E_{N}\left( \alpha x+\beta y\right) =\varphi _{0}\varphi
_{1}...\varphi _{N}\left( \alpha x+\beta y\right) $

$\ \ \ \ \ \ \ \ \ \ \ \ \ \ \ \ \ \ \ \ \ \ \ \ \ \ \ =\varphi _{0}\varphi
_{1}...\varphi _{N-1}\left( \alpha \varphi _{N}\left( x\right) +\beta
\varphi _{N}(y)\right) $

\strut

by the bimodule map property of $\varphi _{N}$

\strut

$\ \ \ \ \ \ \ \ \ \ \ \ \ \ \ \ \ \ \ \ \ \ \ \ \ \ \ =\varphi _{0}\varphi
_{1}...\varphi _{N-2}\left( \alpha \varphi _{N-1}\varphi _{N}(x)+\beta
\varphi _{N-1}\varphi _{N}(y)\right) $

\strut

by the bimodule map property of $\varphi _{N-1}$

$\ \ \ \ \ \ \ \ \ \ \ \ \ \ \ \ \ \ \ \ \ \ \ \ \ \ \ =...=\alpha \varphi
_{0}\varphi _{1}...\varphi _{N}(x)+\beta \varphi _{0}\varphi _{1}...\varphi
_{N}(y)$

\strut

by the linearity of $\varphi _{0}.$

\strut

So, $E_{N}=\varphi _{0}\varphi _{1}...\varphi _{N}$ is a well-determined
linear functional on $A_{N},$ for each $N.$

\strut

\begin{definition}
Let $B$ be a unital algebra$.$ Let $\left( \varphi _{k}:A_{k}\rightarrow
A_{k-1}\right) _{k=1}^{\infty },$ with $A_{0}=B$ be a sequence of
conditional expectations, where $\varphi _{0}:B\rightarrow \Bbb{C}$ is a
fixed linear functional. i.e, we have a tower of NCPSpaces,

\strut 

\begin{center}
$\Bbb{C}\subset ^{\varphi _{0}}B\subset ^{\varphi _{1}}A_{1}\subset
^{\varphi _{2}}A_{2}\subset ^{\varphi _{3}}A_{3}\subset ^{\varphi _{4}}\cdot
\cdot \cdot .$
\end{center}

\strut 

For the fixed $j\in \Bbb{N},$ if $x_{1},...,x_{s}\in (A_{j+1},\varphi _{j+1})
$ are $A_{j}$-valued random variables ($s\in \Bbb{N}$), then the moment
series of them is defined by

\strut 

\begin{center}
$M_{x_{1},...,x_{s}}^{(j+1)}(z_{1},...,z_{s})=\sum_{n=1}^{\infty }\underset{%
i_{1},...,i_{n}\in \{1,...,s\}}{\sum }\varphi _{j+1}\left(
x_{i_{1}}a_{i_{2}}x_{i_{2}}...a_{i_{n}}x_{i_{n}}\right)
\,z_{i_{1}}...z_{i_{n}},$
\end{center}

\strut 

as a $A_{j}$-formal series in $\Theta _{A_{j}}^{s},$ where $%
b_{i_{2}},...,b_{i_{b}}\in A_{j}$ are arbitrary for each $%
(i_{1},...,i_{n})\in \{1,...,s\}^{n},$ $n\in \Bbb{N}.$ Similarly, in $\Theta
_{A_{j}}^{s},$ we can define a R-transform of $x_{1},...,x_{s}$ by

\strut 

\begin{center}
$R_{x_{1},...,x_{s}}^{(j+1)}(z_{1},...,z_{s})=\sum_{n=1}^{\infty }\underset{%
i_{1},...,i_{n}\in \{1,...,s\}}{\sum }k_{n}^{(j+1)}\left(
x_{i_{1}},...,x_{i_{n}}\right) \,z_{i_{1}}...z_{i_{n}},$
\end{center}

\strut 

with its $(i_{1},...,i_{n})$-th coefficients

\strut 

\begin{center}
$k_{n}^{(j+1)}\left( x_{i_{1}},...,x_{i_{n}}\right) =c_{j+1}^{(n)}\left(
x_{i_{1}}\otimes a_{i_{2}}x_{i_{2}}\otimes ...\otimes
a_{i_{n}}x_{i_{n}}\right) ,$
\end{center}

\strut 

where $a_{i_{2}},...,a_{i_{n}}\in A_{j}$ are arbitrary for all $%
(i_{1},...,i_{n})\in \{1,...,s\}^{n},$ $n\in \Bbb{N},$ and where $\widehat{%
c_{j+1}}=\left( c_{j+1}^{(k)}\right) _{k=1}^{\infty }\in I^{c}\left(
A_{j+1},A_{j}\right) $ is the cumulant multiplicative bimodule map induced
by $\varphi _{j+1}:A_{j+1}\rightarrow A_{j}.$ Similar to the $[B\subset
^{\varphi }A]$-case in Chapter 1, we can define trivial moment series and
trivial R-transform of $x_{1},...,x_{s}$ by

\strut 

\begin{center}
$M_{x_{1},...,x_{s}}^{(j+1)\,\,:\,\,t}(z_{1},...,z_{s})=\sum_{n=1}^{\infty }%
\underset{i_{1},...,i_{n}\in \{1,...,s\}}{\sum }\varphi _{j+1}\left(
x_{i_{1}}x_{i_{2}}...x_{i_{n}}\right) \,z_{i_{1}}...z_{i_{n}}$
\end{center}

and

\begin{center}
$R_{x_{1},...,x_{s}}^{(j+1)\,\,:\,\,t}(z_{1},...,z_{s})=\sum_{n=1}^{\infty }%
\underset{i_{1},...,i_{n}\in \{1,...,s\}}{\sum }k_{n}^{(j+1)\,:\,\,t}\left(
x_{i_{1}},...,x_{i_{n}}\right) \,z_{i_{1}}...z_{i_{n}},$
\end{center}

\strut 

in $\Theta _{A_{j}}^{s},$ as $A_{j}$-formal series, respectively. Let $%
E_{j+1}=\varphi _{0}\varphi _{1}\cdot \cdot \cdot \varphi
_{j+1}:A_{j+1}\rightarrow \Bbb{C}$ be a linear functional. Then a moment
series of $x_{1},...,x_{s}$ and an R-transform of $x_{1},...,x_{s}$ are
denoted by

\strut 

\begin{center}
$m_{x_{1},...,x_{s}}(z_{1},...,z_{s})=\sum_{n=1}^{\infty }\underset{%
i_{1},...,i_{n}\in \{1,...,s\}}{\sum }E_{j+1}\left(
x_{i_{1}}x_{i_{2}}...x_{i_{n}}\right) \,z_{i_{1}}...z_{i_{n}}$
\end{center}

\strut and

\begin{center}
$r_{x_{1},...,x_{s}}(z_{1},...,z_{s})=\sum_{n=1}^{\infty }\underset{%
i_{1},...,i_{n}\in \{1,...,s\}}{\sum }k_{n}^{(E_{j+1})}\left(
x_{i_{1}},...,x_{i_{n}}\right) z_{i_{1}}...z_{i_{n}},$
\end{center}

\strut 

as formal series in $\Theta _{\Bbb{C}}^{s}\overset{denote}{=}\Theta _{s},$
respectively. Here, the $(i_{1},...,i_{n})$-th $E_{j+1}$-moments and $%
(i_{1},...,i_{n})$-th $E_{j+1}$-cumulants are defined in the sense of
Speicher and Nica (See [6] and [7]). i.e, by using our notation introduced
in Chapter 1, they are just scalar-valued trivial moments and trivial
cumulants.
\end{definition}

\strut

Since we can regard $E_{j+1}$ as a linear functional and $k_{n}^{(E_{j+1})}$
as a scalar-valued cumulants, we will use the notation used in [6] and [7],
for the convenience. (For example, instead of using the notation for $%
\widehat{E_{j+1}}(\pi )(...)$ and $\widehat{c_{E_{j+1}}}(\pi )(...),$ we
will use $E_{j+1\,\,:\,\,\pi }(...)$ and $k_{\pi }^{(E_{j+1})}(...),$ like
in [6] and [7].)

\strut \strut

Notice that we can apply all our results in Chapter 1 to $[A_{j}\subset
^{\varphi _{j+1}}A_{j+1}],$ for the fixed $j\in \Bbb{N},$ in the tower of
amalgamated NCPSpaces,

\strut

\begin{center}
$\Bbb{C}\subset ^{\varphi _{0}}B\subset ^{\varphi _{1}}A_{1}\subset
^{\varphi _{2}}A_{2}\subset ^{\varphi _{3}}A_{3}\subset ^{\varphi _{4}}\cdot
\cdot \cdot .$
\end{center}

\strut

By \frame{*}$_{A_{j}},$ we will denote the $A_{j}$-valued boxed convolution
on $\Theta _{A_{j}}^{s},$ for all $s\in \Bbb{N}.$ Also, we can consider the $%
A_{j}$-freeness for each $j.$ In this section, we will concentrate on
observing $B$-freeness of $E_{N}$ and considering $E_{N}$-cumulants. But
first, we can get the following R-transform calculus for each step $%
[A_{j}\subset ^{\varphi _{j+1}}A_{j+1}].$ By definition and by Chapter 1, we
have the following $A_{j}$-valued R-transform calculus, for each $j$-th step
of a tower of NCPSpaces.

\strut \strut

\begin{proposition}
Let $B$ be a unital algebra and assume that we have a tower of algebras $%
T_{B}\left( (A_{i})_{i=1}^{\infty }\right) $ and a sequence of conditional
expectations $\left( \varphi _{j}:A_{j}\rightarrow A_{j-1}\right)
_{j=1}^{\infty },$ with $A_{0}=B.$ i.e, we have a tower of NCPSpaces

\strut 

\begin{center}
$\Bbb{C}\subset ^{\varphi _{0}}B\subset ^{\varphi _{1}}A_{1}\subset
^{\varphi _{2}}A_{2}\subset ^{\varphi _{3}}A_{3}\subset ^{\varphi _{4}}\cdot
\cdot \cdot .$
\end{center}

\strut 

Fix $j\in \Bbb{N}.$ Let $x_{1},...,x_{s},y_{1},...,y_{s}\in (A_{j+1},\varphi
_{j+1})$ be $A_{j}$-valued random variables ($s\in \Bbb{N}$). Assume that
two subsets of $A_{j+1},$ $X=\{x_{1},...,x_{s}\}$ and $Y=\{y_{1},...,y_{s}\}$
are free over $A_{j}$ (in short, $A_{j}$-free). Then

\strut 

(1) $%
R_{x_{1},...,x_{s},y_{1},...,y_{s}}^{(j+1)}(z_{1},...,z_{2s})=R_{x_{1},...,x_{s}}^{(j+1)}(z_{1},...,z_{s})+R_{y_{1},...,y_{s}}^{(j+1)}(z_{s+1},...,z_{2s}).
$

\strut 

(2) $R_{x_{1}+y_{1},...,x_{s}+y_{s}}^{(j+1)}(z_{1},...,z_{s})=\left(
R_{x_{1},...,x_{s}}^{(j+!)}+R_{y_{1},...,y_{s}}^{(j+1)}\right)
(z_{1},...,z_{s}).$

\strut 

(3) $R_{x_{1}y_{1},...,x_{s}y_{s}}^{(j+1)}(z_{1},...,z_{s})=\left(
R_{x_{1},...,x_{s}}^{(j+1)}\,\,\,\frame{*}_{A_{j}}\,%
\,R_{y_{1},...,y_{s}}^{(j+1)\,\,:\,\,t}\right) (z_{1},...,z_{s}).$ \ $%
\square $
\end{proposition}

\strut \strut

From now, we will observe the $E_{j+1}$-cumulants of random variables.
Notice that $\left( A_{j+1},\,E_{j+1}\right) $ is a NCPSpace (over $\Bbb{C}$%
), since $E_{j+1}:A_{j+1}\rightarrow \Bbb{C}$ is a well-determined linear
functional. Hence if $x_{1},...,x_{s}\in A_{j+1}$ are operators, then they
can be regarded as (scalar-valued) random variables ($s\in \Bbb{N}$). So, $%
E_{j+1}$-cumulants of them are well-defined with respect to the
multiplicative bimodule map $\widehat{c_{E_{j+1}}},$ induced by $E_{j+1}.$

\strut \strut \strut \strut

\begin{lemma}
Let $B$ be a unital algebra and let

\strut 

\begin{center}
$\Bbb{C}\subset ^{\varphi _{0}}B\subset ^{\varphi _{1}}A_{1}\subset
^{\varphi _{2}}A_{2}\subset ^{\varphi _{3}}A_{3}\subset ^{\varphi _{4}}\cdot
\cdot \cdot $
\end{center}

\strut 

be a tower of amalgamated NCPSpaces. Let $E_{j+1}=\varphi _{0}\varphi
_{1}....\varphi _{j+1}:A_{j+1}\rightarrow B$ be a $B$-functional, for each
\thinspace $j=1,...,n.$ Let $x_{1},...,x_{s}\in (A_{j+1},E_{j+1})$ be $B$%
-valued random variables ($s\in \Bbb{N}$). Then

\strut 

$k_{n}^{(E_{j+1})}\left( x_{i_{1}},...,x_{i_{n}}\right) $

\strut 

\begin{center}
$=\underset{\pi \in NC(n)}{\sum }\left( \underset{V=(v_{1},...,v_{k})\in \pi 
}{\prod }\left( \underset{\theta \in NC(k)}{\sum }E_{j}(k_{k\,\,:\,\,\theta
}^{(j+1)\,\,:\,\,t}(x_{v_{1}},...,x_{v_{k}}))\right) \right) \,\mu (\pi
,1_{n}),$
\end{center}

\strut 

for all $(i_{1},...,i_{n})\in \{1,...,s\}^{n},$ $n\in \Bbb{N}.$
\end{lemma}

\strut

\begin{proof}
Fix $n\in \Bbb{N}$ and $(i_{1},...,i_{n})\in \{1,...,s\}^{n}.$ Then

\strut

$k_{n}^{(E_{j+1})}\left( x_{i_{1}},...,x_{i_{n}}\right)
=c_{E_{j+1}}^{(n)}\left( x_{i_{1}}\otimes x_{i_{2}}\otimes ...\otimes
x_{i_{n}}\right) $

\strut

by using of our notation (See Chapter 1)

\strut

$\ \ \ =\underset{\pi \in NC(n)}{\sum }\widehat{E_{j+1}}(\pi )\left(
x_{i_{1}}\otimes b_{i_{2}}x_{i_{2}}\otimes ...\otimes
b_{i_{n}}x_{i_{n}}\right) \,\mu (\pi ,1_{n})$

\strut

$\ \ \ =\underset{\pi \in NC(n)}{\sum }E_{j+1\,\,:\,\,\pi }\left(
x_{i_{1}},...,x_{i_{n}}\right) \,\mu (\pi ,1_{n})$

\strut

by using notation of Nica and Speicher (See [6] and [7])

\strut

$\ \ \ =\underset{\pi \in NC(n)}{\sum }\left( \underset{V=(v_{1},...,v_{k})%
\in \pi }{\prod }E_{j+1}\left( x_{v_{1}}...x_{v_{k}}\right) \right) \,\mu
(\pi ,1_{n}).$

\strut

(Since $E_{j+1}$ is a linear functional, we don't need to consider the
insertion property with respect to inner partitions and outer partitions.
So, we can conclude the above last equality.)

\strut \strut

Now, observe $E_{j+1}(x_{v_{1}}...x_{v_{k}}),$ for the fixed block $%
V=(v_{1},...,v_{k})\in \pi ,$ for $\pi \in NC(n)$ ;

\strut

$E_{j+1}\left( x_{v_{1}}...x_{v_{n}}\right) =\varphi _{0}\varphi
_{1}...\varphi _{j}\varphi _{j+1}(x_{v_{1}}...x_{v_{k}})$

$\ \ \ \ \ \ \ \ =\varphi _{0}\varphi _{1}...\varphi _{j}\left( \varphi
_{j+1}(x_{v_{1}}...x_{v_{k}})\right) $

\strut

(2.2.1)

\strut

$\ \ \ \ \ \ \ \ =\varphi _{0}\varphi _{1}...\varphi _{j}\left( \underset{%
\theta \in NC(k)}{\sum }\widehat{c_{j+1}}(\theta )\left( x_{v_{1}}\otimes
...\otimes x_{v_{k}}\right) \right) .$

$\strut $

From now, we will denote a partition-dependent amalgamated \textbf{trivial}
cumulant having the form, $\widehat{c_{j+1}}(\theta )\left( x_{v_{1}}\otimes
...\otimes x_{v_{k}}\right) ,$ by $k_{k\,\,:\,\,\theta
}^{(j+1)\,:\,\,t}(x_{v_{1}},...,x_{v_{k}}).$ By using this new notation, we
can re-express the formula (2.2.1) ;

\strut

(2.2.2)

\strut

$\ \ \ \ \ \ \ \ =\varphi _{0}\varphi _{1}...\varphi _{j}\left( \underset{%
\theta \in NC(k)}{\sum }k_{k\,\,:\,\,\,\theta
}^{(j+1)\,\,:\,\,t}(x_{v_{1}},...,x_{v_{k}})\right) $

\strut

$\ \ \ \ \ \ \ \ =\varphi _{0}\varphi _{1}...\varphi _{j-1}\left( \underset{%
\theta \in NC(k)}{\sum }\varphi _{j}\left( k_{k\,\,:\,\,\theta
}^{(j+1)\,\,:\,\,t}(x_{v_{1}},...,x_{v_{k}})\right) \right) $

\strut

$\ \ \ \ \ \ \ \ =.......$

\strut

$\ \ \ \ \ \ \ \ =\underset{\theta \in NC(k)}{\sum }\varphi _{0}\varphi
_{1}...\varphi _{j}\left( k_{k\,\,:\,\,\theta
}^{(j+1)\,\,:\,\,t}(x_{v_{1}},...,x_{v_{k}})\right) $

\strut

$\ \ \ \ \ \ \ \ =\underset{\theta \in NC(k)}{\sum }E_{j}\left(
k_{k\,\,:\,\,\theta }^{(j+1)\,:\,t}(x_{v_{1}},...,x_{v_{k}})\right) .$

\strut

Therefore,

\strut

$k_{n}^{(E_{j+1})}\left( x_{i_{1}},...,x_{i_{n}}\right) =\underset{\pi \in
NC(n)}{\sum }\left( \underset{V=(v_{1},...,v_{k})\in \pi }{\prod }%
E_{j+1}\left( x_{v_{1}}...x_{v_{k}}\right) \right) \,\mu (\pi ,1_{n})$

\strut

$\ \ \ \ \ \ =\underset{\pi \in NC(n)}{\sum }\left( \underset{%
V=(v_{1},...,v_{k})\in \pi }{\prod }\left( \underset{\theta \in NC(k)}{\sum }%
E_{j}(k_{k\,\,:\,\,\theta
}^{(j+1)\,\,:\,\,t}(x_{v_{1}},...,x_{v_{k}}))\right) \right) \,\mu (\pi
,1_{n}).$
\end{proof}

\strut \strut

The above lemma shows the relation between $E_{j+1}$-cumulants and $E_{j}$%
-moments. Also, this shows that to compute the $E_{j+1}$-cumulants, we only
need to consider the M\"{o}bius inversion for $A_{j}$-valued moments with
respect to $\varphi _{j+1},$ for each block.

\strut

\begin{corollary}
Let $B$ be a unital algebra and let

\strut 

\begin{center}
$\Bbb{C}\subset ^{\varphi _{0}}B\subset ^{\varphi _{1}}A_{1}\subset
^{\varphi _{2}}A_{2}\subset ^{\varphi _{3}}A_{3}\subset ^{\varphi _{4}}\cdot
\cdot \cdot $
\end{center}

\strut 

be a tower of amalgamated NCPSpaces. Define $E_{j+1}=\varphi _{0}\varphi
_{1}...\varphi _{j}\varphi _{j+1},$ for each $j=0,1,....$ Let $%
x_{1},...,x_{s}\in (A_{j+1},E_{j+1})$ be random variables ($s\in \Bbb{N}$).
Then

\strut 

\begin{center}
$k_{n}^{(E_{j+1})}\left( x_{i_{1}},...,x_{i_{n}}\right) =\underset{\pi \in
NC(n)}{\sum }\left( \underset{V=(v_{1},...,v_{k})\in \pi }{\prod }%
E_{j}\left( \varphi _{j+1}(x_{v_{1}}...x_{v_{k}})\right) \right) \,\mu (\pi
,1_{n}),$
\end{center}

\strut 

for all $(i_{1},...,i_{n})\in \{1,...,s\}^{n},$ $n\in \Bbb{N}.$ \ $\square $
\end{corollary}

\strut \strut

\begin{example}
\strut Let $j=2$ and let $x,y\in (A_{2},E_{2})$ be $B$-valued random
variables. Then, by the straightforward computation, we have that

\strut 

$k_{3}^{(E_{2})}\left( x,x,y\right) =c^{(3)}\left( x\otimes x\otimes
y\right) $

\strut 

$\ \ \ \ =\underset{\pi \in NC(3)}{\sum }E_{2\,\,:\,\,\pi }\left(
x,x,y\right) \mu (\pi ,1_{3})$

\strut 

$\ \ \ \ =E_{2}(xxy)-E_{2}(xx)E_{2}(y)-E_{2}\left( xy\right) E_{2}(x)$

$\ \ \ \ \ \ \ \ \ \ \ \ \ \ \ \ \ \ \ \ \ \ \ \ \ \ \ \ \
-E_{2}(x)E_{2}(xy)+E_{2}(x)E_{2}(x)E_{2}(y)$

\strut 

$\ \ \ \ =\varphi _{0}\varphi _{1}\varphi _{2}\left( xxy\right) -\varphi
_{0}\varphi _{1}\varphi _{2}(xx)\varphi _{0}\varphi _{1}\varphi
_{2}(y)-\varphi _{0}\varphi _{1}\varphi _{2}\left( xy\right) \varphi
_{0}\varphi _{1}\varphi _{2}(x)$

$\ \ \ \ \ \ \ \ \ \ \ \ \ \ \ \ \ \ -\varphi _{0}\varphi _{1}\varphi
_{2}(x)\varphi _{0}\varphi _{1}\varphi _{2}(xy)+\varphi _{0}\varphi
_{1}\varphi _{2}(x)\cdot \varphi _{0}\varphi _{1}\varphi _{2}(x)\cdot
\varphi _{0}\varphi _{1}\varphi _{2}(y)$

\strut 

$\ \ \ =E_{1}\left( \varphi _{2}(x^{2}y)\right) -E_{1}\left( \varphi
_{2}(x^{2})\right) E_{1}\left( \varphi _{2}(y)\right) -E_{1}\left( \varphi
_{2}(xy)\right) E_{1}\left( \varphi _{2}(x)\right) $

$\ \ \ \ \ \ \ \ \ \ \ \ \ \ \ \ -E_{1}\left( \varphi _{2}(x)\right)
E_{1}\left( \varphi _{2}(xy)\right) +E_{1}\left( \varphi _{2}(x)\right)
E_{1}\left( \varphi _{2}(x)\right) E_{1}\left( \varphi _{2}(y)\right) .$

\strut 

So, the above lemma and corollary are applied well in the above example.
Here, by the above corollary, we have the informal way how to find such $%
E_{2}$-cumulant ;

\strut 

(i) Consider the trivial $A_{1}$-valued cumulant ;

\strut 

$\ \ \ \ \ k_{3}^{(2)\,\,:\,\,t}\left( x,x,y\right) =\underset{\pi \in NC(3)%
}{\sum }\widehat{\varphi _{2}}(\pi )\left( x,x,y\right) \mu (\pi ,1_{3})$

\strut 

$\ \ \ \ \ \ \ \ \ \ \ \ \ =\varphi _{2}\left( x^{2}y\right) -\varphi
_{2}(x)\varphi _{2}\left( xy\right) -\varphi _{2}\left( x\varphi
_{2}(x)y\right) $

$\ \ \ \ \ \ \ \ \ \ \ \ \ \ \ \ \ \ \ \ \ \ \ \ \ \ \ \ \ \ \ \ \ \ \ \ \ \
\ -\varphi _{2}\left( x^{2}\right) \varphi _{2}(y)+\varphi _{2}(x)\varphi
_{2}(x)\varphi _{2}(y)$.

\strut \strut 

(ii) Ignore the insertion property and act $E_{1}=\varphi _{0}\varphi _{1}.$
Then we can get that

\strut 

$\ \ \ \ \ \ k_{3}^{(E_{2})}(x,x,y)=E_{1}\left( \varphi _{2}\left(
x^{2}y\right) \right) -E_{1}\left( \varphi _{2}(x)\right) E_{1}\left(
\varphi _{2}\left( xy\right) \right) $

$\ \ \ \ \ \ \ \ \ \ \ \ \ \ \ \ \ \ \ \ \ \ \ \ \ \ \ \ -E_{1}\left(
\varphi _{2}\left( xy\right) \right) E_{1}\left( \varphi _{2}(x)\right)
-E_{1}\left( \varphi _{2}(x^{2})\right) E_{1}\left( \varphi _{2}(y)\right) $

$\ \ \ \ \ \ \ \ \ \ \ \ \ \ \ \ \ \ \ \ \ \ \ \ \ \ \ \ +E_{1}\left(
\varphi _{2}(x)\right) E_{1}\left( \varphi _{2}(x)\right) E_{1}\left(
\varphi _{2}(y)\right) $.

\strut $\strut $
\end{example}

\strut

\begin{corollary}
Let $B$ be a unital algebra and let

\strut 

\begin{center}
$\Bbb{C}\subset ^{\varphi _{0}}B\subset ^{\varphi _{1}}A_{1}\subset
^{\varphi _{2}}A_{2}\subset ^{\varphi _{3}}A_{3}\subset ^{\varphi _{4}}\cdot
\cdot \cdot $
\end{center}

\strut 

be a tower of amalgamated NCPSpaces. Let $E_{j+1}=\varphi _{0}\varphi
_{1}...\varphi _{j+1}$ be a linear functional on $A_{j+1},$ for all $j.$ Let 
$x\in (A_{j+1},\,E_{j+1})$ be a random varialbe. Then

\strut 

\begin{center}
$\ k_{n}^{(E_{j+1})}\left( \underset{n-times}{\underbrace{x,........,x}}%
\right) =\underset{\pi \in NC(n)}{\sum }\left( \underset{V=(v_{1},...,v_{k})%
\in \pi }{\prod }\,E_{j}\left( \varphi _{j+1}(x^{k})\right) \right) \,\mu
(\pi ,1_{n}).$
\end{center}

$\square $
\end{corollary}

\strut \strut

Now, we will observe the R-transform calculus for $E_{j+1}.$ The results are
came from Chapter 1, [6] and [7] ;

\strut \strut

\begin{proposition}
(See [6] and [7]) Let $B$ be a unital algebra and let

\strut 

\begin{center}
$\Bbb{C}\subset ^{\varphi _{0}}B\subset ^{\varphi _{1}}A_{1}\subset
^{\varphi _{2}}A_{2}\subset ^{\varphi _{3}}A_{3}\subset ^{\varphi _{4}}\cdot
\cdot \cdot $
\end{center}

\strut 

be a tower of amalgamated NCPSpaces. Then, for any fixed $j\in \Bbb{N},$ the
following R-transform calculus holds true ; suppose that $%
x_{1},...,x_{s},y_{1},...,y_{s}\in (A_{j+1},E_{j+1})$ are random variables ($%
s\in \Bbb{N}$) and assume that $X=\{x_{1},...,x_{s}\}$ and $%
Y=\{y_{1},...,y_{s}\}$ are free. Then

\strut 

(1) $
r_{x_{1},...,x_{s},y_{1},...,y_{s}}(z_{1},...,z_{2s})=r_{x_{1},...,x_{s}}(z_{1},...,z_{s})+r_{y_{1},...,y_{s}}(z_{s+1},...,z_{2s}).
$

\strut 

(2) $r_{x_{1}+y_{1},...,x_{s}+y_{s}}(z_{1},...,z_{s})=\left(
r_{x_{1},...,x_{s}}+r_{y_{1},...,y_{s}}\right) (z_{1},...,z_{s}).$ $\square $
\end{proposition}

\strut \strut

\begin{proposition}
(See [6] and [7]) Let $B$ be a unital algebra and let

\strut 

\begin{center}
$C\subset ^{\varphi _{0}}B\subset ^{\varphi _{1}}A_{1}\subset ^{\varphi
_{2}}A_{2}\subset ^{\varphi _{3}}A_{3}\subset ^{\varphi _{4}}\cdot \cdot
\cdot $
\end{center}

\strut 

be a tower of amalgamated NCPSpaces. Fix $j\in \Bbb{N}.$ Let $%
(A_{j+1},E_{j+1})$ be a NCPSpace over $B,$ with a linear functional, $%
E_{j+1}=\varphi _{0}\varphi _{1}...\varphi _{j+1}:A_{j+1}\rightarrow \Bbb{C}$
and let $x_{1},...,x_{s},y_{1},...,y_{s}\in (A_{j+1},E_{j+1})$ be $B$-valued
random variables ($s\in \Bbb{N}$). If $X=\{x_{1},...,x_{s}\}$ and $%
Y=\{y_{1},...,y_{s}\}$ are free over $B,$ then

\strut 

\begin{center}
$k_{n}^{(E_{j+1})}\left( x_{i_{1}}y_{i_{1}},...,x_{i_{n}}y_{i_{n}}\right) =%
\underset{\pi \in NC(n)}{\sum }\left( k_{\pi
}^{(E_{j+1})}(x_{i_{1}},...,x_{i_{n}})\right) \left( k_{Kr(\pi
)}^{(E_{j+1})}(y_{i_{1}},...,y_{i_{n}})\right) ,$
\end{center}

\strut 

for all $(i_{1},...,i_{n})\in \{1,...,s\}^{n},$ $n\in \Bbb{N}.$ $\square $
\end{proposition}

\strut

In the above proposition, $k_{\pi }^{(E_{j+1})}(...)$ and $k_{Kr(\pi
)}^{(E_{j+1})}(...)$ are used as partition-dependent $E_{j+1}$-cumulants, in
the sense of Speicher and Nica (See [6] and [7]). The above result can be
proved by [6] and [7]. Also, it can be proved by using the result in Chapter
1.

\strut

\begin{corollary}
(See [6] and [7]) Let $B$ be a unital algebra and let

\strut 

\begin{center}
$\Bbb{C}\subset ^{\varphi _{0}}B\subset ^{\varphi _{1}}A_{1}\subset
^{\varphi _{2}}A_{2}\subset ^{\varphi _{3}}A_{3}\subset ^{\varphi _{4}}\cdot
\cdot \cdot $
\end{center}

\strut 

be a tower of amalgamated NCPSpaces. Fix $j\in \Bbb{N}.$ Let $%
X=\{x_{1},...,x_{s}\}$ and $Y=\{y_{1},...,y_{s}\}$ be free subsets of random
variables ($s\in \Bbb{N}$) in a NCPSpace, $\left( A_{j+1},\,E_{j+1}\right) .$
Then

\strut 

\begin{center}
$r_{x_{1}y_{1},...,x_{s}y_{s}}(z_{1},...,z_{s})=\left(
r_{x_{1},...,x_{s}}\,\,\,\frame{*}_{s}\,\,\,r_{y_{1},...,y_{s}}\right)
(z_{1},...,z_{s}),$
\end{center}

\strut 

where \ \frame{*}$_{s}\,:\,\Theta _{s}\times \Theta _{s}\rightarrow \Theta
_{s}$ is a (scalar-valued) boxed convolution in the sense of Nica and
Speicher. $\square $
\end{corollary}

\strut

Define a conditional expectation

\strut

\begin{center}
$E_{k,\,j}:A_{j}\rightarrow A_{k-1},$
\end{center}

\strut by

\begin{center}
$E_{k,\,j}=\varphi _{k}\varphi _{k+1}...\varphi _{j}\varphi _{j},$
\end{center}

\strut

for all $k<j$ in $\Bbb{N}.$ Indeed, it is a conditional expectation, for the
fixed $k<j$ ;

\strut

(i) \ Let $a\in A_{k-1}.$ Then

\strut

$\ \ \ \ \ \ \ E_{k,\,j}(a)=\varphi _{k}\varphi _{k+1}...\varphi
_{j-1}\left( \varphi _{j}(a)\right) =\varphi _{k}...\varphi _{j-1}(a)$

\strut

\ \ \ \ since $a\in A_{k-1}\subset A_{j-1}$

\strut

$\ \ \ \ \ \ \ \ \ \ \ \ \ \ \ \ \ \ =...=\varphi _{k}(a)=a.$

\strut

\ \ \ \ \ So,

\begin{center}
$E_{k,\,j}(a)=a,$ for all $a\in A_{k-1}.$
\end{center}

\strut

(ii) Let $a,a^{\prime }\in A_{k-1}$ and $x\in A_{j}.$ Then

\strut

$\ \ \ \ \ \ E_{k,\,j}(axa^{\prime })=\varphi _{k}\varphi _{k+1}...\varphi
_{j-1}\left( \varphi _{j}(axa^{\prime })\right) $

$\ \ \ \ \ \ \ \ \ \ \ \ \ \ \ \ \ \ \ \ \ \ =\varphi _{k}\varphi
_{k+1}...\varphi _{j-1}\left( a\varphi _{j}(x)a^{\prime }\right) $

$\ \ \ \ \ \ \ \ \ \ \ \ \ \ \ \ \ \ \ \ \ \ =\varphi _{k}...\varphi
_{j-2}\left( \varphi _{j-1}(a\varphi _{j}(x)a^{\prime })\right) $

$\ \ \ \ \ \ \ \ \ \ \ \ \ \ \ \ \ \ \ \ \ \ =\varphi _{k}...\varphi
_{j-2}\left( a\cdot \varphi _{j-1}\varphi _{j}(x)\cdot a^{\prime }\right) $

$\ \ \ \ \ \ \ \ \ \ \ \ \ \ \ \ \ \ \ \ \ \ =.....=a\cdot \varphi
_{k}\varphi _{k+1}...\varphi _{j}(x)\cdot a^{\prime }=aE_{k,\,j}(x)a^{\prime
}.$

\strut

\ \ \ \ \ Thus we have that

\begin{center}
$E_{k,\,\,j}\left( axa^{\prime }\right) =aE_{k,\,\,j}(x)a^{\prime },$ for
all $a,a^{\prime }\in A_{k-1}$ and $x\in A_{j}.$
\end{center}

\strut

By (i) and (ii), $E_{k,\,\,j}:A_{j}\rightarrow A_{k-1}$ is a conditional
expectation (or a $A_{k-1}$-functional on $A_{j}$). And hence $\left(
A_{j},\,\,E_{k,\,j}\right) $ is a NCPSpace over $A_{k-1}.$ We can observe
the following R-transform calculus, by Chapter 1 ; Defien, for $A_{k-1}$%
-valued random variables $x_{1},...,x_{s}$ in $\left(
A_{j},\,E_{k,\,j}\right) ,$

\strut

\begin{center}
$M_{x_{1},...,x_{s}}^{(E_{k,\,\,j})}(z_{1},...,z_{s})=\sum_{n=1}^{\infty }%
\underset{i_{1},...,i_{n}\in \{1,...,s\}}{\sum }E_{k,\,j}\left(
x_{i_{1}}a_{i_{2}}x_{i_{2}}...a_{i_{n}}x_{i_{n}}\right)
\,z_{i_{1}}...z_{i_{n}}$
\end{center}

and

\begin{center}
$R_{x_{1},...,x_{s}}^{(E_{k\,,\,j})}(z_{1},...,z_{s})=\sum_{n=1}^{\infty }%
\underset{i_{1},...,i_{n}\in \{1,...,s\}}{\sum }k_{n}^{(E_{k,\,j})}\left(
x_{i_{1}},...,x_{i_{n}}\right) \,z_{i_{1}}...z_{i_{n}},$
\end{center}

\strut

as the moment series of $x_{1},...,x_{s}$ and the R-transform of $%
x_{1},...,x_{s},$ respectively, in $\Theta _{A_{k-1}}^{s}.$ Here, the $%
(i_{1},...,i_{n})$-th cumulants of $x_{1},..,x_{s}$ are determined by

\strut

\begin{center}
$k_{n}^{(E_{k,\,j})}\left( x_{i_{1}},...,x_{i_{n}}\right)
=c_{E_{k,\,j}}^{(n)}\left( x_{i_{1}}\otimes a_{i_{2}}x_{i_{2}}\otimes
...\otimes a_{i_{n}}x_{i_{n}}\right) ,$
\end{center}

\strut

where $a_{i_{2}},...,a_{i_{n}}\in A_{k-1}$ are arbitrary and $\widehat{%
c_{E_{k,\,j}}}=\left( c_{E_{k,\,j}}^{(m)}\right) _{m=1}^{\infty }\in
I^{c}\left( A_{j},\,\,A_{k-1}\right) $ is the cumulant multiplicative
bimodule map induced by a $A_{k-1}$-functional, $E_{k,\,\,j}.$ A trivial
moment series of $x_{1},...,x_{s}$, denoted by $M_{x_{1},...,x_{s}}^{(E_{k,%
\,j})\,:\,t},$ and a trivial R-transform of $x_{1},...,x_{s},$ denoted by $%
R_{x_{1},...,x_{s}}^{(E_{k,\,j})\,:\,t},$ are defined usually.

\strut \strut

\begin{proposition}
Let $B$ be a unital algebra and let

\strut 

\begin{center}
$\Bbb{C}\subset ^{\varphi _{0}}B\subset ^{\varphi _{1}}A_{1}\subset
...A_{k-1}\subset ^{\varphi _{k}}A_{k}\subset ...\subset ^{\varphi
_{j}}A_{j}\subset ^{\varphi _{j+1}}\cdot \cdot \cdot $
\end{center}

\strut 

be a tower of amalgamated NCPSpaces. Fix $k<j$ in $\Bbb{N}.$ Define a
conditional expectation $E_{k,\,j}:A_{j}\rightarrow A_{k-1}$, where $A_{0}=B.
$ Let $X=\{x_{1},...,x_{s}\}$ and $Y=\{y_{1},...,y_{s}\}$ be two subsets of $%
A_{k-1}$-valued random variables ($s\in \Bbb{N}$) in $\left(
A_{j},\,\,E_{k,\,j}\right) .$ If $X$ and $Y$ are free over $A_{k-1},$ in $%
\left( A_{j},\,E_{k,\,j}\right) ,$ then we have the following R-transform
calculus ;

\strut 

(1) $R_{x_{1},...,x_{s},y_{1},...,y_{s}}^{(E_{k,%
\,j})}(z_{1},...,z_{2s})=R_{x_{1},...,x_{s}}^{(E_{k,%
\,j})}(z_{1},...,z_{s})+R_{y_{1},...,y_{s}}^{(E_{k,%
\,j})}(z_{s+1},...,z_{2s}).$

\strut 

(2) $R_{x_{1}+y_{1},...,x_{s}+y_{s}}^{(E_{k,\,\,j})}(z_{1},...,z_{s})=\left(
R_{x_{1},...,x_{s}}^{(E_{k,\,j})}+R_{y_{1},...,y_{s}}^{(E_{k,\,j})}\right)
(z_{1},...,z_{s}).$

\strut 

(3) $R_{x_{1}y_{1},...,x_{s}y_{s}}^{(E_{k,\,\,j})}(z_{1},...,z_{s})=\left(
R_{x_{1},...,x_{s}}^{(E_{k,\,\,j})}\,\,\frame{*}_{A_{k-1}}^{(E_{k,\,\,j})}\,%
\,\,\,R_{y_{1},...,y_{s}}^{(E_{k,\,\,j})}\right) (z_{1},...,z_{s}),$

\strut 

where \frame{*}$_{A_{k-1}}:\mathcal{\Theta }_{A_{k-1}}^{s}\times \mathcal{%
\Theta }_{A_{k-1}}^{s}\rightarrow \mathcal{\Theta }_{A_{k-1}}^{s}$ is a $%
A_{k-1}$-valued boxed convolution on $\mathcal{\Theta }_{A_{k-1}}^{s}.$ $%
\square $
\end{proposition}

\strut \strut

Clearly, $E_{k,\,k}=\varphi _{k}:A_{k}\rightarrow A_{k-1}$ is the $A_{k-1}$%
-functional.

\strut

\begin{definition}
Let $B$ be a unital algebra and let

\strut 

\begin{center}
$\Bbb{C}\subset ^{\varphi _{0}}B\subset ^{\varphi _{1}}A_{1}\subset
^{\varphi _{2}}A_{2}\subset ^{\varphi _{3}}A_{3}\subset ^{\varphi _{4}}\cdot
\cdot \cdot $
\end{center}

\strut 

be a tower of amalgamated NCPSpaces. Define an enveloping algebra $A_{\infty
}$ and a linear functional $E_{\infty }:A_{\infty }\rightarrow \Bbb{C}$ by

\strut \strut 

\begin{center}
$E_{\infty }(x)=E_{j}(x),$ \ for $x\in A_{j}\subset A_{\infty },$ for all $\
j\in \Bbb{N}.$
\end{center}

\strut 

We will say the NCPSpace, $\left( A_{\infty },E_{\infty }\right) ,$ is the
enveloping NCPSpace of the given tower.
\end{definition}

\strut

By definition, for any $x\in A_{N}$ and $y\in A_{M}$ and for all $\alpha
,\beta \in \Bbb{C},$ we have that

\strut

\begin{center}
$E_{\infty }\left( \alpha x+\beta y\right) =\alpha E_{\infty }(x)+\beta
E_{\infty }(y)=\alpha E_{N}(x)+\beta E_{M}(y),$
\end{center}

\strut

where $N,M\in \Bbb{N}.$ Then above observation shows that there exits
cumulants of random variables in $(A_{\infty },E_{\infty })$ as follows ;

\strut

\begin{definition}
Let $B$ be a unital algebra and let

\strut 

\begin{center}
$\Bbb{C}\subset ^{\varphi _{0}}B\subset ^{\varphi _{1}}A_{1}\subset
^{\varphi _{2}}A_{2}\subset ^{\varphi _{3}}A_{3}\subset ^{\varphi _{4}}\cdot
\cdot \cdot $
\end{center}

\strut 

be a tower of amalgamated NCPSpaces and $(A_{\infty },E_{\infty })$ an
enveloping NCPSpace of the tower. Let $x_{1},...,x_{s}\in (A_{\infty
},E_{\infty })$ be random variables ($s\in \Bbb{N}$). Define $%
(i_{1},...,i_{n})$-th cumulants of $x_{1},...,x_{s}$ by

\strut 

\begin{center}
$k_{n}^{(E_{\infty })}\left( x_{i_{1}},...,x_{i_{n}}\right) =\underset{\pi
\in NC(n)}{\sum }\left( E_{\infty }\right) _{\pi }\left(
x_{i_{1}},...,x_{i_{n}}\right) \mu (\pi ,1_{n}),$
\end{center}

\strut 

for all $(i_{1},...,i_{n})\in \{1,...,s\}^{n},$ $n\in \Bbb{N},$ where $%
k_{n}^{(E_{\infty })}(...)$ and $(E_{\infty })_{\pi }(...)$ are
scalar-valued cumulants and partition-dependent moments in the sense of
Speicher and Nica.
\end{definition}

\strut

Suppose that $x,y\in (A_{\infty },E_{\infty })$ are random variables and
assume that $x\in A_{j}$ and $y\in A_{k},$ for some $j\leq \,k\in \Bbb{N}.$
Then

\strut

$k_{n}^{(E_{\infty })}\left( x,y,x\right) =\underset{\pi \in NC(3)}{\sum }%
(E_{\infty })_{\pi }\left( x,y,x\right) \mu (\pi ,1_{3})$

$\ \ \ \ \ \ \ \ =E_{\infty }(xyx)-E_{\infty }(x)E_{\infty }(yx)-E_{\infty
}(xy)E_{\infty }(x)$

\begin{center}
$-E_{\infty }\left( xE_{\infty }(y)x\right) +E_{\infty }(x)E_{\infty
}(y)E_{\infty }(x)$
\end{center}

$\ \ \ \ \ \ \ \ =E_{\infty }(xyx)-E_{\infty }(x)E_{\infty }(yx)-E_{\infty
}(xy)E_{\infty }(x)$

\begin{center}
$-E_{\infty }\left( x^{2}\right) E_{\infty }(y)+E_{\infty }(x)E_{\infty
}(y)E_{\infty }(x)$
\end{center}

$\ \ \ \ \ \ \ \ =E_{k}(xyx)-E_{j}(x)E_{k}(yx)-E_{k}(xy)E_{j}(x)$

\begin{center}
$-E_{j}(x^{2})E_{k}(y)+E_{j}(x)E_{k}(y)E_{j}(x).$
\end{center}

\strut

Since $j\leq k$ and $A_{j}\subset ^{E_{j,k}}A_{k},$ a random variable in $%
(A_{\infty },E_{\infty }),$

\begin{center}
$xyx\in A_{k},$ $\ xy\in A_{k}$ \ and \ $yx\in A_{k}.$
\end{center}

\strut \strut

So, we can get the last equality, by definition of $E_{\infty }.$

\strut

\begin{proposition}
(See Chapter 1 or see [6] and [7]) Let $B$ be a unital algebra and let

\strut 

\begin{center}
$\Bbb{C}\subset ^{\varphi _{0}}B\subset ^{\varphi _{1}}A_{1}\subset
^{\varphi _{2}}A_{2}\subset ^{\varphi _{3}}A_{3}\subset ^{\varphi _{4}}\cdot
\cdot \cdot $
\end{center}

\strut 

be a tower of amalgamated NCPSpaces and $\left( A_{\infty },E_{\infty
}\right) $, an enveloping NCPSpace of this tower. Let $X=\{x_{1},...,x_{s}\}$
and $Y=\{y_{1},...,y_{s}\}$ be two free subsets of $(A_{\infty },E_{\infty
}),$ then

\strut 

(2.9.1) $r_{x_{1},...,x_{s}y_{1},...,y_{s}}^{(E_{\infty
})}(z_{1},...,z_{2s})=r_{x_{1},...,x_{s}}^{(E_{\infty
})}(z_{1},...,z_{s})+r_{y_{1},...,y_{s}}^{(E_{\infty })}(z_{s+1},...,z_{2s}).
$

\strut 

(2.9.2) $r_{x_{1}+y_{1},...,x_{s}+y_{s}}^{(E_{\infty
})}(z_{1},...,z_{s})=\left( r_{x_{1},...,x_{s}}^{(E_{\infty
})}+r_{y_{1},...,y_{s}}^{(E_{\infty })}\right) (z_{1},...,z_{s}).$

\strut 

(2.9.3) $r_{x_{1}y_{1},...,x_{s}y_{s}}^{(E_{\infty
})}(z_{1},...,z_{s})=\left( r_{x_{1},...,x_{s}}^{(E_{\infty })}\,\,\frame{*}%
\,\,\,r_{y_{1},...,y_{s}}^{(E_{\infty })}\right) (z_{1},...,z_{s}),$

\strut 

where \frame{*}$\,:\Theta _{s}\times \Theta _{s}\rightarrow \Theta _{s}$ is
the (scalar-valued) boxed convolution defined by Nica and Speicher. $\square 
$
\end{proposition}

\strut

Now, we will observe the freeness on $(A_{\infty },E_{\infty }),$ more in
detail.

\strut

\begin{theorem}
Let $B$ be a unital algebra and let

\strut 

\begin{center}
$\Bbb{C}\subset ^{\varphi _{0}}B\subset ^{\varphi _{1}}A_{1}\subset
^{\varphi _{2}}A_{2}\subset ^{\varphi _{3}}A_{3}\subset ^{\varphi _{4}}\cdot
\cdot \cdot $
\end{center}

\strut 

be a tower of amalgamated NCPSpaces and $(A_{\infty },E_{\infty }),$ an
enveloping NCPSpace. Let $X=\{x_{1},...,x_{s}\}\subset \cup _{j=1}^{k}A_{j}$
and $Y=\{y_{1},...,y_{s}\}\subset \cup _{j=k+1}^{l}A_{j}$ be subsets of
operators ($s\in \Bbb{N}$), where $k<l,$ in $\Bbb{N}.$ If $X$ and $Y$ are
free in $(A_{l},\,E_{l})$, then (2.9.1), (2.9.2) and (2.9.3) holds true in $%
\left( A_{\infty },E_{\infty }\right) .$ More precisely,

\strut 

(1) $r_{x_{1},...,x_{s},y_{1},...,y_{s}}^{(E_{\infty
})}(z_{1},...,z_{2s})=r_{x_{1},...,x_{s}}^{(E_{k})}(z_{1},...,z_{s})+r_{y_{1},...,y_{s}}^{(E_{l})}(z_{s+1},...,z_{2s}).
$

\strut 

(2) $r_{x_{1}+y_{1},...,x_{s}+y_{s}}^{(E_{\infty })}(z_{1},...,z_{s})=\left(
r_{x_{1},...,x_{s}}^{(E_{k})}+r_{y_{1},...,y_{s}}^{(E_{l})}\right)
(z_{1},...,z_{s}).$

\strut 

(3) $r_{x_{1}y_{1},...,x_{s}y_{s}}^{(E_{\infty })}(z_{1},...,z_{s})=\left(
r_{x_{1},...,x_{s}}^{(E_{k})}\,\,\frame{*}_{\Bbb{C}}\,\,%
\,r_{y_{1},...,y_{s}}^{(E_{l})}\right) (z_{1},...,z_{s}).$
\end{theorem}

\strut

\begin{proof}
We can regards $X$ and $Y$ in $(A_{\infty },E_{\infty })$ as subsets of
random variables in $(A_{l},E_{l}).$ We have that

\strut

$\ \ \ \ r_{x_{1},...,x_{s},y_{1},...,y_{s}}^{(E_{\infty
})}(z_{1},...,z_{2s})=r_{x_{1},...,x_{s},y_{1}...,y_{s}}^{(E_{l})}(z_{1},...,z_{2s}) 
$

\strut

by definition of $E_{\infty }$

\strut

$\ \ \ \ \ \ \ \ \ \ \ \ \ \ \ \ \ \ \ \ \ \ \ \ \
=r_{x_{1},...,x_{s}}^{(E_{l})}(z_{1},...,z_{s})+r_{y_{1},...,y_{s}}^{(E_{l})}(z_{s+1},...,z_{2s}) 
$

\strut

by the freeness of $X$ and $Y,$ in $(A_{l},E_{l})$

\strut

$\ \ \ \ \ \ \ \ \ \ \ \ \ \ \ \ \ \ \ \ \ \ \ \ \
=r_{x_{1},...,x_{s}}^{(E_{k})}(z_{1},...,z_{s})+r_{y_{1},...,y_{s}}^{(E_{l})}(z_{s+1},...,z_{2s}) 
$

\strut

since $E_{l}=\varphi _{0}\varphi _{1}...\varphi _{k}\varphi _{k+1}...\varphi
_{l}$ and $X\subset (A_{k},E_{k}).$ Similarly, we have that

\strut

\begin{center}
$
\begin{array}{ll}
r_{x_{1}+y_{1},...,x_{s}+y_{s}}^{(E_{\infty })}(z_{1},...,z_{s}) & 
=r_{x_{1}+y_{1},...,x_{s}+y_{s}}^{(E_{l})}(z_{1},...,z_{s}) \\ 
& =\left( r_{x_{1},...,x_{s}}^{(E_{l})}+r_{y_{1},...,y_{s}}^{(E_{l})}\right)
(z_{1},...,z_{s}) \\ 
& =\left( r_{x_{1},...,x_{s}}^{(E_{k})}+r_{y_{1},...,y_{s}}^{(E_{l})}\right)
(z_{1},...,z_{s}).
\end{array}
$
\end{center}

\strut

Now, fix $(i_{1},...,i_{n})\in \{1,...,s\}^{n},$ $n\in \Bbb{N}.$ Notice that 
$x_{1}y_{1},...,x_{s}y_{s}$ can be regarded as random variables in $%
(A_{l},E_{l}).$ Consider the $(i_{1},...,i_{n})$-th coefficient of $%
r_{x_{1}y_{1},...,x_{s}y_{s}}^{(E_{\infty })}$ ;

\strut

$\ \ \ \ coef_{i_{1},...,i_{n}}\left(
r_{x_{1}y_{1},...,x_{s}y_{s}}^{(E_{\infty })}\right)
=coef_{i_{1},...,i_{n}}\left( r_{x_{1}y_{1},...,x_{s}y_{s}}^{(E_{l})}\right) 
$

\strut

by the fact that $x_{j}y_{j}\in A_{l},$ for all $j$ and by the definition of 
$E_{\infty }$

\strut

$\ \ \ \ \ \ \ \ \ \ \ \ \ \ =k_{n}^{(E_{l})}\left(
x_{i_{1}}y_{i_{1}},...,x_{i_{n}}y_{i_{n}}\right) $

$\ \ \ \ \ \ \ \ \ \ \ \ \ \ =\underset{\pi \in NC(n)}{\sum }k_{\pi }\left(
x_{i_{1}},...,x_{i_{n}}\right) k_{Kr(\pi )}\left(
y_{i_{1}},...,y_{i_{n}}\right) $

\strut

by the freeness of $X$ and $Y$ in $(A_{l},E_{l})$, in the sense of Speicher
and Nica (See Chapter 1 or see [6] and [7])

\strut

$\ \ \ \ \ \ \ \ \ \ \ \ \ \ =coef_{i_{1},...,i_{n}}\left(
r_{x_{1},...,x_{s}}^{(E_{l})}\,\,\,\frame{*}\,\,%
\,r_{y_{1},...,y_{s}}^{(E_{l})}\right) $

$\ \ \ \ \ \ \ \ \ \ \ \ \ \ =coef_{i_{1},...,i_{n}}\left(
r_{x_{1},...,x_{s}}^{(E_{k})}\,\,\,\frame{*}\,\,%
\,r_{y_{1},...,y_{s}}^{(E_{l})}\right) ,$

\strut

since $x_{1},...,x_{s}\in \cup _{j=1}^{k}A_{j}$ \ and $E_{l}=E_{k}\varphi
_{k+1}...\varphi _{l}.$
\end{proof}

\strut

\begin{corollary}
Let $B$ be a unital algebra and suppose that we have a tower of amalgamated
NCPSpaces, $\Bbb{C}\subset ^{\varphi _{0}}B\subset ^{\varphi
_{1}}A_{1}\subset ^{\varphi _{2}}A_{2}\subset ^{\varphi _{3}}\cdot \cdot
\cdot .$ Let $X$ and $Y$ be subsets in $A_{k}.$ If $X$ and $Y$ are free in $%
(A_{k},E_{k}),$ then $X$ and $Y$ are free in $(A_{N},E_{N}),$ for all $N>k.$
\ $\square $
\end{corollary}

\strut

\begin{remark}
Suppose that we choose $x_{1},...,x_{s}\in \cup _{j=l+1}^{k}A_{j}$ and $%
y_{1},...,y_{s}\in \cup _{j=1}^{l}A_{j},$ where $k>l,$ in this time. Then we
can get the same results as above, but it means different.
\end{remark}

\strut

\begin{remark}
If we take a sufficiently large $N\in \Bbb{N},$ and if we replace $\infty $
to such $N,$ then we can get the same results as in the previous theorem.
\end{remark}

\strut

More generally, we have that ;

\strut

\begin{theorem}
Let $B$ be a unital algebra and let

\strut 

\begin{center}
$\Bbb{C}\subset ^{\varphi _{0}}B\subset ^{\varphi _{1}}A_{1}\subset
^{\varphi _{2}}A_{2}\subset ^{\varphi _{3}}A_{3}\subset ^{\varphi _{4}}\cdot
\cdot \cdot $
\end{center}

\strut 

be a tower of amalgamated NCPSpaces. Fix a sufficiently large $N\in \Bbb{N}$
and fix $k<N$. Consider the conditional expectation (or a $A_{k-1}$%
-functional), $E_{k,\,N}:A_{N}\rightarrow A_{k-1}.$ \strut Let $%
X=\{x_{1},...,x_{s}\}\subset \cup _{j=k}^{m}A_{j}$ and $Y=\{y_{1},...,y_{s}%
\}\subset \cup _{j=m+1}^{n}A_{j}$ be two subsets, in $(A_{N},\,E_{k,N})$ of $%
A_{k-1}$-valued random variables ($s\in \Bbb{N}$), where $k<m<n<N.$ If $X$
and $Y$ are free over $A_{k-1},$ in $(A_{N},E_{k,N}),$ then

\strut 

(1) $%
R_{x_{1},...,x_{s},y_{1},...,y_{s}}^{(E_{k,N})}(z_{1},...,z_{2s})=R_{x_{1},...,x_{s}}^{(E_{k,m})}(z_{1},...,z_{s})+R_{y_{1},...,y_{s}}^{(E_{k,n})}(z_{s+1},...,z_{2s}).
$

\strut 

(2) $R_{x_{1}+y_{1},...,x_{s}+y_{s}}^{(E_{k,N})}(z_{1},...,z_{s})=\left(
R_{x_{1},...,x_{s}}^{(E_{k,m})}+R_{y_{1},...,y_{s}}^{(E_{k,n})}\right)
(z_{1},...,z_{s}).$

\strut 

(3) $R_{x_{1}y_{1},...,x_{s}y_{s}}^{(E_{k,N})}(z_{1},...,z_{s})=\left(
R_{x_{1},...,x_{s}}^{(E_{k,m})}\,\,\frame{*}_{A_{k-1}}\,%
\,R_{y_{1},...,y_{s}}^{(E_{k,n})\,:\,t}\right) (z_{1},...,z_{s}).$
\end{theorem}

\strut

\begin{proof}
We can regard $X$ and $Y$ are subsets of $A_{N}$ and consider them as $%
A_{k-1}$-valued random variables in $(A_{N},\,E_{k,N}).$ Since $X$ and $Y$
are free over $A_{k-1},$ we have that

\strut

\begin{center}
$%
R_{x_{1},...,x_{s},y_{1},...,y_{s}}^{(E_{k,N})}(z_{1},...,z_{2s})=R_{x_{1},...,x_{s}}^{(E_{k,N})}(z_{1},...,z_{s})+R_{y_{1},...,y_{s}}^{(E_{k,N})}(z_{s+1},...,z_{2s}) 
$
\end{center}

\strut

Observe that

\strut

$\ \ \ \ coef_{i_{1},...,i_{n}}\left( R_{x_{1},...,x_{s}}^{(E_{k,N})}\right)
=k_{n}^{(E_{k,N})}\left( x_{i_{1}},...,x_{i_{n}}\right) $

$\ \ \ \ \ \ \ \ \ \ \ \ \ \ \ \ \ \ \ \ \ =c_{E_{k,N}}^{(n)}\left(
x_{i_{1}}\otimes a_{i_{2}}x_{i_{2}}\otimes ...\otimes
a_{i_{n}}x_{i_{n}}\right) $

$\ \ \ \ \ \ \ \ \ \ \ \ \ \ \ \ \ \ \ \ \ =\underset{\pi \in NC(n)}{\sum }%
\widehat{E_{k,N}}(\pi )\left( x_{i_{1}}\otimes a_{i_{2}}x_{i_{2}}\otimes
...\otimes a_{i_{n}}x_{i_{n}}\right) \mu (\pi ,1_{n})$

$\ \ \ \ \ \ \ \ \ \ \ \ \ \ \ \ \ \ \ \ \ =\underset{\pi \in NC(n)}{\sum }%
\widehat{E_{k,m}}(\pi )\left( x_{i_{1}}\otimes a_{i_{2}}x_{i_{2}}\otimes
...\otimes a_{i_{n}}x_{i_{n}}\right) \,\mu (\pi ,1_{n}),$

\strut

since $E_{k,N}=E_{k,m}\varphi _{m+1}...\varphi _{n}...\varphi _{N}$ and $%
x_{i_{j}}$'s are in $A_{m},$ where $a_{i_{2}},...,a_{i_{n}}\in A_{k-1}$ are
arbitrary. Thus

\strut

\begin{center}
$coef_{i_{1},...,i_{n}}\left( R_{x_{1},...,x_{s}}^{(E_{k,N})}\right)
=coef_{i_{1},...,i_{n}}\left( R_{x_{1},...,x_{s}}^{(E_{k,m})}\right) .$
\end{center}

\strut

Similarly,

\strut

\begin{center}
$coef_{i_{1},...,i_{n}}\left( R_{y_{1},...,y_{s}}^{(E_{k,N})}\right)
=coef_{i_{1},...,i_{n}}\left( R_{y_{1},...,y_{s}}^{(E_{k,n})}\right) .$
\end{center}

\strut

So, (1) is proved. By using the same idea, we can easily prove (2).

\strut

Now, fix $(i_{1},...,i_{n})\in \{1,...,s\}^{n},$ $n\in \Bbb{N}.$ Then

\strut

$\ \ \ coef_{i_{1},...,i_{p}}\left(
R_{x_{1}y_{1},...,x_{s}y_{s}}^{(E_{k,N})}\right) =k_{p}^{(E_{k,N})}\left(
x_{i_{1}}y_{i_{1}},...,x_{i_{p}}y_{i_{p}}\right) $

$\ \ \ \ \ \ \ \ \ =c_{E_{k,N}}^{(p)}\left( x_{i_{1}}y_{i_{1}}\otimes
a_{i_{2}}x_{i_{2}}y_{i_{2}}\otimes ...\otimes
a_{i_{p}}x_{i_{p}}y_{i_{p}}\right) $

$\ \ \ \ \ \ \ \ \ =\underset{\pi \in NC(p)}{\sum }\widehat{E_{k,N}}(\pi
)\left( x_{i_{1}}y_{i_{1}}\otimes a_{i_{2}}x_{i_{2}}y_{i_{2}}\otimes
...\otimes a_{i_{p}}x_{i_{p}}y_{i_{p}}\right) \,\mu (\pi ,1_{p})$

$\ \ \ \ \ \ \ \ \ =\underset{\pi \in NC(p)}{\sum }\left( \widehat{%
c_{E_{k,m}}}\oplus \widehat{c_{E_{k,n}}}\right) \left( x_{i_{1}}\otimes
y_{i_{1}}\otimes a_{i_{2}}x_{i_{2}}\otimes y_{i_{2}}\otimes ...\otimes
a_{i_{p}}x_{i_{p}}\otimes y_{i_{p}}\right) $

\strut

by the $A_{k-1}$-freeness of $X$ and $Y$ in $(A_{N},\,E_{k,N})$

\strut

$\ \ \ \ \ \ \ \ \ =coef_{i_{1},...,i_{n}}\left(
R_{x_{1},...,x_{s}}^{(E_{n,m})}\,\,\,\frame{*}_{A_{k-1}}\,\,%
\,R_{y_{1},...,y_{s}}^{(E_{k,n})\,\,:\,\,t}\right) .$

\strut

Therefore,

\strut

\begin{center}
$R_{x_{1}y_{1},...,x_{s}y_{s}}^{(E_{k,N})}=R_{x_{1},...,x_{s}}^{(E_{k,m})}\,%
\,\,\frame{*}_{B}\,\,\,R_{y_{1},...,y_{s}}^{(E_{k,n})}.$
\end{center}
\end{proof}

\strut

Fix $k\in \Bbb{N}.$ We can define an enveloping NCPSpace over $A_{k-1},$ $%
\left( A_{\infty },\,E_{k,\,\infty }\right) ,$ where $A_{\infty }$ is an
enveloping algebra of the tower and $E_{k,\,\infty }:A_{\infty }\rightarrow
A_{k-1}$ is defined by

\strut

\begin{center}
$E_{k,\,\infty }(x)\overset{def}{=}E_{k,\,j}(x),$ for all $x\in A_{j},$
\end{center}

\strut

where $j>k.$

\strut

\begin{corollary}
Let $B$ be a unital algebra and let

\strut 

\begin{center}
$\Bbb{C}\subset ^{\varphi _{0}}B\subset ^{\varphi _{1}}A_{1}\subset
^{\varphi _{2}}A_{2}\subset ^{\varphi _{3}}A_{3}\subset ^{\varphi _{4}}\cdot
\cdot \cdot $
\end{center}

\strut 

be a tower of amalgamated NCPSpaces. Let $X$ and $Y$ be a two subsets of $%
A_{k}.$ If $X$ and $Y$ are free over $A_{k-1},$ in $(A_{k},$\thinspace $%
E_{k,k}),$ then $X$ and $Y$ are free over $A_{k-1},$ in $(A_{N},E_{k,N}),$
for all $N>k,$ in $\Bbb{N}.$ $\square $
\end{corollary}

\strut

\begin{remark}
Let $(A_{\infty },\,E_{k,\,\infty })$ be an enveloping NCPSpace over $%
A_{k-1}.$ Then the same results with the previous theorem holds true when we
replace $N$ by $\infty .$
\end{remark}

\strut

\strut

\strut

\subsection{Distributions On a Tower of Algebras}

\strut

\strut

\strut

In this section, we will consider (scalar-valued or operator-valued)
distributions when the tower of amalgamated NCPSpaces is given. Throughout
this section, we let

\strut

\begin{center}
$\Bbb{C}\subset ^{\varphi _{0}}B\subset ^{\varphi _{1}}A_{1}\subset
^{\varphi _{2}}A_{2}\subset ^{\varphi _{3}}A_{3}\subset ^{\varphi _{4}}\cdot
\cdot \cdot $
\end{center}

\strut

be a tower of amalgamated NCPSpaces. Also, like in the previous section, we
will use the following notations ;

\strut

\begin{center}
$E_{N}=\varphi _{0}\varphi _{1}...\varphi _{N}:A_{N}\rightarrow \Bbb{C},$ a
linear functional.
\end{center}

For $k<N,$

\begin{center}
$E_{k,N}=\varphi _{k}\varphi _{k+1}...\varphi _{N}:A_{N}\rightarrow A_{k-1},$
a $A_{k-1}$-functional
\end{center}

with

\begin{center}
$E_{1,N}:A_{N}\rightarrow B,$ a $B$-functional.
\end{center}

\strut

In the previous section, we considered the R-transform calculus of $E_{N}$
and $E_{k,\,N}.$ Also, when we consider the enveloping algebra of the given
tower, $A_{\infty },$ we could define a linear functional and $A_{k-1}$%
-functional, $E_{\infty }:A_{\infty }\rightarrow \Bbb{C}$ and $E_{k,\,\infty
}:A_{\infty }\rightarrow A_{k-1},$ respectively. Also, R-transform calculus
related to them are given in the previous section. Define

\strut

\begin{center}
$\sum_{s}=\{\sigma :\Bbb{C}[X_{1},...,X_{s}]\rightarrow \Bbb{C}:\sigma $ is
a $\Bbb{C}$-valued distribution$\}$
\end{center}

and

\begin{center}
$\sum_{A_{k-1}}^{s}=\{\sigma
^{(A_{k-1})}:A_{k-1}[X_{1},...,X_{s}]\rightarrow A_{k-1}:\sigma ^{(A_{k-1})}$
is a $A_{k-1}$-functional$\},$
\end{center}

\strut

for the fixed $s\in \Bbb{N},$ where $X_{1},...,X_{s}$ are noncommutative
indeterminents. If $\sigma \in \sum_{s},$ then there exists an arbitrary
NCPSpace $(A,\varphi )$ and random variables $x_{1},...,x_{s}\in (A,\varphi
) $ such that

\strut

\begin{center}
$\sigma \left( P\right) =\varphi \left( P(x_{1},...,x_{s})\right) \in \Bbb{C}%
,$ for all $P\in C[X_{1},...,X_{s}].$
\end{center}

\strut

And if $\sigma ^{(A_{k-1})}\in \sum_{A_{k-1}}^{s},$ then there exsits an
arbitrary NCPSpace over $A_{k-1},$ $(D,\psi )$, and $A_{k-1}$-valued random
variables $y_{1},...,y_{s}\in (D,\psi )$ such that

\strut

\begin{center}
$\sigma ^{(A_{k-1})}(Q)=\psi \left( Q(y_{1},...,y_{s})\right) \in A_{k-1},$
for all $Q\in A_{k-1}[X_{1},...,X_{s}].$
\end{center}

\strut

In particular, we denote such $\sigma $ and $\sigma ^{(A_{k-1})}$ by

\strut

\begin{center}
$\sigma _{x_{1},...,x_{s}\in (A,\varphi )}$ \ \ and \ \ $\sigma
_{y_{1},...,y_{s}\in (D,\psi )}^{(A_{k-1})},$
\end{center}

\strut

respectively. Notice that, for the given $\sigma $ and $\sigma ^{(A_{k-1})},$
$x_{1},...,x_{s}\in (A,\varphi )$ and $y_{1},...,y_{s}\in (D,\psi )$ are Not
uniquely determined. So, we can regard them as equivalence class under the
identically-distributedness. i.e

\strut

$\sigma =\sigma _{x_{1},...,x_{s}\in (A,\varphi )}\sim \sigma
_{x_{1}^{\prime },...,x_{s}^{\prime }\in (A^{\prime },\varphi ^{\prime })}$

\begin{center}
$\Longleftrightarrow \varphi \left( P(x_{1},...,x_{s})\right) =\varphi
^{\prime }\left( P(x_{1}^{\prime },...,x_{s}^{\prime })\right) ,$ $\forall $ 
$P\in \Bbb{C}[X_{1},...,X_{s}]$
\end{center}

and

\strut

$\sigma ^{(A_{k-1})}=\sigma _{y_{1},...,y_{s}\in (D,\psi )}^{(A_{k-1})}\sim
\sigma _{y_{1}^{\prime },...,y_{s}^{\prime }\in (D^{\prime },\psi ^{\prime
})}^{(A_{k-1})}$

\begin{center}
$\Longleftrightarrow \psi \left( Q(y_{1},...,y_{s})\right) =\psi ^{\prime
}\left( Q(y_{1}^{\prime },...,y_{s}^{\prime })\right) ,$ $\forall Q\in
A_{k-1}[X_{1},...,X_{s}].$
\end{center}

\strut

Recall the definitions ;

\ \ \ \ \ \ \ \ \ \ \ \ the free additive convolution $\boxplus
:\sum_{s}\times \sum_{s}\rightarrow \sum_{s},$

\ \ \ \ \ \ \ \ \ \ \ \ the free multiplicative convolution $\boxtimes
:\sum_{s}\times \sum_{s}\rightarrow \sum_{s}$

and

the oprator-valued free additive convolution $\boxplus
_{A_{k-1}}:\sum_{A_{k-1}}^{s}\times \sum_{A_{k-1}}^{s}\rightarrow
\sum_{A_{k-1}}^{s}.$

\strut

(See [8] for $\boxplus $ and $\boxtimes $. See [1] for $\boxplus _{A_{k-1}}$%
) Now, we will define the oprerator-valued free multiplicative convolution
of operator-valued distributions

\strut

\begin{center}
$\boxtimes _{A_{k-1}}:\sum_{A_{k-1}}^{s}\times \sum_{A_{k-1}}^{s}\rightarrow
\sum_{A_{k-1}}^{s}$
\end{center}

by

\begin{center}
$\sigma _{y_{1},...,y_{s}\in (D,\psi )}^{(A_{k-1})}\,\,\boxtimes
_{A_{k-1}}\,\,\sigma _{p_{1},...,p_{s}\in (D^{\prime },\psi ^{\prime
})}^{(A_{k-1})}=\sigma _{y_{1}p_{1},...,y_{s}p_{s}\in \left(
D*_{A_{k-1}}D^{\prime },\,\,\psi *\psi ^{\prime }\right) }^{(A_{k-1})},$
\end{center}

\strut \strut

where $\left( D*_{A_{k-1}}D^{\prime },\,\psi *\psi \right) $ is the free
product of $(D,\psi )$ and $(D^{\prime },\psi ^{\prime }),$ with
amalgamation over $A_{k-1}.$

\strut

Rest of this section, we will consider the $B$-valued distributions
depending on our tower.

\strut \strut

\begin{proposition}
Let $\sigma \in \sum_{s}$ and $\sigma ^{(A_{k-1})}\in \sum_{A_{k-1}}^{s}.$
Then

\strut 

(1) $\sigma _{x_{1},...,x_{s}\in (A_{\infty },E_{\infty })}\,\,\boxplus
\,\sigma _{y_{1},...,y_{s}\in (A_{\infty },E_{\infty })}=\sigma
_{x_{1}+y_{1},...,x_{s}+y_{s}\in (A_{\infty },E_{\infty })}.$

\strut 

(2) $\sigma _{x_{1},...,x_{s}\in (A_{\infty },E_{k,\infty
})}^{(A_{k-1})}\,\,\boxplus _{A_{k-1}}\,\,\sigma _{y_{1},...,y_{s}\in
(A_{\infty },E_{k,\infty })}^{(A_{k-1})}=\sigma
_{x_{1}+y_{1},...,x_{s}+y_{s}\in (A_{\infty },E_{k,\infty })}^{(A_{k-1})}.$

\strut 

(3) $\sigma _{x_{1},...,x_{s}\in (A_{j},E_{j})}\,\,\boxplus \,\,\sigma
_{y_{1},...,y_{s}\in (A_{k},E_{k})}=\sigma _{x_{1}+y_{1},...,x_{s}+y_{s}\in
(A_{k},E_{k})},$ if $\,j<k.$

$\strut $

(4) $\sigma _{x_{1},...,x_{s}\in (A_{i},\,E_{k,i})}^{(A_{k-1})}\,\,\boxplus
_{A_{k-1}}\sigma _{y_{1},...,y_{s}\in (A_{j},\,E_{k,j})}^{(A_{k-1})}=\sigma
_{x_{1}+y_{1},...,x_{s}+y_{s}\in (A_{j},E_{k,j})}^{(A_{k-1})},$

$\strut $

\ \ \ \ \ if $K<i<j.\,\,\ \ \square $
\end{proposition}

\strut

\begin{proposition}
Let $\sigma \in \sum_{s}$ and $\sigma ^{(A_{k-1})}\in \sum_{A_{k-1}}.$ Then

\strut 

(1) $\sigma _{x_{1},...,x_{s}\in (A_{\infty },E_{\infty })}\,\,\boxtimes
\,\sigma _{y_{1},...,y_{s}\in (A_{\infty },E_{\infty })}=\sigma
_{x_{1}y_{1},...,x_{s}y_{s}\in (A_{\infty },E_{\infty })}.$

\strut 

(2) $\sigma _{x_{1},...,x_{s}\in (A_{\infty },E_{k,\infty
})}^{(A_{k-1})}\,\,\boxtimes _{A_{k-1}}\,\sigma _{y_{1},...,y_{s}\in
(A_{\infty },E_{k,\infty })}^{(A_{k-1})}=\sigma
_{x_{1}y_{1},...,x_{s}y_{s}\in (A_{\infty },E_{k,\infty })}^{(A_{k-1})}.$

\strut 

(3) $\sigma _{x_{1},...,x_{s}\in (A_{i},E_{i})}\,\,\,\boxtimes \,\sigma
_{y_{1},...,y_{s}\in (A_{j},E_{j})}=\sigma _{x_{1}y_{1},...,x_{s}y_{s}\in
(A_{j},E_{j})},$ if $i<j.$

\strut 

(4) $\sigma _{x_{1},...,x_{s}\in (A_{i},E_{k,i})}^{(A_{k-1})}\,\,\,\boxtimes
_{A_{k-1}}\,\sigma _{y_{1},...,y_{s}\in (A_{j},E_{k,j})}^{(A_{k-1})}=\sigma
_{x_{1}y_{1},...,x_{s}y_{s}\in (A_{j},\,E_{k,j})}^{(A_{k-1})},$

\strut 

\ \ \ \ if $k<i<j.$
\end{proposition}

\strut

\begin{proof}
Suppose that $X=\{x_{1},...,x_{s}\}$ and $Y=\{y_{1},...,y_{s}\}$ are two
subsets of $(A_{\infty },E_{\infty })$ which are free. Then, there exists $%
i,j\in N$ such that

\strut

\begin{center}
$X\subset \cup _{l=1}^{m}A_{l}\subseteq A_{i}$ \ \ and \ \ $Y\subset \cup
_{l=1}^{n}A_{l}\subseteq A_{j}.$
\end{center}

\strut

Assume that $m<n.$ Then $\{x_{1}y_{1},...,x_{s}y_{s}\}\subset (A_{j},E_{j}).$
By freeness,

\strut

\begin{center}
$\sigma _{x_{1},...,x_{s}\in (A_{i},E_{i})}\,\,\boxtimes \,\sigma
_{y_{1},...,y_{s}\in (A_{j},E_{j})}=\sigma _{x_{1}y_{1},...,x_{s}y_{s}\in
(A_{j},E_{j})}.$
\end{center}

\strut

So, (1) and (3) are proved. Similarly, $\{x_{1}y_{1},...,x_{s}y_{s}\}\subset
(A_{j},E_{k,j}),$ if $k<i<j.$ So,

\strut

\begin{center}
$\sigma _{x_{1},...,x_{s}\in (A_{i},E_{k,i})}^{(A_{k-1})}\,\,\,\boxtimes
_{A_{k-1}}\,\,\sigma _{y_{1},...,y_{s}\in
(A_{j},E_{k,\,j})}^{(A_{k-1})}=\sigma _{x_{1}y_{1},...,x_{s}y_{s}\in
(A_{j},E_{k,j})}^{(A_{k-1})}.$
\end{center}

\strut

Therefore, we can prove (2) and (4).
\end{proof}

\strut \strut

Now, we will observe the following concept introduced in [18], so-called
''compatibility''. The R-transform theory, under this compatibility is
studied in [12].

\strut

\begin{definition}
Let $B$ be a unital algebra and $A,$ an algebra over $B.$ Let $(A,\varphi )$
be a NCPSpace and let $(A,E)$ be a NCPSpace over $B,$ with its $B$%
-functional $E:A\rightarrow B.$ We say that $(A,\varphi )$ and $(A,E)$ are
compatible if

\strut 

\begin{center}
$\varphi (x)=\varphi \left( E(x)\right) ,$ \ for all $x\in A.$
\end{center}
\end{definition}

\strut

As we can see in [12], even under the compatibility, to compute
scalar-valued R-transforms from the operator-valued R-transforms or
operator-valued moment series is complicated. But [12] provides the method.

\strut

\begin{theorem}
Let $B$ be a unital algebra and let

\strut 

\begin{center}
$\Bbb{C}\subset ^{\varphi _{0}}B\subset ^{\varphi _{1}}A_{1}\subset
^{\varphi _{2}}A_{2}\subset ^{\varphi _{3}}A_{3}\subset ^{\varphi _{4}}\cdot
\cdot \cdot $
\end{center}

\strut 

be a tower of amalgamated NCPSpaces. Then, for any fixed $k\in \Bbb{N},$ $%
\left( A_{N},\,E_{k,N}\right) $ and $\left( A_{N},\,E_{N}\right) $ are
compatible. In particular, if we have the enveloping algebra of the tower, $%
A_{\infty },$ then $\left( A_{\infty },E_{k,\infty }\right) $ and $\left(
A_{\infty },E_{\infty }\right) $ are compatible, for all $k.$
\end{theorem}

\strut

\begin{proof}
Fix $N\in \Bbb{N}.$ By definition, we have that

\strut

\begin{center}
$
\begin{array}{ll}
E_{N}(x) & =\varphi _{0}\varphi _{1}...\varphi _{N}(x) \\ 
& =(\varphi _{0}\varphi _{1}...\varphi _{k-1})(\varphi _{k}...\varphi
_{N})(x) \\ 
& =E_{k-1}\left( E_{k,N}(x)\right)
\end{array}
$
\end{center}

and

\begin{center}
$
\begin{array}{ll}
E_{N}\left( E_{k,N}(x)\right) & =\left( \varphi _{0}\varphi _{1}...\varphi
_{N}\right) \left( \varphi _{k}...\varphi _{N}(x)\right) \\ 
& =\varphi _{0}\varphi _{1}...\varphi _{N}(x)=E_{k-1}\left(
E_{k,N}(x)\right) ,
\end{array}
$
\end{center}

\strut

for all $x\in A_{N}.$ Therefore, for all $x$ in $A_{N},$

\strut

\begin{center}
$E_{\infty }(x)=E_{N}\left( E_{k,N}(x)\right) .$
\end{center}
\end{proof}

\strut

The above theorem shows that our tower of amalgamated NCPSpaces has nice
compatible properties.

\strut

\strut

\strut

\subsection{Compressed R-transform Theory For the Chain of Projections}

\strut

\strut

\strut

In this section, we will observe a tower of amalgamated NCPSpaces, generated
by a chain of projections $(p_{n})_{n=1}^{\infty }$ of a NCPSpace $%
(A,\varphi )$. The $B$-valued R-transform theory of compressed algebras
(i.e, $B\subset pAp,$ where $p\in A$ is a projection) and $B$-valued
freeness of them are studied in [10]. Also, by [35], there is the
compatibility of $(A_{k}\subset A,\,E_{k+1,\,\infty }\,)$ and $%
(A,\,E_{k}\circ E_{k+1,\,\infty }).$ It is possible that there are only
finitely many projections in the given chain (for example, an algebra $A$
can be a finite dimensional algebra over a finite dimensional algebra $B$).

\strut

\begin{definition}
Let $A$ be a unital algebra. A nonzero element $p\in A\setminus \{0_{A}\}$
is called a projection (or an idempotent element) if $p$ satisfies

\strut 

\begin{center}
$p^{2}=p.$
\end{center}

\strut 

By $A_{pro},$ we will denote a subset of all projections in an algebra, $A.$
It is said that $(p_{n})_{n=1}^{\infty }\subset A_{pro}$ is a chain of
projections if

\strut 

\begin{center}
$\Bbb{C}\subset p_{1}Ap_{1}\subset p_{2}Ap_{2}\subset p_{3}Ap_{3}\subset
\cdot \cdot \cdot $
\end{center}

\strut 

is a tower of compressed subalgebras of $A$ satisfying the property that $%
p_{j}Ap_{j}$ are subalgebras of $p_{j+1}Ap_{j+1},$ for all $\ j\in \Bbb{N}.$
For the convinience of using notation, we will denote maximal projection in
a chain of projections $(p_{n})_{n=1}^{\infty },$ by $p_{\infty }.$ If the
given chain $(p_{n})_{n=1}^{\infty }$ is a finite sequence, then $p_{\infty
}=p_{N},$ for some $N\in \Bbb{N}.$
\end{definition}

\strut \strut

\begin{lemma}
Let $A$ be a unital algebra and let $p\in A_{pro}.$ Then $px=xp,$ for all $%
x\in pAp.$
\end{lemma}

\strut

\begin{proof}
Let $x\in pAp.$ Then there exists an element $a\in A$ such that $x=pap.$ Then

\strut

\begin{center}
$px=p(pap)=p^{2}ap=pap=x$
\end{center}

and

\begin{center}
$xp=(pap)p=pap^{2}=pap=x.$
\end{center}
\end{proof}

\strut

\begin{lemma}
Let $A$ be a unital algebra and let $p,q\in A_{pro}$ satisfying that $pAp$
is a subalgebra of $qAq.$ Then $pq=qp=p.$
\end{lemma}

\strut

\begin{proof}
By assumption, $pAp\subset qAq.$ Let $x\in pAp.$ Then, by the previous lemma,

\strut

\begin{center}
$xq=qx$.
\end{center}

\strut

Since $p\in pAp,$ we can get that

\strut

\begin{center}
$pq=qp.$
\end{center}

\strut

Now, observe that

\strut

\begin{center}
$pq=p^{3}q=p\left( p^{2}q\right) =p\left( p(pq)\right) =p\left( pqp\right)
=p $
\end{center}

\strut

since $pAp\subset qAq.$
\end{proof}

\strut

\begin{corollary}
Let $B$ be a unital algebra and $A,$ an algebra over $B.$ Let $%
(p_{n})_{n=1}^{\infty }\subset A_{pro}$ be a chain of projections. Then, for
any fixed $j<N\in \Bbb{N},$ $\left( p_{j}p_{j+1}...p_{N}\right) A\left(
p_{j}p_{j+1}...p_{N}\right) =p_{j}Ap_{j}.$
\end{corollary}

\strut

\begin{proof}
By the previous lemmas, we have that

\strut

\begin{center}
$
\begin{array}{ll}
\left( p_{j}p_{j+1}...p_{N}\right) A\left( p_{j}p_{j+1}...p_{N}\right) & 
=p_{j}p_{j+1}...p_{N}Ap_{N}...p_{j+1}p_{j} \\ 
& =p_{j}p_{j+1}...\left( p_{N}Ap_{N}\right) ...p_{j+1}p_{j} \\ 
& =p_{j}...\left( p_{N-1}Ap_{N-1}\right) ...p_{j} \\ 
& =p_{j}Ap_{j}.
\end{array}
$
\end{center}
\end{proof}

\strut

\begin{definition}
Let $B$ be a unital algebra and $A,$ an algebra over $B.$ Let $(A,\varphi )$
be a NCPSpace over $B,$ with its $B$-functional $\varphi :A\rightarrow B.$
Suppose that we have a chain of projections $(p_{k})_{k=1}^{\infty }\subset
A_{pro}$ satisfying

\strut 

\begin{center}
$\Bbb{C}\subset ^{\varphi _{0}}B\subset p_{1}Ap_{1}\subset
p_{2}Ap_{2}\subset ....\subset p_{\infty }Ap_{\infty }\subset A.$
\end{center}

and

\begin{center}
$\varphi (p_{k})\overset{denote}{=}\alpha _{k}\cdot 1_{B}\in \Bbb{C}\cdot
1_{B},$ for all $k\in \Bbb{N}.$
\end{center}

\strut 

(It is possible that this sequence of projections is a finite sequence, i.e, 
$p_{\infty }=p_{N},$ for some $N$) We will call such chains of projections,
scalar-valued chains of projections.
\end{definition}

\strut

Suppose that we have a scalar-valued chain of projections $%
(p_{k})_{k=1}^{\infty }\subset A_{pro}.$ We can naturally define conditional
expectations

\strut

\begin{center}
$\varphi _{j+1}:p_{j+1}Ap_{j+1}\rightarrow p_{j}Ap_{j}$
\end{center}

by

\begin{center}
$\varphi _{j+1}\left( p_{j+1}ap_{j+1}\right) \overset{def}{=}p_{j}ap_{j},$
\end{center}

\strut \strut

for all $a\in A,$ \ $j\in \Bbb{N}.$ Notice that, by the previous lemmas,

\strut

\begin{center}
$
\begin{array}{ll}
\varphi _{j+1}\left( p_{j+1}ap_{j+1}\right) & =p_{j}\left(
p_{j+1}ap_{j+1}\right) p_{j} \\ 
& =p_{j}p_{j+1}ap_{j+1}p_{j} \\ 
& =\left( p_{j}p_{j+1}\right) a\left( p_{j}p_{j+1}\right) \\ 
& =p_{j}ap_{j},
\end{array}
$
\end{center}

\strut

for all $a\in A.$ i.e,

\begin{center}
$\varphi _{j+1}(x)=p_{j}xp_{j},$ for all $x\in p_{j+1}Ap_{j+1},$
\end{center}

for all $j\in \Bbb{N}.$

\strut

Hence we can have a tower of amalgamated NCPSpaces induced by a
scalar-valued chain of projections, $(p_{k})_{k=1}^{\infty }\subset A_{pro},$

\strut

\begin{center}
$\Bbb{C}\subset ^{\varphi _{0}}B\subset ^{\varphi _{1}}p_{1}Ap_{1}\subset
^{\varphi _{2}}p_{2}Ap_{2}\subset ^{\varphi _{3}}\cdot \cdot \cdot ,$
\end{center}

\strut

where a linear functional $\varphi _{0}:B\rightarrow \Bbb{C}$ is arbitrary
given.

\strut

Let $B$ be a unital algebra and $(A,\varphi ),$ a NCPSpace over $B.$ Suppose
that we have a tower of amalgamated NCPSpaces induced by a scalar-valued
chain of projections $(p_{k})_{k=1}^{\infty }\subset A_{pro}$ and the given
linear functional $\varphi _{0}:B\rightarrow \Bbb{C}.$ If we define a map

\strut

\begin{center}
$E_{j+1}:p_{j+1}Ap_{j+1}\rightarrow \Bbb{C}$
\end{center}

by

\begin{center}
$E_{j+1}(x)=\frac{1}{\alpha _{j+1}}\cdot \varphi _{0}\varphi _{1}...\varphi
_{j+1}(x),$ for all $x\in A_{j+1},$
\end{center}

\strut

then it is a well-defined linear functional. Recall that $\alpha _{j+1}$ is
the scalar-part of $\varphi (p_{j+1}).$

\strut \strut \strut

Now, we will introduce compressed R-transform theory in [10] ; Let $B$ be a
unital algebra and $(A,\varphi )$, a NCPSpace over $B,$ with its $B$%
-functional $\varphi :A\rightarrow B.$ Let $p\in A_{pro}$ such that

\begin{center}
$B\subset pAp\subset A$
\end{center}

and

\begin{center}
$\varphi (p)\overset{denote}{=}b_{0}\in C_{A}(B)\cap B_{inv}.$
\end{center}

\strut

Note that if we have a scalar-valued chain of projections $%
(p_{k})_{k=1}^{\infty }\subset A_{pro},$ then every $p_{j}$ satisfies

\strut

\begin{center}
$\varphi (p_{j})=\alpha _{j}\cdot 1_{B}\in C_{A}(B)\cap B_{inv},$
\end{center}

\strut

for all $j\in \Bbb{N}$ and $\alpha _{j}\in \Bbb{C}.$

\strut

We can define the conditional expectation

\begin{center}
$\varphi _{p}\overset{def}{=}b_{0}^{-1}\cdot \varphi \mid
_{pAp}:pAp\rightarrow B.$
\end{center}

\strut

So, we can define a NCPSpace over $B,$ $\left( pAp,\,\varphi _{p}\right) $
and we call it a compressed NCPSpace by $p\in A_{pro}.$ In [10], we observed
the amalgamated R-transform theory on this compressed NCPSpace over $B.$

\strut \strut \strut

\begin{definition}
Let $B$ be a unital algebra and $b\in B.$ Let $(A,\varphi )$ be a NCPSpace
over $B.$ Define a symmetric cumulants of $B$-valued random variables $%
x_{1},...,x_{s}\in (A,\varphi )$ by

\strut 

\begin{center}
$k_{n}^{symm(b)}\left( x_{i_{1}},...,x_{i_{n}}\right) =c^{(n)}\left(
x_{i_{1}}\otimes bx_{i_{2}}\otimes ...\otimes bx_{i_{n}}\right) ,$
\end{center}

\strut 

for all $(i_{1},...,i_{n})\in \{1,...,s\}^{n},$ $n\in \Bbb{N}.$ Then we can
define the symmetric R-transform of $x_{1},...,x_{s},$ by $b\in B,$ by

\strut 

\begin{center}
$R_{x_{1},...,x_{s}}^{symm(b)}(z_{1},...,z_{s})=\sum_{n=1}^{\infty }%
\underset{i_{1},...,i_{n}\in \{1,...,s\}}{\sum }%
k_{n}^{symm(b)}(x_{i_{1}},...,x_{i_{n}})\,z_{i_{1}}...z_{i_{n}},$
\end{center}

\strut 

as a $B$-formal series in $\Theta _{B}^{s}.$
\end{definition}

\strut \strut

\begin{theorem}
(See Theorem 2.1 in [10]) Let $(A,\varphi ),$ $p\in A_{pro}$ and $b_{0}\in B$
be given as before. Let $x_{1},...,x_{s}\in (A,\varphi )$ be $B$-valued
random variables ($s\in \Bbb{N}$) and assume that $\{p\}$ and $%
\{x_{1},...,x_{s}\}$ are free over $B.$ Then

\strut 

\begin{center}
$R_{px_{1}p,...,px_{s}p}^{(\varphi
_{p})\,:%
\,t}(z_{1},...,z_{s})=R_{x_{1},...,x_{s}}^{symm(b_{0})}(z_{1},...,z_{s})$
\end{center}
\end{theorem}

\strut

\begin{proof}
Fix $n\in \Bbb{N}$ and $(i_{1},...,i_{n})\in \{1,...,s\}^{n}.$ Then

\strut

$coef_{i_{1},...,i_{n}}\left( R_{px_{1}p,...,px_{s}p}^{(\varphi
_{p})\,\,:\,\,t}\right) =c_{\varphi _{p}}^{(n)}\left( px_{i_{1}}p\otimes
px_{i_{2}}p\otimes ...\otimes px_{i_{n}}p\right) $

$\ \ \ \ \ \ \ \ \ =\underset{\pi \in NC(n)}{\sum }\widehat{\varphi _{p}}%
(\pi )\left( px_{i_{1}}p\otimes ...\otimes px_{i_{n}}p\right) \mu (\pi
,1_{n})$

$\ \ \ \ \ \ \ \ \ =\underset{\pi \in NC(n)}{\sum }\widehat{%
b_{0}^{-1}\varphi }\,\,(\pi )\left( px_{i_{1}}p\otimes ...\otimes
px_{i_{n}}p\right) \mu (\pi ,1_{n})$

\strut

(2.20.1)

\strut

$\ \ \ \ \ \ \ \ \ =\underset{\pi \in NC^{\prime }(n+1)}{\sum }b_{0}^{-1}\,\,%
\widehat{\varphi _{p}}\,\,(\pi )\left( p\otimes x_{i_{1}}p\otimes
x_{i_{2}}p\otimes ...\otimes x_{i_{n}}p\right) \mu (\pi ,1_{n+1}),$

\strut

where $NC^{\prime }(n+1)=\{\theta \in NC(n+1):(1)\in \theta \}$ which is
lattice-isomorphic to $NC(n).$ By the $B$-freeness of $\{p\}$ and $%
\{x_{1},...,x_{s}\},$ we have that (2.20.1) is equivalent to (2.20.2) (For
this squivalence, see the proof of Theorem 2.1 in [10]) ;

\strut

\strut (2.20.2)

\strut

$\underset{\pi \in NC^{\prime }(n+1)}{\sum }b_{0}^{-1}\cdot \widehat{%
c_{\varphi }}\left( \pi \cup Kr(\pi )\right) \left( 1_{B}\otimes p\otimes
x_{i_{1}}\otimes p\otimes ...\otimes x_{i_{n}}\otimes p\right) $

$\ \ \ \ \ =\underset{\pi \in NC^{\prime }(n+1)}{\sum }b_{0}^{-1}\cdot 
\widehat{c_{\varphi }}\left( \pi \right) (1_{B}\otimes x_{i_{1}}\otimes
...\otimes x_{i_{n}})\cdot \widehat{c_{\varphi _{p}}}(Kr(\pi ))\left(
p\otimes ...\otimes p\right) $

\strut

since $p$ is $B$-central, in the sense of [9]

\strut

\ $\ \ \ =\underset{\pi \in NC^{\prime }(n+1)}{\sum }b_{0}^{-1}\cdot 
\widehat{c_{\varphi }}(\pi )\left( 1_{B}\otimes x_{i_{1}}\otimes ...\otimes
x_{i_{n}}\right) b_{0}^{\left| Kr(\pi )\right| }$

$\ \ \ \ \ \ =\underset{\pi \in NC^{\prime }(n+1)}{\sum }b_{0}^{-1}\cdot 
\widehat{c_{\varphi }}(\pi )\left( 1_{B}\otimes x_{i_{1}}\otimes ...\otimes
x_{i_{n}}\right) b_{0}^{(n+1)-\left| \pi \right| },$

\strut

since $\left| \pi \right| +\left| Kr(\pi )\right| =n+1,$ for all $\pi \in
NC(n)$. (Notice that $NC^{\prime }(n+1)=NC(n).$) Therefore,

\strut

\begin{center}
$k_{n}^{(\varphi _{p})\,\,:\,\,t}\left( px_{i_{1}}p,...,px_{i_{n}}p\right)
=k_{n}^{(\varphi )\,:\,symm(b_{0})}\left( x_{i_{1}},...,x_{i_{n}}\right) .$
\end{center}
\end{proof}

\strut

Remark that $p\in A_{pro}$ commutes with $B.$ Since $b=pbp,$ we have that

$\strut \strut $

\begin{center}
$pb=p(pbp)=p^{2}bp=pbp=pbp^{2}=(pbp)p=bp,$
\end{center}

\strut

for all $b\in B.$ i.e, $p\in (A,\varphi )$ is $B$-central, in the sense of
[9].

\strut \strut

\strut Let $(p_{k})_{k=1}^{\infty }\subset A_{pro}$ be a scalar-valued chain
of projections in $(A,\varphi ).$ Then each $p_{j}\in A_{pro}$ satisfies $%
\varphi (p_{j})=\alpha _{j}\cdot 1_{B},$ for some $\alpha _{j}\in \Bbb{C}.$
Therefore, each $\varphi (p_{j})\in C_{A_{j}}(A_{j-1})\cap (A_{j-1})_{inv},$
for all $j\in \Bbb{N},$ with $A_{0}=B.$

\strut

\begin{corollary}
Let $B$ be a unital algebra and $(A,\varphi ),$ a NCPSpace over $B.$ Let $%
x_{1},...,x_{s}\in (A,\varphi )$ be $B$-valued random variables ($s\in \Bbb{N%
}$). Suppose that we have $p\in A_{pro}$ such that $B\subset pAp\subset A$
and $\varphi (p)\overset{denote}{=}b_{0}\in C_{A}(B)\cap B_{inv}.$ If $\{p\}$
and $\{x_{1},...,x_{s}\}$ are free over $B,$ then

\strut 

\begin{center}
$
\begin{array}{ll}
k_{n}^{(\varphi _{p})}\left( px_{i_{1}}p,...,px_{i_{n}}p\right)  & 
:=c_{\varphi _{p}}^{(n)}\left( px_{i_{1}}p\otimes
b_{i_{2}}(px_{i_{2}}p)\otimes ...\otimes b_{i_{n}}(px_{i_{n}}p)\right)  \\ 
& =c_{\varphi }^{(n)}\left( x_{i_{1}}\otimes b_{0}b_{i_{2}}x_{i_{2}}\otimes
...\otimes b_{0}b_{i_{n}}x_{i_{n}}\right) ,
\end{array}
$
\end{center}

\strut 

for all $(i_{1},...,i_{n})\in \{1,...,s\}^{n},$ $n\in \Bbb{N},$ where $%
b_{i_{2}},...,b_{i_{n}}\in B$ are arbitrary.
\end{corollary}

\strut

\begin{proof}
The proof is similar to the proof of the previous theorem. Fix $n\in \Bbb{N}$
and $(i_{1},...,i_{n})\in \{1,...,s\}^{n}.$ Observe that

\strut

\begin{center}
$k_{n}^{(\varphi _{p})}\left( px_{i_{1}}p,...,px_{i_{n}}p\right) =c_{\varphi
_{p}}^{(n)}\left( px_{i_{1}}p\otimes b_{i_{2}}px_{i_{2}}p\otimes ...\otimes
pb_{i_{n}}x_{i_{n}}p\right) $
\end{center}

\strut

where $b_{i_{2}},...,b_{i_{n}}\in B$ are arbitrary. Since $%
pb_{i_{j}}p=b_{i_{j}},$ we have that

\strut

\begin{center}
$pb_{i_{j}}=p\left( pb_{i_{j}}p\right)
=p^{2}b_{i_{j}}p=pb_{i_{j}}p=pb_{i_{j}}p^{2}=(pb_{i_{j}}p)p=b_{i_{j}}p,$
\end{center}

\strut

for all $j=1,...,n.$ So, we have that

\strut

\begin{center}
$k_{n}^{(\varphi _{p})}\left( px_{i_{1}}p,...,px_{i_{n}}p\right) =c_{\varphi
_{p}}^{(n)}\left( px_{i_{1}}p\otimes px_{i_{2}}^{\prime }p\otimes ...\otimes
px_{i_{n}}^{\prime }p\right) ,$
\end{center}

\strut

where $x_{i_{j}}^{\prime }=b_{i_{j}}x_{i_{j}},$ for all $\ j=1,...,n.$
Similar to the proof of the previous theorem, we have that

\strut

\begin{center}
$c_{\varphi _{p}}^{(n)}\left( px_{i_{1}}p\otimes px_{i_{2}}^{\prime
}p\otimes ...\otimes px_{i_{n}}^{\prime }p\right) =c_{\varphi }^{(n)}\left(
x_{i_{1}}\otimes b_{0}x_{i_{2}}^{\prime }\otimes ...\otimes
b_{0}x_{i_{n}}^{\prime }\right) ,$
\end{center}

\strut

where $b_{0}=\varphi (p)\in C_{A}(B)\cap B_{inv}.$ Therefore,

\strut

\begin{center}
$k_{n}^{(\varphi _{p})}\left( x_{i_{1}},...,x_{i_{n}}\right) =c_{\varphi
}^{(n)}\left( x_{i_{1}}\otimes b_{0}x_{i_{2}}^{\prime }\otimes ...\otimes
b_{0}x_{i_{n}}^{\prime }\right) .$
\end{center}
\end{proof}

\strut

By the previous corollary, we can get the following theorem ;

\strut

\begin{theorem}
(Also see Theorem 2.5 and Corollary 2.7 in [10]) Let $B$ be a unital algebra
and let $(A,\varphi )$ be a NCPSpace over $B.$ Let $p\in A_{pro}$ satisfy $%
\varphi (p)\in C_{A}(B)\cap B_{inv}$ and suppose that $X=\{x_{1},...,x_{s}\}$
and $Y=\{y_{1},...,y_{s}\}$ are two subsets of $B$-valued random variables
in $(A,\varphi ).$ Assume that $\{p\}$ and $X\cup Y$ are free over $B,$ in $%
(A,\varphi ).$ If $X$ and $Y$ are free over $B,$ in $(A,\varphi ),$ then $pXp
$ and $pYp$ are free over $B,$ in $\left( pAp,\varphi _{p}\right) .$
\end{theorem}

\strut

\begin{proof}
Suppose that $\{p\}$ and $X\cup Y$ are free over $B.$ Then, by the previous
corollary, we can get that

\strut

$k_{n}^{(\varphi _{p})}\left( a_{i_{1}},...,a_{i_{n}}\right) =c_{\varphi
_{p}}^{(n)}\left( pa_{i_{1}}p\otimes b_{i_{2}}pa_{i_{2}}p\otimes ...\otimes
b_{i_{n}}pa_{i_{n}}p\right) $

\strut

$\ \ \ \ \ \ \ \ \ \ \ \ \ \ \ \ \ =c_{\varphi _{p}}^{(n)}\left(
pa_{i_{1}}p\otimes pb_{i_{2}}a_{i_{2}}p\otimes ...\otimes
pb_{i_{n}}a_{i_{n}}p\right) $

\strut

(2.22.1)\strut

\strut $\ \ \ \ \ \ \ \ \ \ \ \ \ \ \ \ =c_{\varphi }^{(n)}\left(
a_{i_{1}}\otimes b_{0}(b_{i_{2}}a_{i_{2}})\otimes ...\otimes
b_{0}(b_{i_{n}}a_{i_{n}})\right) $

\strut

where $b_{0}=\varphi (p)\in C_{A}(B)\cap B_{inv},$ $b_{i_{2}},...,b_{i_{n}}%
\in B$ are arbitrary and where $a_{i_{1}},...,a_{i_{n}}\in X\cup Y.$ We can
rewrite (2.22.1) as

\strut

(2.22.2) $\ \ \ \ \ \ \ \ \ \ \ \ \ \ \ \ \ \ \ \ \ c_{\varphi }^{(n)}\left(
a_{i_{1}}\otimes b_{i_{2}}^{\prime }a_{i_{2}}\otimes ...\otimes
b_{i_{n}}^{\prime }a_{i_{n}}\right) ,$

\strut

where $b_{i_{2}}^{\prime },...,b_{i_{n}}^{\prime }\in B$ are arbitrary such
that $b_{i_{j}}^{\prime }=b_{0}b_{i_{j}},$ for all $j=1,...,n.$ Since $X$
and $Y$ are free over $B,$ in $(A,\varphi ),$ we have that

\strut

(2.22.2) \ \ \ $=\left\{ 
\begin{array}{lll}
c_{\varphi }^{(n)}\left( x_{i_{1}}\otimes b_{i_{2}}^{\prime
}x_{i_{2}}\otimes ...\otimes b_{i_{n}}^{\prime }x_{i_{n}}\right) &  & \text{%
or} \\ 
&  &  \\ 
c_{\varphi }^{(n)}\left( y_{i_{1}}\otimes b_{i_{2}}^{\prime
}y_{i_{2}}\otimes ...\otimes b_{i_{n}}^{\prime }y_{i_{n}}\right) . &  & 
\end{array}
\right. $

\strut

$\ \ \ \ \ \ \ \ \ \ \ \ \ \ =\left\{ 
\begin{array}{lll}
k_{n}^{(\varphi _{p})}\left( px_{i_{1}}p,...,px_{i_{n}}p\right) &  & \text{or%
} \\ 
&  &  \\ 
k_{n}^{(\varphi _{p})}\left( py_{i_{1}}p,...,py_{i_{n}}p\right) . &  & 
\end{array}
\right. $

\strut

Therefore, $pXp$ and $pYp$ are also free over $B,$ in $\left( pAp,\,\varphi
_{p}\right) .$
\end{proof}

\strut

\begin{corollary}
Let $B$ be a unital algebra and $(A,\varphi ),$ a NCPSpace over $B.$ Let $%
X=\{x_{1},...,x_{s}\}$ and $Y=\{y_{1},...,y_{s}\}$ be subsets of $B$-valued
random variables in $(A,\varphi ).$ Assume that $X$ and $Y$ are free over $B,
$ in $(A,\varphi ).$ Now, let $p\in A_{pro}$ such that $\{p\}$ and $X\cup Y$
are free over $B,$ in $(A,\varphi ).$ By $b_{0}\in C_{A}(B)\cap B_{inv},$ we
will denote $\varphi (p)\in C_{A}(B)\cap B_{inv}.$ Then

\strut 

(1) $R_{px_{1}p,...,px_{s}p}^{(\varphi
_{p})\,\,:\,\,t}(z_{1},...,z_{s})=R_{x_{1},...,x_{s}}^{(\varphi
)\,\,:\,\,symm(b_{0})}(z_{1},...,z_{s})$ \ and

\strut 

\ \ \ \ \ $R_{px_{1}p,...,px_{s}p}^{(\varphi _{p})}(z_{1},...,z_{s})$

$\ \ \ \ \ \ \ \ \ \ =\sum_{n=1}^{\infty }\underset{i_{1},...,i_{n}\in
\{1,...,s\}}{\sum }c_{\varphi }^{(n)}\left( x_{i_{1}}\otimes
b_{0}b_{i_{2}}x_{i_{2}}\otimes ...\otimes b_{0}b_{i_{n}}x_{i_{n}}\right)
\,z_{i_{1}}...z_{i_{n}},$

\strut 

where $b_{i_{2}},...,b_{i_{n}}\in B$ are arbitrary.

\strut 

(2) $R_{px_{1}p,...,px_{s}p,py_{1}p,...,py_{s}p}^{(\varphi
_{p})}(z_{1},...,z_{2s})$

\begin{center}
$=R_{px_{1}p,...,px_{s}p}^{(\varphi
_{p})}(z_{1},...,z_{s})+R_{py_{1}p,...,py_{s}p}^{(\varphi
_{p})}(z_{s+1},...,z_{2s}).$
\end{center}

\strut 

(3) $R_{px_{1}p+py_{1}p,...,px_{s}p+py_{1}p}^{(\varphi
_{p})}(z_{1},...,z_{s})$

\begin{center}
$=\left( R_{px_{1}p,...,px_{s}p}^{(\varphi
_{p})}+R_{py_{1}p,...,py_{s}p}^{(\varphi _{p})}\right) (z_{1},...,z_{s}).$
\end{center}

\strut 

(4) $R_{(px_{1}p)(py_{1}p),...(px_{s}p)(py_{s}p)}^{(\varphi
_{p})}(z_{1},...,z_{s})$

$\ \ \ \ \ \ \ \ \ \ \ \ \ \ \ \ \ \ \ \ =\left(
R_{px_{1}p,...,px_{s}p}^{(\varphi _{p})}\,\,\frame{*}_{B}\,%
\,R_{py_{1}p,...,py_{s}p}^{(\varphi _{p})\,\,:\,\,t}\right) (z_{1},...,z_{s})
$

$\ \ \ \ \ \ \ \ \ \ \ \ \ \ \ \ \ \ \ \ =\left(
R_{px_{1}p,...,px_{s}p}^{(\varphi _{p})}\,\,\frame{*}_{B}\,%
\,R_{y_{1},...,y_{s}}^{(\varphi )\,\,:\,\,t}\right) (z_{1},...,z_{s})$
\end{corollary}

\strut

\begin{proof}
(1), (2) and (3) are proved in the previous theorems and corollary. Since $%
pXp$ and $pYp$ are free over $B,$ in $\left( pAp,\varphi _{p}\right) ,$ we
have that

\strut

\begin{center}
$R_{(px_{1}p)(py_{1}p),...,(px_{s}p)(py_{s}p)}^{(\varphi
_{p})}=R_{px_{1}p,...,px_{s}p}^{(\varphi _{p})}\,\,\,\frame{*}%
_{B}\,\,\,R_{py_{1}p,...,py_{s}p}^{(\varphi _{p})}.$
\end{center}

\strut

Also, we can get the following result ; now, fix $n\in \Bbb{N}$ and $%
(i_{1},...,i_{n})\in \{1,...,s\}^{n}.$

\strut

$coef_{i_{1},...,i_{n}}\left(
R_{(px_{1}p)(py_{1}p),...,(px_{s}p)(py_{1}p)}^{(\varphi _{p})}\right) $

$\ \ \ =c_{\varphi _{p}}^{(n)}\left( px_{i_{1}}y_{i_{1}}p\otimes
b_{i_{2}}px_{i_{2}}y_{i_{2}}p\otimes ...\otimes
b_{i_{n}}px_{i_{n}}y_{i_{n}}p\right) ,$

\strut

since $%
(px_{i_{j}}p)(py_{i_{j}}p)=px_{i_{j}}p^{2}y_{i_{j}}p=px_{i_{j}}py_{i_{j}}p=px_{i_{j}}y_{i_{j}}p, 
$ for all $j=1,...,n,$ where $b_{i_{2}},...,b_{i_{n}}\in B$ are arbitrary

\strut

$\ \ \ =c_{\varphi }^{(n)}\left( x_{i_{1}}y_{i_{1}}\otimes
b_{0}b_{i_{2}}x_{i_{2}}y_{i_{2}}\otimes ...\otimes
b_{0}b_{i_{n}}x_{i_{n}}y_{i_{n}}\right) $

$\ \ \ =c_{\varphi }^{(n)}\left( x_{i_{1}}y_{i_{1}}\otimes b_{i_{2}}^{\prime
}x_{i_{2}}y_{i_{2}}\otimes ...\otimes b_{i_{n}}^{\prime
}x_{i_{n}}y_{i_{n}}\right) $

\strut

where $b_{i_{j}}^{\prime }=b_{0}b_{i_{j}}\in B,$ for all $j=1,...,n$

\strut

$\ \ \ =\underset{\pi \in NC(n)}{\sum }\widehat{c_{\varphi }}\,(\pi \cup
Kr(\pi ))\left( x_{i_{1}}\otimes y_{i_{1}}\otimes b_{i_{2}}^{\prime
}x_{i_{2}}\otimes y_{i_{2}}\otimes ...\otimes b_{i_{n}}^{\prime
}x_{i_{n}}\otimes y_{i_{n}}\right) $

\strut

by the freeness of $\{x_{1},...,x_{s}\}$ and $\{y_{1},...,y_{s}\}.$ So,

\strut

\begin{center}
$coef_{i_{1},..,i_{n}}\left(
R_{(px_{1}p)(py_{1}p),...,(px_{s}p)(py_{s}p)}^{(\varphi _{p})}\right)
=coef_{i_{1},...,i_{n}}\left( R_{px_{1}p,...,px_{s}p}^{(\varphi _{p})}\,\,%
\frame{*}_{B}\,\,R_{y_{1},...,y_{s}}^{(\varphi )\,\,:\,\,t}\right) .$
\end{center}
\end{proof}

\strut

\begin{quote}
\frame{\textbf{Notations}} Let $B$ be a unital algebra and $(A,\varphi ),$ a
NCPSpace over $B.$ Now consider the following tower of amalgamated NCPSpaces,

\strut
\end{quote}

\begin{center}
$\Bbb{C}\subset ^{\varphi _{0}}B\subset ^{\varphi _{1}}p_{1}Ap_{1}\subset
^{\varphi _{2}}p_{2}Ap_{2}\subset ^{\varphi _{3}}\cdot \cdot \cdot ,$
\end{center}

\begin{quote}
\strut

induced by a scalar-valued chain of projections in $A,$ $(p_{k})_{k=1}^{%
\infty }\subset A_{pro},$ satisfying $\varphi (p_{k})=\alpha _{k}\cdot
1_{B}\in C_{A}(B)\cap B_{inv},$ $\forall k,$ and the given linear functional 
$\varphi _{0}:B\rightarrow \Bbb{C}.$ As we defined before,

\strut
\end{quote}

\begin{center}
$\varphi _{j+1}(p_{j+1}ap_{j+1})=p_{j}ap_{j},$ for each $j\in \Bbb{N},$ for
all $a\in A.$
\end{center}

\begin{quote}
\strut \strut \strut
\end{quote}

As in Section 2.1, we can apply the R-transform theory on a tower of
amalgamated NCPSpaces to this compressed amalgamated NCPSpaces induced by a
scalar-valued chain of projections $(p_{k})_{k=1}^{\infty }\subset A_{pro}$
and the given linear functional $\varphi _{0}.$\strut

\strut \strut

\begin{theorem}
Let $B$ be a unital algebra and $(A,\varphi ),$ a NCPSpace over $B.$ Let $%
(p_{k})_{k=1}^{\infty }\subset A_{pro}$ be a scalar-valued chain of
projections. Consider the tower of amalgamated NCPSpaces induced by the
chain of projections $(p_{k})_{k=1}^{\infty }$ and a linear functional $%
\varphi _{0}:B\rightarrow \Bbb{C}.$ Then, for any $j\in \Bbb{N},$ if two
subsets of $B$-valued random variables $X$ and $Y$ are free over $B,$ in $%
(A,\varphi ),$ and if $p_{\infty }$ and $X\cup Y$ are free over $B,$ in $%
(A,\varphi ),$ then $p_{j}Xp_{j}$ and $p_{j}Yp_{j}$ are free over $B,$ in $%
\left( p_{j}Ap_{j},\,\frac{1}{\alpha _{1}}\cdot E_{1,\,j}\right) .$ $\square 
$
\end{theorem}

\strut

Notice that

\strut

\begin{center}
$\frac{1}{\alpha _{1}}=\frac{1}{E_{1,\,j}(p_{1}...p_{j})}=\frac{1}{%
E_{1,\,j}(p_{1})}=\frac{1}{\varphi _{1}...\varphi _{j}(p_{1})}=\frac{1}{%
\varphi _{1}(p_{1})}=\frac{1}{\varphi (p_{1})}.$
\end{center}

\strut \strut \strut

As an application of of the above theorem, we can get the following
R-transform calculus ;

\strut \strut

\begin{theorem}
Let $B$ be a unital algebra and $(A,\varphi ),$ a NCPSpace over $B.$ Let $%
(p_{k})_{k=1}^{\infty }\subset A_{pro}$ be a scalar-valued chain of
projections. Suppose that we have a tower of amalgamated NCPSpaces induced
by the chain of projections and a linear functional $\varphi
_{0}:B\rightarrow \Bbb{C}.$ Let $X=\{x_{1},...,x_{s}\}$ and $%
Y=\{y_{1},...,y_{s}\}$ be two $B$-free subsets in $(A,\varphi ).$ If $%
p_{\infty }$ and $X\cup Y$ are free over $B,$ in $(A,\varphi ),$ then, for
any fixed $\ j\in \Bbb{N},$ we have that

\strut 

(1) $%
R_{p_{j}x_{1}p_{j},...,p_{j}x_{s}p_{j},p_{j}y_{1}p_{j},...,p_{j}y_{s}p_{j}}^{(%
\frac{1}{\alpha _{1}}E_{1,\,j})}(z_{1},...,z_{2s})$

\begin{center}
$=R_{p_{j}x_{1}p_{j},...,p_{j}x_{s}p_{j}}^{(\frac{1}{\alpha _{1}}%
E_{1,\,\,j})}$\strut $%
(z_{1},...,z_{s})+R_{p_{j}yp_{j},...,p_{j}y_{s}p_{j}}^{(\frac{1}{\alpha _{1}}%
E_{1,\,j})}(z_{s+1},...,z_{2s}).$
\end{center}

\strut 

(2) $%
R_{p_{j}x_{1}p_{j}+p_{j}y_{1}p_{j},...,p_{j}x_{s}p_{j}+p_{j}y_{s}p_{j}}^{(%
\frac{1}{\alpha _{1}}E_{1,\,j})}(z_{1},...,z_{2s})$

\begin{center}
$=R_{p_{j}x_{1}p_{j},...,p_{j}x_{s}p_{j}}^{(\frac{1}{\alpha _{1}}%
E_{1,\,\,j})}$\strut $%
(z_{1},...,z_{s})+R_{p_{j}yp_{j},...,p_{j}y_{s}p_{j}}^{(\frac{1}{\alpha _{1}}%
E_{1,\,j})}(z_{1},...,z_{s}).$
\end{center}

\strut \strut 

(3) $%
R_{(p_{j}x_{1}p_{j})(p_{j}y_{1}p_{j}),...,(p_{j}x_{s}p_{j})(p_{j}y_{s}p_{j})}^{(%
\frac{1}{\alpha _{1}}E_{1,\,j})}(z_{1},...,z_{s})$

\begin{center}
$=\left( R_{p_{j}x_{1}p_{j},...,p_{j}x_{s}p_{j}}^{(\frac{1}{\alpha _{1}}%
E_{1,\,j})}\,\,\,\frame{*}_{B}\,\,%
\,R_{p_{j}y_{1}p_{j},...,p_{j}y_{s}p_{j}}^{(\frac{1}{\alpha _{1}}%
E_{1,j})\,\,:\,\,t}\right) (z_{1},...,z_{s})$

$=\left( R_{p_{j}x_{1}p_{j},...,p_{j}x_{s}p_{j}}^{(\frac{1}{\alpha _{1}}%
E_{1,\,j})}\,\,\,\,\,\frame{*}_{B}\,\,\,\,\,\,\,R_{y_{1}\,\,,\,\,\,.\,\,\,.%
\,\,\,.\,\,\,,y_{s}}^{(\varphi )\,\,:\,\,t}\right) \,(z_{1},...,z_{s}).$
\end{center}

$\square $
\end{theorem}

\strut

By using the same idea of the previous two theorems, we have that ;\strut

\strut

\begin{theorem}
Let $B$ be a unital algebra and $(A,\varphi ),$ a NCPSpace over $B.$ Let $%
(p_{k})_{k=1}^{\infty }\subset A_{pro}$ be a scalar-valued chain of
projections. Consider the tower of amalgamated NCPSpaces induced by the
chain of projections $(p_{k})_{k=1}^{\infty }$ and a linear functional $%
\varphi _{0}:B\rightarrow \Bbb{C}.$ Then, for any $k<j\in \Bbb{N},$ if two
subsets of $B$-valued random variables $X$ and $Y$ are free over $A_{k},$ in 
$(A,\varphi ),$ then $p_{j}Xp_{j}$ and $p_{j}Yp_{j}$ are free over $B,$ in $%
\left( p_{j}Ap_{j},\,\frac{1}{\alpha _{k+1}}\cdot E_{k+1,\,j}\right) .$ $%
\square $
\end{theorem}

\strut \strut \strut

\begin{theorem}
Let $B$ be a unital algebra and $(A,\varphi ),$ a NCPSpace over $B.$ Let $%
(p_{k})_{k=1}^{\infty }\subset A_{pro}$ be a scalar-valued chain of
projections. Suppose that we have a tower of amalgamated NCPSpaces induced
by the chain of projections and a linear functional $\varphi
_{0}:B\rightarrow \Bbb{C}.$ Let $X=\{x_{1},...,x_{s}\}$ and $%
Y=\{y_{1},...,y_{s}\}$ be two $A_{k}$-free subsets in $(A,\varphi ).$ Then,
for any fixed $k<\ j\in \Bbb{N},$ we have that

\strut 

(1) $%
R_{p_{j}x_{1}p_{j},...,p_{j}x_{s}p_{j},p_{j}y_{1}p_{j},...,p_{j}y_{s}p_{j}}^{(%
\frac{1}{\alpha _{k+1}}E_{k+1,\,j})}(z_{1},...,z_{2s})$

\begin{center}
$=R_{p_{j}x_{1}p_{j},...,p_{j}x_{s}p_{j}}^{(\frac{1}{\alpha _{k+1}}%
E_{k+1,\,\,j})}$\strut $%
(z_{1},...,z_{s})+R_{p_{j}yp_{j},...,p_{j}y_{s}p_{j}}^{(\frac{1}{\alpha
_{k+1}}E_{k+1,\,j})}(z_{s+1},...,z_{2s}).$
\end{center}

\strut 

(2) $%
R_{p_{j}x_{1}p_{j}+p_{j}y_{1}p_{j},...,p_{j}x_{s}p_{j}+p_{j}y_{s}p_{j}}^{(%
\frac{1}{\alpha _{k+1}}E_{k+1,\,j})}(z_{1},...,z_{2s})$

\begin{center}
$=R_{p_{j}x_{1}p_{j},...,p_{j}x_{s}p_{j}}^{(\frac{1}{\alpha _{k+1}}%
E_{k+1,\,\,j})}$\strut $%
(z_{1},...,z_{s})+R_{p_{j}yp_{j},...,p_{j}y_{s}p_{j}}^{(\frac{1}{\alpha
_{k+1}}E_{k+1,\,j})}(z_{1},...,z_{s}).$
\end{center}

\strut \strut 

(3) $%
R_{(p_{j}x_{1}p_{j})(p_{j}y_{1}p_{j}),...,(p_{j}x_{s}p_{j})(p_{j}y_{s}p_{j})}^{(%
\frac{1}{\alpha _{k+1}}E_{k+1,\,j})}(z_{1},...,z_{s})$

\begin{center}
$=\left( R_{p_{j}x_{1}p_{j},...,p_{j}x_{s}p_{j}}^{(\frac{1}{\alpha _{k+1}}%
E_{k+1,\,j})}\,\,\,\frame{*}_{B}\,\,%
\,R_{p_{j}y_{1}p_{j},...,p_{j}y_{s}p_{j}}^{(\frac{1}{\alpha _{k+1}}%
E_{k+1,j})\,\,:\,\,t}\right) (z_{1},...,z_{s})$

$=\left( R_{p_{j}x_{1}p_{j},...,p_{j}x_{s}p_{j}}^{(\frac{1}{\alpha _{k+1}}%
E_{k+1,\,j})}\,\,\,\,\,\frame{*}_{B}\,\,\,\,\,\,\,R_{y_{1}\,\,,\,\,\,.\,\,%
\,.\,\,\,.\,\,\,,y_{s}}^{(\varphi )\,\,:\,\,t}\right) \,(z_{1},...,z_{s}).$
\end{center}

$\square $
\end{theorem}

\strut

\strut

\strut

\strut

\section{Compatibility of towers of Amalgamated Noncommutative Probability
Spaces}

\strut

\strut

\strut

Throughout this chapter, we will consider the tower of amalgamated NCPSpaces,

\strut

\begin{center}
$\Bbb{C}\subset ^{\varphi _{0}}B\subset ^{\varphi _{1}}A_{1}\subset
^{\varphi _{2}}A_{2}\subset ^{\varphi _{3}}A_{3}\subset \cdot \cdot \cdot .$
\end{center}

\strut

In this chapter, we will consider the compatibility of amalgamated NCPSpaces.

\strut

\begin{definition}
Let $(A,\varphi )$ be a NCPSpace over $B$ and let $(A,\varphi _{0})$ be a
NCPSpace (over $\Bbb{C}$). We say that $(A,\varphi )$ and $(A,\varphi _{0})$
are compatible if

\strut 

\begin{center}
$\varphi _{0}(x)=\varphi _{0}\left( \varphi (x)\right) ,$ \ for all $x\in A.$
\end{center}

\strut 

Now let $D$ be a unital subalgebra of $A$ such that $1_{D}=1_{A}(=1_{B})$
and assume that there is a conditional expectation $\varphi ^{\prime
}:A\rightarrow D$ and hence we have a NCPSpace over $D,$ $(A,\varphi
^{\prime }).$ (i.e, $B\subset D\subset A.$) We say that two amalgamated
NCPSpaces $(A,\varphi ^{\prime }:A\rightarrow D)$ and $(A,\varphi
:A\rightarrow B)$ are compatible if

\strut 

\begin{center}
$\varphi (x)=\varphi \left( \varphi ^{\prime }(x)\right) ,$ \ for all $x\in
A.$
\end{center}
\end{definition}

\strut

\begin{proposition}
Let $B$ be a unital algebra and let

\strut 

\begin{center}
$\Bbb{C}\subset ^{\varphi _{0}}B\subset ^{\varphi _{1}}A_{1}\subset
^{\varphi _{2}}A_{2}\subset ^{\varphi _{3}}\cdot \cdot \cdot $
\end{center}

\strut 

be a tower of amalgamated NCPSpaces. Then, for any fixed $\ j\in \Bbb{N},$ $%
\left( A_{j+1},\varphi _{j+1}\right) $ and $\left( A_{j+1},\varphi
_{j}\varphi _{j+1}\right) $ are compatible.
\end{proposition}

\strut

\begin{proof}
For any $x\in A_{j+1},$ we have that

\strut

\begin{center}
$\varphi _{j}\varphi _{j+1}(x)=\varphi _{j}\varphi _{j+1}^{2}(x)=\varphi
_{j}\varphi _{j+1}\left( \varphi _{j+1}(x)\right) ,$
\end{center}

\strut

since $\varphi _{j+1}^{2}(x)=\varphi _{j+1}\left( \varphi _{j+1}(x)\right)
=\varphi _{j+1}(x).$
\end{proof}

\strut

\strut In general, we can get the following facts ;

\strut

\begin{theorem}
Let $B$ be a unital algebra and let

\strut 

\begin{center}
$\Bbb{C}\subset ^{\varphi _{0}}B\subset ^{\varphi _{1}}A_{1}\subset
^{\varphi _{2}}A_{2}\subset ^{\varphi _{3}}A_{3}\subset ^{\varphi _{4}}\cdot
\cdot \cdot $
\end{center}

\strut 

be a tower of amalgamated NCPSpaces. For any fixed $k<i<j$ in $\Bbb{N},$ $%
\left( A_{j},\,\,E_{i,\,j}\right) $ and $\left( A_{j},\,\,E_{k,\,j}\right) $
are compatible.
\end{theorem}

\strut \strut

\begin{proof}
For any $x\in A_{j},$ we have that

\strut

$\ \ E_{k,\,j}(x)=\varphi _{k}...\varphi _{i}...\varphi _{j}(x)=\varphi
_{k}...\varphi _{i-1}\left( \varphi _{i}...\varphi _{j}(x)\right) $

$\ \ \ \ \ \ \ \ \ \ \ \ \ =E_{k,\,\,i-1}\left( E_{i,\,j}(x)\right)
=E_{k,\,i-1}\left( E_{i,\,j}^{2}(x)\right) $

\strut

since $E_{i,\,j}:A_{j}\rightarrow A_{i-1}$ is a conditional expectation

\strut

$\ \ \ \ \ \ \ \ \ \ \ \ \ =E_{k,\,i-1}E_{i,\,j}\left( E_{i,\,j}(x)\right)
=E_{k,\,j}\left( E_{i,\,j}(x)\right) .$

\strut

Therefore,

\strut

\begin{center}
$E_{k,\,j}(x)=E_{k,\,j}\left( E_{i,\,j}(x)\right) \in A_{k-1},$ for all $%
x\in A_{j}.$
\end{center}
\end{proof}

\strut

\begin{theorem}
Let $B$ be a unital algebra and suppose that we have a tower of amalgamated
NCPSpaces as before. Then, for any fixed $k<j$ in $\Bbb{N},$ $\left(
A_{j},\,E_{j}\right) $ and $\left( A_{j},\,E_{k,\,j}\right) $ are compatible.
\end{theorem}

\strut

\begin{proof}
For any $x\in A_{j},$ we have that

\strut

\begin{center}
$
\begin{array}{ll}
E_{j}(x) & =\varphi _{0}\varphi _{1}..\varphi _{k}...\varphi
_{j}(x)=E_{k-1}\left( E_{k,\,j}(x)\right) \\ 
& =E_{k-1}\left( E_{k,\,j}^{2}(x)\right) =E_{k-1}E_{k,\,j}\left(
E_{k,\,j}(x)\right) \\ 
& =E_{j}\left( E_{k,\,j}(x)\right) \in \Bbb{C}.
\end{array}
$
\end{center}

\strut
\end{proof}

\strut

The above two theorems shows us that we can freely use the compatibility on
the tower.

\strut

\begin{lemma}
(See [18] and [12]) Let $B$ be a unital algebra and $A,$ an algebra over $B.$
Suppose a NCPSpace $(A,\varphi _{0})$ and a NCPSpace over $B,$ $(A,\varphi )$
are compatible. Let $x_{1},...,x_{s}\in A$ be operators such that $%
\{x_{1},...,x_{s}\}$ and $B$ are free in $(A,\varphi _{0}).$ Then

\strut 

\begin{center}
$
\begin{array}{ll}
k_{n}^{(\varphi )}\left( x_{i_{1}},...,x_{i_{n}}\right)  & =c^{(n)}\left(
x_{i_{1}}\otimes b_{i_{2}}x_{i_{2}}\otimes ...\otimes
b_{i_{n}}x_{i_{n}}\right)  \\ 
& =\left( \varphi _{0}(b_{i_{2}})...\varphi _{0}(b_{i_{n}})\right) \cdot
k_{n}^{(\varphi _{0})\,}\left( x_{i_{1}},...,x_{i_{n}}\right) \cdot 1_{B},
\end{array}
$
\end{center}

\strut 

in $B,$ where $k_{n}^{(\varphi _{0})}(...)$ is the scalar-valued cumulants
of $x_{1},...,x_{s},$ in the sense of Speicher and Nica, and where $%
b_{i_{2}},...,b_{i_{n}}\in B$ are arbitrary, for all $(i_{1},...,i_{n})\in
\{1,...,s\}^{n},$ $n\in \Bbb{N}.$ $\square $
\end{lemma}

\strut

In the above lemma, $k_{n}^{(\varphi _{0})}(...)$ is our $k_{n}^{(\varphi
_{0})\,:\,t}(...).$ (See Chapter 1 and [6], [7])\strut 

\strut

\begin{proposition}
Let $B$ be a unital algebra and let

\strut 

\begin{center}
$\Bbb{C}\subset ^{\varphi _{0}}B\subset ^{\varphi _{1}}A_{1}\subset
^{\varphi _{2}}A_{2}\subset ^{\varphi _{3}}\cdot \cdot \cdot $
\end{center}

\strut 

be a tower of amalgamated NCPSpaces. Fix $k<\ j\in \Bbb{N}.$ Suppose that $%
X=\{x_{1},...,x_{s}\},$ $Y=\{y_{1},...,y_{s}\}\subset A_{j+1}$ are two
subsets of operators and assume that $X\cup Y$ is free from $A_{k},$ in $%
\left( A_{j+1},\,E_{j+1}\right) .$ If $X$ and $Y$ are free in $\left(
A_{j+1},\,E_{j+1}\right) ,$ then $X$ and $Y$ are free over $A_{k},$ in $%
\left( A_{j+1},\,E_{k+1,\,j+1}\right) .$
\end{proposition}

\strut

\begin{proof}
Now take the following (sub)tower of the given tower

\strut

\begin{center}
$\Bbb{C}\subset ^{E_{k}}A_{k}\subset ^{E_{k+1,\,\,j+1}}A_{j+1}.$
\end{center}

\strut

Then, by the previous theorems, $\left( A_{j+1},\,E_{j+1}\right) $ and $%
\left( A_{j+1},\,E_{k+1,\,j+1}\right) $ are compatible. By hypothesis, for
any $p_{1},...,p_{n}\in X\cup Y,$ we have that

\strut

$\ k_{n}^{(E_{k+1,\,j+1})}\left( p_{1},...,p_{n}\right)
=c_{E_{k+1,\,j+1}}^{(n)}\left( p_{1}\otimes b_{2}p_{2}\otimes ...\otimes
b_{n}p_{n}\right) $

\strut

where $b_{2},...,b_{n}\in B$ are arbitrary

\strut

(3.5.1) $\ \ \ \ \ \ \ \ \ \ \ =\left( E_{j+1}(b_{2})\cdot \cdot \cdot
E_{j+1}(b_{n})\right) \cdot k_{n}^{(E_{j+1})}\left( p_{1},...,p_{n}\right)
\cdot 1_{B},$

\strut

by the previous lemma, where $k_{n}^{(E_{j+1})}(...)$ is the scalar-valued
cumulant, in the sense of Speicher and Nica. By the freeness of $X$ and $Y$
in $\left( A_{j+1},\,E_{j+1}\right) ,$ (3.5.1) goes to

\strut

$\ \ \ \ \ \ \ =\left\{ 
\begin{array}{lll}
\left( E_{j+1}(b_{2})\cdot \cdot \cdot E_{j+1}(b_{n})\right) \cdot
k_{n}^{(E_{j+1})}\left( x_{i_{1}},...,x_{i_{n}}\right) \cdot 1_{B} &  & 
\text{or} \\ 
&  &  \\ 
\left( E_{j+1}(b_{2})\cdot \cdot \cdot E_{j+1}(b_{n})\right) \cdot
k_{n}^{(E_{j+1})}\left( y_{i_{1}},...,y_{i_{n}}\right) \cdot 1_{B} &  & 
\end{array}
\right. $

\strut

$\ \ \ \ \ \ \ =\left\{ 
\begin{array}{lll}
k_{n}^{(\varphi )}(x_{i_{1}},...,x_{i_{n}}) &  & \text{or} \\ 
&  &  \\ 
k_{n}^{(\varphi )}(y_{i_{1}},...,y_{i_{n}}), &  & 
\end{array}
\right. $

\strut

for all $(i_{1},...,i_{n})\in \{1,...,s\}^{n},$ $n\in \Bbb{N}.$ Therefore, $%
X $ and $Y$ are free over $A_{k},$ in $\left( A_{j+1},\,E_{k+1,\,j+1}\right)
.$
\end{proof}

\strut

\begin{proposition}
Let $B$ be a unital algebra and suppose that we have a tower of amalgamated
NCPSpaces as before. Fix $k<j$ in $\Bbb{N}.$ Let $X=\{x_{1},...,x_{s}\}$ and 
$Y=\{y_{1},...,y_{s}\}$ be two subsets of operators in $A_{j+1}.$ Assume
that $X\cup Y$ is free from $A_{k}$ in $\left( A_{j+1},\,E_{j+1}\right) .$
If $X$ and $Y$ are free in $\left( A_{j+1},\,E_{j+1}\right) ,$ then

\strut 

(1) $R_{x_{1},...,x_{s},y_{1},...,y_{s}}^{(E_{k+1,\,j%
\,})}(z_{1},...,z_{s})=R_{x_{1},...,x_{s}}^{(E_{k+1,%
\,j})}(z_{1},...,z_{s})+R_{y_{1},...,y_{s}}^{(E_{k+1,%
\,j})}(z_{s+1},...,z_{2s}).$

\strut 

(2) $R_{x_{1}+y_{1},...,x_{s}+y_{s}}^{(E_{k+1,\,j})}(z_{1},...,z_{s})=\left(
R_{x_{1},...,x_{s}}^{(E_{k+1,\,j})}+R_{y_{1},...,y_{s}}^{(E_{k+1,\,j})}%
\right) (z_{1,}...,z_{s}).$

\strut 

(3) $R_{x_{1}y_{1},...,x_{s}y_{s}}^{(E_{k+1,\,j})}(z_{1},...,z_{s})=\left(
R_{x_{1},...,x_{s}}^{(E_{k+1,\,j})}\,\,\,\frame{*}_{A_{k}}\,\,%
\,R_{y_{1},...,y_{s}}^{(E_{k+1,\,\,j})}\right) (z_{1},...,z_{s}).$

\strut 

(Recall that $E_{k+1,\,k+1}=\varphi _{k+1}:A_{k+1}\rightarrow A_{k},$ for
all $k.$)
\end{proposition}

\strut

\begin{proof}
Remember that $\left( A_{j+1},\,E_{j+1}\right) $ and $\left(
A_{j+1},\,E_{k+1,\,j}\right) $ are compatible, for any $k<j$ in $\Bbb{N}.$
Since $X\cup Y$ and $A_{k}$ are free in $\left( A_{j+1},\,E_{j+1}\right) ,$
if $X$ and $Y$ are free in $\left( A_{j+1},\,E_{j+1}\right) ,$ then $X$ and $%
Y$ are free over $A_{k},$ in $\left( A_{j+1},\,E_{k+1,\,j\,}\right) ,$ by
the previous proposition.
\end{proof}

\strut

\strut

\strut \strut

\section{Ladders of Amalgamated Noncommutative Probability Spaces}

\strut

\strut \strut

\strut

In this chapter, we will consider new algebraic structure, so-called a
ladder of amalgamated NCPSpaces. Suppose that we have two towers of
amalgamated NCPSpaces,

\strut

\begin{center}
$\Bbb{C}\subset ^{\varphi _{0}}B\subset ^{\varphi _{1}}A_{1}\subset
^{\varphi _{2}}A_{2}\subset ^{\varphi _{3}}A_{3}\subset ^{\varphi
_{4}}A_{4}\subset ^{\varphi _{5}}\cdot \cdot \cdot $
\end{center}

and

\begin{center}
$\Bbb{C}\subset ^{\varphi _{0}^{\prime }}B^{\prime }\subset ^{\varphi
_{1}^{\prime }}A_{1}^{\prime }\subset ^{\varphi _{2}^{\prime }}A_{2}^{\prime
}\subset ^{\varphi _{3}^{\prime }}A_{3}^{\prime }\subset ^{\varphi
_{4}^{\prime }}A_{4}^{\prime }$ $\subset ^{\varphi _{5}^{\prime }}\cdot
\cdot \cdot .$
\end{center}

\strut

Assume that there is a relation between this two towers ;

\strut

\begin{center}
$A_{j}^{\prime }\subset A_{j}$ is a subalgebra, for all $j\in \Bbb{N}\cup
\{0\},$
\end{center}

\strut

with $1_{A_{j}^{\prime }}=1_{A_{j}},$ where $A_{0}^{\prime }=B^{\prime }$
and $A_{0}=B.$ Then, for any $k,$ we have the following square as a part of
the ladder,

\strut

\begin{center}
$
\begin{array}{lll}
A_{k} & \subset ^{E_{k+1,\,j}} & A_{j} \\ 
\cup ^{i_{k}} &  & \cup ^{i_{j}} \\ 
A_{k}^{\prime } & \subset ^{E_{k+1,\,j}^{\prime }} & A_{j}^{\prime },
\end{array}
$
\end{center}

\strut

where $E_{k+1,\,j}^{\prime }=\varphi _{k+1}...\varphi _{j}\mid
_{A_{j}^{\prime }}:A_{j}^{\prime }\rightarrow A_{k}^{\prime }$ is a
conditional expectation and $i_{k}:A_{k}\rightarrow A_{k}^{\prime },$ $%
i_{j}:A_{j}\rightarrow A_{j}^{\prime }$ are conditional expectations. i.e,
we have the following ladder generated by those two towers of amalgamated
NCPSpaces,

\strut

\begin{center}
$
\begin{array}{lllllllllll}
\Bbb{C} & \subset ^{\varphi _{0}} & B & \subset ^{\varphi _{1}} & A_{1} & 
\subset ^{\varphi _{2}} & A_{2} & \subset ^{\varphi _{3}} & A_{3} & \subset
^{\varphi _{4}} & \cdot \cdot \cdot \\ 
&  & \cup ^{i_{0}} &  & \cup ^{i_{1}} &  & \cup ^{i_{2}} &  & \cup ^{i_{3}}
&  & \cdot \cdot \cdot \\ 
\Bbb{C} & \subset ^{\varphi _{0}^{\prime }} & B^{\prime } & \subset
^{\varphi _{1}^{\prime }} & A_{1}^{\prime } & \subset ^{\varphi _{2}^{\prime
}} & A_{2}^{\prime } & \subset ^{\varphi _{3}^{\prime }} & A_{3}^{\prime } & 
\subset ^{\varphi _{4}^{\prime }} & \cdot \cdot \cdot
\end{array}
$
\end{center}

\strut

Notice that $i_{k}\varphi _{k+1}:A_{k+1}\rightarrow A_{k}^{\prime }$ is a
conditional expectation, for all $k\in \Bbb{N}.$ Indeed,

\strut

\begin{center}
$i_{k}\varphi _{k+1}(a)=i_{k}(a)=a,$ for all $a\in A_{k}^{\prime },$
\end{center}

\strut

since $a\in A_{k}^{\prime }\subset A_{k}.$ Also,

\strut

\begin{center}
$i_{k}\varphi _{k+1}\left( axa^{\prime }\right) =i_{k}\left( a\varphi
_{k+1}(x)a^{\prime }\right) =a\left( i_{k}\varphi _{k+1}(x)\right) a^{\prime
},$
\end{center}

\strut

for all $a,a^{\prime }\in A_{k}^{\prime }$ and $x\in A_{k+1}.$

\strut

\begin{definition}
We say that the following inclusions

\strut 

\begin{center}
$
\begin{array}{lllllllllll}
\Bbb{C} & \subset ^{\varphi _{0}} & B & \subset ^{\varphi _{1}} & A_{1} & 
\subset ^{\varphi _{2}} & A_{2} & \subset ^{\varphi _{3}} & A_{3} & \subset
^{\varphi _{4}} & \cdot \cdot \cdot  \\ 
&  & \cup ^{i_{0}} &  & \cup ^{i_{1}} &  & \cup ^{i_{2}} &  & \cup ^{i_{3}}
&  & \cdot \cdot \cdot  \\ 
\Bbb{C} & \subset ^{\varphi _{0}^{\prime }} & B^{\prime } & \subset
^{\varphi _{1}^{\prime }} & A_{1}^{\prime } & \subset ^{\varphi _{2}^{\prime
}} & A_{2}^{\prime } & \subset ^{\varphi _{3}^{\prime }} & A_{3}^{\prime } & 
\subset ^{\varphi _{4}^{\prime }} & \cdot \cdot \cdot 
\end{array}
$
\end{center}

\strut 

is a commuting ladder of amalgamated NCPSpaces if conditions (i) and (ii)
are satisfied ;

\strut 

(i) Two towers of inclusions

\strut 

\begin{center}
$
\begin{array}{lllllllllll}
\Bbb{C} & \subset ^{\varphi _{0}} & B & \subset ^{\varphi _{1}} & A_{1} & 
\subset ^{\varphi _{2}} & A_{2} & \subset ^{\varphi _{3}} & A_{3} & \subset
^{\varphi _{4}} & \cdot \cdot \cdot , \\ 
&  &  &  &  &  &  &  &  &  &  \\ 
\Bbb{C} & \subset ^{\varphi _{0}^{\prime }} & B^{\prime } & \subset
^{\varphi _{1}^{\prime }} & A_{1}^{\prime } & \subset ^{\varphi _{2}^{\prime
}} & A_{2}^{\prime } & \subset ^{\varphi _{3}^{\prime }} & A_{3}^{\prime } & 
\subset ^{\varphi _{4}^{\prime }} & \cdot \cdot \cdot 
\end{array}
$
\end{center}

\strut 

are towers of amalgamated NCPSpaces.

\strut 

(ii) For any choice of a square of amalgamated NCPSaces,

\strut 

\begin{center}
$
\begin{array}{lll}
A_{k} & \subset ^{\varphi _{k+1}} & A_{k+1} \\ 
\cup ^{i_{k}} &  & \cup ^{i_{k+1}} \\ 
A_{k}^{\prime } & \subset ^{\varphi _{k+1}^{\prime }} & A_{k+1}^{\prime }
\end{array}
$
\end{center}

\strut 

is a commuting square of amalgamated NCPSpaces, i.e conditional expectations
on that square satisfy that

\strut \strut 

\begin{center}
$i_{k}\varphi _{k+1}=\varphi _{k+1}^{\prime }i_{k+1}:A_{k+1}\rightarrow
A_{k}^{\prime }.$
\end{center}
\end{definition}

\strut

\begin{example}
Let $B$ be a unital algebra and $(A,\varphi ),$ a NCPSpace over $B.$ Suppose
that we have a linear functional $\varphi _{0}:B\rightarrow \Bbb{C}.$ Also
assume that there is a scalar-valued chain of projections $%
(p_{k})_{k=1}^{\infty }\subset A_{pro}.$ Suppose that there is a NCPSpace
over $B,$ $(D,\psi )$ such that $D\subset A$ and there is a coditional
expectation $E_{D}^{A}:A\rightarrow D$ defined by $E_{D}^{A}(x)=x,$ if $x\in
D\subset A$ and $E_{D}^{A}(x)=0_{D},$ otherwise. Then we can get a commuting
ladder of (compressed) amalgamated NCPSpaces,

\strut 

\begin{center}
$
\begin{array}{lllllllll}
\Bbb{C} & \subset ^{\varphi _{0}} & B & \subset ^{\varphi _{1}} & p_{1}Ap_{1}
& \subset ^{\varphi _{2}} & p_{2}Ap_{2} & \subset ^{\varphi _{3}} & \cdot
\cdot \cdot  \\ 
&  & \cup ^{E_{D}^{A}} &  & \text{ \ }\cup ^{E_{D}^{A}} &  & \text{ \ }\cup
^{E_{D}^{A}} &  &  \\ 
\Bbb{C} & \subset ^{\varphi _{0}} & B & \subset ^{\varphi _{1}} & p_{1}Dp_{1}
& \subset ^{\varphi _{2}} & p_{2}Dp_{2} & \subset ^{\varphi _{3}} & \cdot
\cdot \cdot ,
\end{array}
$
\end{center}

\strut 

where $E_{D}^{A}=E_{D}^{A}\mid _{p_{j}Ap_{j}}$ and $\varphi _{j}$'s are
defined in Section 2.3 (i.e, $\varphi _{j+1}(p_{j+1}ap_{j+1})=p_{j}ap_{j},$ $%
\forall a\in A$), for all $j\in \Bbb{N}.$ Indeed, for any $j,$ we have a
commuting square of amalgamated NCPSpaces,

\strut 

\begin{center}
$
\begin{array}{lll}
p_{j}Ap_{j} & \subset ^{\varphi _{j+1}} & p_{j+1}Ap_{j+1} \\ 
\,\,\,\cup ^{E_{D}^{A}} &  & \,\,\,\,\,\,\,\,\cup ^{E_{D}^{A}} \\ 
p_{j}Dp_{j} & \subset ^{\varphi _{j+1}} & p_{j+1}Dp_{j+1},
\end{array}
$
\end{center}

\strut 

since

\strut 

\begin{center}
$
\begin{array}{ll}
E_{D}^{A}\,\,\varphi _{j+1}\left( p_{j+1}ap_{j+1}\right)  & =E_{D}^{A}\left(
p_{j}ap_{j}\right)  \\ 
&  \\ 
& =\left\{ 
\begin{array}{lll}
p_{j}ap_{j} &  & \text{if }a\in D \\ 
&  &  \\ 
0_{D} &  & \text{otherwise.}
\end{array}
\right. 
\end{array}
$
\end{center}

and

\strut 

\begin{center}
$
\begin{array}{ll}
\varphi _{j+1}\,\,E_{D}^{A}\left( p_{j+1}ap_{j+1}\right)  & =\left\{ 
\begin{array}{lll}
\varphi _{j+1}\left( p_{j+1}ap_{j+1}\right)  &  & \text{if }a\in D \\ 
&  &  \\ 
\varphi _{j+1}\left( 0_{D}\right)  &  & \text{otherwise}
\end{array}
\right.  \\ 
&  \\ 
& =\left\{ 
\begin{array}{lll}
p_{j}ap_{j} &  & \text{if }a\in D \\ 
&  &  \\ 
0_{D} &  & \text{otherwise.}
\end{array}
\right. 
\end{array}
$
\end{center}

\strut 

Hence the given ladder satisfies the condition (ii) of Definition 4.1. Since
two towers in the ladder are towers of compressed amalgamated NCPSPaces, our
ladder satisfies the condition (i) of Definition 4.1. Thus the ladder is a
commuting ladder of amalgamated NCPSpaces.
\end{example}

\strut \strut

\begin{proposition}
Let

\strut 

\begin{center}
$
\begin{array}{lllllllllll}
\Bbb{C} & \subset ^{\varphi _{0}} & B & \subset ^{\varphi _{1}} & A_{1} & 
\subset ^{\varphi _{2}} & A_{2} & \subset ^{\varphi _{3}} & A_{3} & \subset
^{\varphi _{4}} & \cdot \cdot \cdot  \\ 
&  & \cup ^{i_{0}} &  & \cup ^{i_{1}} &  & \cup ^{i_{2}} &  & \cup ^{i_{3}}
&  & \cdot \cdot \cdot  \\ 
\Bbb{C} & \subset ^{\varphi _{0}^{\prime }} & B^{\prime } & \subset
^{\varphi _{1}^{\prime }} & A_{1}^{\prime } & \subset ^{\varphi _{2}^{\prime
}} & A_{2}^{\prime } & \subset ^{\varphi _{3}^{\prime }} & A_{3}^{\prime } & 
\subset ^{\varphi _{4}^{\prime }} & \cdot \cdot \cdot 
\end{array}
$
\end{center}

\strut 

be a commuting ladder of amalgamated NCPSpaces. Then, for any $k<j$ in $\Bbb{%
N},$ a square of inclusions

\strut 

\begin{center}
$
\begin{array}{lll}
A_{k} & \subset ^{E_{k+1,\,j\,}} & A_{j} \\ 
\cup ^{i_{k}} &  & \cup ^{i_{j}} \\ 
A_{k}^{\prime } & \subset ^{E_{k+1,\,j}^{\prime }} & A_{j}^{\prime }
\end{array}
$
\end{center}

\strut 

is a commuting square of amalgamated NCPSpaces.
\end{proposition}

\strut

\begin{proof}
Fix $k<j$. By the definition of commuting squares of amalgamated NCPSpaces,
it suffices to show that

\strut

\begin{center}
$i_{k}E_{k+1,\,j}=E_{k+1,\,j}^{\prime }i_{j}:A_{j}\rightarrow A_{k}^{\prime
} $
\end{center}

\strut

as a conditional expectation on $A_{j}$ over $A_{k}$. Observe that, for any $%
x\in A_{j},$

\strut

\begin{center}
$
\begin{array}{ll}
i_{k}E_{k+1,\,j}(x) & =i_{k}\varphi _{k+1}...\varphi _{j}(x)=\left(
i_{k}\varphi _{k+1}\right) \varphi _{k+2}...\varphi _{j}(x) \\ 
& =\left( \varphi _{k+1}^{\prime }i_{k+1}\right) \varphi _{k+2}...\varphi
_{j}(x) \\ 
& =\varphi _{k+1}^{\prime }\left( i_{k+1}\varphi _{k+2}\right) \varphi
_{k+3}...\varphi _{j}(x) \\ 
& =\varphi _{k+1}^{\prime }\left( \varphi _{k+2}^{\prime }i_{k+2}\right)
\varphi _{k+3}...\varphi _{j}(x) \\ 
& =...=\varphi _{k+1}^{\prime }...\varphi _{j-1}^{\prime }\left(
i_{j-1}\varphi _{j}\right) (x). \\ 
& =\varphi _{k+1}^{\prime }...\varphi _{j-1}^{\prime }\left( \varphi
_{j}^{\prime }i_{j}\right) (x)=\left( \varphi _{k+1}^{\prime }...\varphi
_{j}^{\prime }\right) i_{j}(x) \\ 
& =E_{k+1,\,j}^{\prime }i_{j}(x).
\end{array}
$
\end{center}

\strut

Therefore, the square of inclusions

\strut

\begin{center}
$
\begin{array}{lll}
A_{k} & \subset ^{E_{k+1,\,j\,}} & A_{j} \\ 
\cup ^{i_{k}} &  & \cup ^{i_{j}} \\ 
A_{k}^{\prime } & \subset ^{E_{k+1,\,j}^{\prime }} & A_{j}^{\prime }
\end{array}
$
\end{center}

\strut

is a commuting square of amalgamated NCPSpaces.
\end{proof}

\strut

\begin{corollary}
Let

\strut 

\begin{center}
$
\begin{array}{lllllllllll}
\Bbb{C} & \subset ^{\varphi _{0}} & B & \subset ^{\varphi _{1}} & A_{1} & 
\subset ^{\varphi _{2}} & A_{2} & \subset ^{\varphi _{3}} & A_{3} & \subset
^{\varphi _{4}} & \cdot \cdot \cdot  \\ 
&  & \cup ^{i_{0}} &  & \cup ^{i_{1}} &  & \cup ^{i_{2}} &  & \cup ^{i_{3}}
&  & \cdot \cdot \cdot  \\ 
\Bbb{C} & \subset ^{\varphi _{0}^{\prime }} & B^{\prime } & \subset
^{\varphi _{1}^{\prime }} & A_{1}^{\prime } & \subset ^{\varphi _{2}^{\prime
}} & A_{2}^{\prime } & \subset ^{\varphi _{3}^{\prime }} & A_{3}^{\prime } & 
\subset ^{\varphi _{4}^{\prime }} & \cdot \cdot \cdot 
\end{array}
$
\end{center}

\strut 

be a commuting ladder of amalgamated NCPSpaces. Let $(j_{k})_{k=1}^{\infty }$
be any (sub)sequence of $\Bbb{N}.$ Then we can construct a commuting
(sub)ladder of amalgamated NCPSpaces

\strut 

(4.2.1)\strut 

\begin{center}
$
\begin{array}{lllllllll}
\Bbb{C} & \subset ^{E_{j_{1}}} & A_{j_{1}} & \subset ^{E_{j_{1}+1,\,j_{2}}}
& A_{j_{2}} & \subset ^{E_{j_{2}+1,\,j_{3}}} & A_{j_{3}} & \subset  & \cdot
\cdot \cdot  \\ 
&  & \cup ^{i_{j_{1}}} &  & \cup ^{i_{j_{2}}} &  & \cup ^{i_{j_{3}}} &  & 
\cdot \cdot \cdot  \\ 
\Bbb{C} & \subset ^{E_{j_{1}}^{\prime }} & A_{j_{1}}^{\prime } & \subset
^{E_{j_{1}+1,\,j_{2}}^{\prime }} & A_{j_{2}}^{\prime } & \subset
^{E_{j_{2}+1,\,j_{3}}^{\prime }} & A_{j_{3}}^{\prime } & \subset  & \cdot
\cdot \cdot .
\end{array}
$
\end{center}
\end{corollary}

\strut

\begin{proof}
By the previous proposition, we know that, for each chosen $j_{p}$ ( $p\in 
\Bbb{N}$),

\strut

\begin{center}
$
\begin{array}{lll}
A_{j_{p}} & \subset ^{E_{j_{p}+1,\,j_{p+1}}} & A_{j_{p+1}} \\ 
\cup ^{i_{j_{p}}} &  & \cup ^{i_{j_{p+1}}} \\ 
A_{j_{p}}^{\prime } & \subset ^{E_{j_{p}+1,\,j_{p+1}}^{\prime }} & 
A_{j_{p+1}}^{\prime }
\end{array}
$
\end{center}

\strut

is a commuting square of amalgamated NCPSpaces. So the ladder (4.2.1)
satisfies the condition (ii) of Definition 4.1. By the previous chapter (or
by compatibility), we know that the given two towers are all towers of
amalgamated NCPSpaces. So the condition (i) of Definition 4.1 is also
satisfied. Thus the ladder (4.2.1) is a commuting ladder of amalgamated
NCPSpaces.
\end{proof}

\strut

The above proposition and corollary show that if we have a \textbf{commuting}
ladder of amalgamated NCPSpaces, then we can naturally define conditional
expectations from the algebras, $A_{j}$'s, in the upper tower of the ladder
to the algebras, $A_{k}^{\prime }$ 's, in the lower tower of the ladder,
whenever $k<j$. The commuting property is crucial. Hence we can observe the
R-transform theory for those conditional expectations. Also, the above
corollary shows us that if we have a commuting ladder of amalgamated
NCPSpaces, then we can choose infinitely many commuting subladders of
amalgamated NCPSpaces from the given ladder.

\strut

\strut We observed, in the previous example, that we have the following
commuting ladder of amalgamated (compressed) NCPSpaces ; Let $(A,\varphi )$
be a NCPSpace over $B$ and let $D$ be a subalgebra of $A$ containing $B$.
Assume that we have a conditional expectation $E_{D}^{A}:A\rightarrow D$
defined by $E_{B}^{A}(a)=a$ if $a\in D$ and $E_{D}^{A}(a)=0_{D},$ otherwise.
Now, let $(p_{k})_{k=1}^{\infty }\subset A_{pro}$ be a scalar-valued chain
of projections such that $\varphi (p_{k})=\alpha _{k},$ $\forall k\in \Bbb{N}%
.$ Then we have the commuting ladder of amalgamated NCPSpaces,\strut

\strut \strut

(*)

\begin{center}
$
\begin{array}{lllllllll}
\Bbb{C} & \subset ^{\varphi _{0}} & B & \subset ^{\varphi _{1}} & p_{1}Ap_{1}
& \subset ^{\varphi _{2}} & p_{2}Ap_{2} & \subset ^{\varphi _{3}} & \cdot
\cdot \cdot \\ 
&  & \cup ^{E_{D}^{A}} &  & \text{ \ }\cup ^{E_{D}^{A}} &  & \text{ \ }\cup
^{E_{D}^{A}} &  &  \\ 
\Bbb{C} & \subset ^{\varphi _{0}} & B & \subset ^{\varphi _{1}} & p_{1}Dp_{1}
& \subset ^{\varphi _{2}} & p_{2}Dp_{2} & \subset ^{\varphi _{3}} & \cdot
\cdot \cdot .
\end{array}
$
\end{center}

\strut

We will call this ladder, a commuting ladder of amalgamated compressed
NCPSpaces induced by $(p_{k})_{k=1}^{\infty }$ and $E_{D}^{A}.$ Remark that
it is possible that $D=p_{N}Ap_{N},$ for some $N\in \Bbb{N}.$ Then the above
ladder can be regarded as

\strut

\begin{center}
$
\begin{array}{lllllllll}
\Bbb{C} & \subset ^{\varphi _{0}} & B & \subset ^{\varphi _{1}} & p_{1}Ap_{1}
& \subset ^{\varphi _{2}}\cdot \cdot \cdot \subset ^{\varphi _{N}} & 
p_{N}Ap_{N} & \subset & \cdot \cdot \cdot \\ 
&  & \cup ^{E_{D}^{A}} &  & \text{ \ }\cup ^{E_{D}^{A}} &  & \text{ \ }%
\shortparallel ^{E_{D}^{A}} &  &  \\ 
\Bbb{C} & \subset ^{\varphi _{0}} & B & \subset ^{\varphi _{1}} & p_{1}Dp_{1}
& \subset ^{\varphi _{2}}\cdot \cdot \cdot \subset ^{\varphi _{N}} & \text{
\ \ }D & = & \cdot \cdot \cdot .
\end{array}
$
\end{center}

\strut \strut

i.e, we take $p_{N}Dp_{N}=p_{N+1}Dp_{N+1}=p_{N+2}Dp_{N+2}=...=D$. In \ this
case, we have that

\strut

\begin{center}
$\varphi _{N+j}\mid _{p_{N+j}Dp_{N+j}}:p_{N+j}Dp_{N+j}\rightarrow
p_{N+(j-1)}Dp_{N+(j-1)}$
\end{center}

\strut \strut

as the identity map, $i_{D}:D\rightarrow D,$ for all $j\in \Bbb{N}.$

\strut

\begin{proposition}
Let $B$ be a unital algebra and $(A,\varphi ),$ a NCPSpace over $B$ and let $%
D$ be a subalgebra of $A$ containing $B$. Suppose that we have a
scaalr-valued chain of projections $(p_{k})_{k=1}^{\infty }\subset A_{pro},$
with $\varphi (p_{k})=\alpha _{k}\cdot 1_{B},$ $\forall k,$ and a
conditional expectation $E_{D}^{A}:A\rightarrow D$ defined by $%
E_{D}^{A}(x)=x,$ if $x\in D$ and $E_{D}^{A}(x)=0_{D}=0_{B},$ otherwise.
Assume that we have a commuting ladder of amalgamated compressed NCPSpaces
(*). Then, for any $k<j\in \Bbb{N},$ we can define a conditional expectation 
$F_{k+1,\,\,j}:p_{j}Ap_{j}\rightarrow p_{k}Dp_{k}$ satisfying the following
property ; let $X$ and $Y$ be subsets of $A$. If $X\cup Y$ and $\{p_{k}\}$
are free over $B,$ in $(A,\varphi ),$ then $X$ and $Y$ are free over $%
p_{k}Ap_{k},$ in $(A,\varphi )$ if and only if $p_{j}Xp_{j}$ and $p_{j}Yp_{j}
$ are free over $p_{k}Dp_{k},$ in $\left( p_{j}Ap_{j},\,\frac{1}{\alpha
_{k+1}}F_{k+1,\,j}\right) .$
\end{proposition}

\strut

\begin{proof}
Suppose that we have a ladder of amalgamated compressed NCPSpaces (*)
induced by $(p_{k})_{k=1}^{\infty }\subset A_{pro}$ and a canonical
conditional expectation, $E_{D}^{A}:A\rightarrow D.$ Then we can define a
conditional expectation,

\strut

\begin{center}
$F_{k+1,\,j}:p_{j}Ap_{j}\rightarrow p_{k}Dp_{k},$
\end{center}

\strut

for any $k<j,$ by

\strut

\begin{center}
$F_{k+1,\,j}=E_{D}^{A}\left( \frac{1}{\alpha _{k+1}}\cdot E_{k+1,\,j}\right)
.$
\end{center}

\strut

Recall that if $X$ and $Y$ are free over $p_{k}Ap_{k},$ in $(A,\varphi ),$
then $p_{j}Xp_{j}$ and $p_{j}Yp_{j}$ are free over $p_{k}Ap_{k},$ in $\left(
p_{j}Ap_{j},\,\frac{1}{\alpha _{k+1}}E_{k+1,\,j}\right) $ \ (See Section
2.3). Also, since

\strut

\begin{center}
$E_{D}^{A}E_{k+1,\,j}=E_{k+1,\,j}\,E_{D}^{A},$

\strut
\end{center}

we need to show that $p_{j}\left( E_{D}^{A}(X)\right) p_{j}$ and $%
p_{j}\left( E_{D}^{A}(Y)\right) p_{j}$ are free over $p_{k}Dp_{k}.$ But,
trivially, if $X$ and $Y$ are free over $p_{k}Ap_{k},$ then $E_{D}^{A}(X)$
and $E_{D}^{A}(Y)$ are free over $p_{k}Dp_{k},$ by the very definition of $%
E_{D}^{A}.$ i.e,

\strut

\begin{center}
$E_{D}^{A}(X)=\{x\in D:x\in D\cap X\}$ and $E_{D}^{A}(Y)=\{y\in D:y\in D\cap
Y\}.$
\end{center}

\strut

Thus they are automaticlly free under the hypothesis. So, we have that, if $%
X $ and $Y$ are free over $p_{k}Ap_{k},$ in $\left( A,\varphi \right) ,$ then

\strut

(i) $\ p_{j}Xp_{j}$ and $p_{j}Yp_{j}$ are free over $p_{k}Ap_{k},$ in $%
\left( p_{j}Ap_{j},\,\frac{1}{\alpha _{k+1}}E_{k+1,\,j}\right) $

\strut

(ii) By (i), $p_{j}E_{D}^{A}(X)p_{j}$ and $p_{j}E_{D}^{A}(Y)p_{j}$ are also
free over $p_{k}Dp_{k}$ in $\left( p_{j}Ap_{j},\,\frac{1}{\alpha _{k+1}}%
E_{D}^{A}E_{k+1,\,j}\right) .$

\strut

So, we can conclude that $p_{k}Xp_{k}$ and $p_{k}Yp_{k}$ are free over $%
p_{k}Dp_{k},$ in $\left( p_{j}Ap_{j},\,F_{k+1,\,j}\right) .$
\end{proof}

\strut

As an application, we have the following theorem ;

\strut

\begin{theorem}
Let $B$ be a unital algebra and $(A,\varphi ),$ a NCPSpace over $B$ and let $%
D$ be a subalgebra of $A$ containing $B$. Suppose that we have a
scaalr-valued chain of projections $(p_{k})_{k=1}^{\infty }\subset A_{pro},$
with $\varphi (p_{k})=\alpha _{k}\cdot 1_{B},$ $\forall k,$ and a
conditional expectation $E_{D}^{A}:A\rightarrow D$ defined by $%
E_{D}^{A}(x)=x,$ if $x\in D$ and $E_{D}^{A}(x)=0_{D}=0_{B},$ otherwise.
Assume that we have a commuting ladder of amalgamated compressed NCPSpaces
(*). Let $F_{k+1,\,j}:p_{j}Ap_{j}\rightarrow p_{k}Ap_{k}$ be $%
E_{D}^{A}E_{k+1,\,j}.$ Assume that we have two $p_{k}Ap_{k}$-free subsets $%
X=\{p_{j}x_{1}p_{j},...,p_{j}x_{s}p_{j}\}$ and $Y=%
\{p_{j}y_{1}p_{j},...,p_{j}y_{s}p_{j}\}.$ If $X\cup Y$ and $\{p_{j}\}$ is
free over $p_{k}Ap_{k},$ in $(p_{j}Ap_{j},E_{k+1,\,j}),$ then, for any fixed 
$\ k<j\in \Bbb{N},$ we have that

\strut 

(1) $%
R_{p_{j}x_{1}p_{j},...,p_{j}x_{s}p_{j},p_{j}y_{1}p_{j},...,p_{j}y_{s}p_{j}}^{(F_{k+1,\,j})}(z_{1},...,z_{2s})
$

\begin{center}
$=R_{p_{j}x_{1}p_{j},...,p_{j}x_{s}p_{j}}^{(F_{k+1,\,\,j})}$\strut $%
(z_{1},...,z_{s})+R_{p_{j}yp_{j},...,p_{j}y_{s}p_{j}}^{(F_{k+1,%
\,j})}(z_{s+1},...,z_{2s}).$
\end{center}

\strut 

(2) $%
R_{p_{j}x_{1}p_{j}+p_{j}y_{1}p_{j},...,p_{j}x_{s}p_{j}+p_{j}y_{s}p_{j}}^{(F_{k+1,\,j})}(z_{1},...,z_{2s})
$

\begin{center}
$=R_{p_{j}x_{1}p_{j},...,p_{j}x_{s}p_{j}}^{(F_{k+1,\,\,j})}$\strut $%
(z_{1},...,z_{s})+R_{p_{j}yp_{j},...,p_{j}y_{s}p_{j}}^{(F_{k+1,%
\,j})}(z_{1},...,z_{s}).$
\end{center}

\strut \strut 

(3) $%
R_{(p_{j}x_{1}p_{j})(p_{j}y_{1}p_{j}),...,(p_{j}x_{s}p_{j})(p_{j}y_{s}p_{j})}^{(F_{k+1,\,j})}(z_{1},...,z_{s})
$

\begin{center}
$=\left( R_{p_{j}x_{1}p_{j},...,p_{j}x_{s}p_{j}}^{(F_{k+1,\,j})}\,\,\,\frame{%
*}_{B}\,\,\,R_{p_{j}y_{1}p_{j},...,p_{j}y_{s}p_{j}}^{(F_{k+1,j})\,\,:\,\,t}%
\right) (z_{1},...,z_{s})$

$=\left( R_{p_{j}x_{1}p_{j},...,p_{j}x_{s}p_{j}}^{(F_{k+1,\,j})}\,\,\,\,\,%
\frame{*}_{B}\,\,\,\,\,\,\,R_{y_{1}\,\,,\,\,\,.\,\,\,.\,\,\,.\,\,%
\,,y_{s}}^{(\varphi )\,\,:\,\,t}\right) \,(z_{1},...,z_{s}).$
\end{center}

$\square $
\end{theorem}

\strut

We can extend this construction of such ladders from two pairs of towers of
amalgamated NCPSpaces to nets of amalgamated NCPSpaces such as

\strut

\begin{center}
$
\begin{array}{lllllllllll}
\Bbb{C} & \subset & B & \subset & A_{1} & \subset & A_{2} & \subset & A_{3}
& \subset & \cdot \cdot \cdot \\ 
&  & \cup &  & \cup &  & \cup &  & \cup &  &  \\ 
\Bbb{C} & \subset & B^{\prime } & \subset & A_{1}^{\prime } & \subset & 
A_{2}^{\prime } & \subset & A_{3}^{\prime } & \subset & \cdot \cdot \cdot \\ 
&  & \cup &  & \cup &  & \cup &  & \cup &  &  \\ 
\Bbb{C} & \subset & B^{\prime \prime } & \subset & A_{1}^{\prime \prime } & 
\subset & A_{2}^{\prime \prime } & \subset & A_{3}^{\prime \prime } & \subset
& \cdot \cdot \cdot \\ 
&  & \cup &  & \cup &  & \cup &  & \cup &  &  \\ 
&  & \,\,\vdots &  & \,\,\vdots &  & \,\,\vdots &  & \,\,\,\vdots &  & \ddots
\end{array}
$
\end{center}

\strut \strut

But even in the commuting ladders, we have seen that it is difficult to
consider the amalgamated freeness between two algebras in that ladderWe only
considered the case when we have a commuting ladder of amalgamated
compressed NCPSpaces induced by a scalar-valued chain of projections and the
given ''good'' conditional expectation.

\strut

\strut

\strut

\strut \textbf{References}

\strut 

\strut 

{\small [1] \ \ A. Nica, R-transform in Free Probability, IHP course note,
available at www.math.uwaterloo.ca/\symbol{126}anica.\strut }

{\small [2] \ \ A. Nica, R-transforms of Free Joint Distributions and
Non-crossing Partitions, J. of Func. Anal, 135 (1996), 271-296.}

{\small [7] \ \ A. Nica, D. Shlyakhtenko and R. Speicher, R-cyclic Families
of Matrices in Free Probability, J. of Funct Anal, 188 (2002),
227-271.\strut }

{\small [4] \ \ A. Nica, D. Shlyakhtenko and R. Speicher, R-diagonal
Elements and Freeness with Amalgamation, Canad. J. Math. Vol 53, Num 2,
(2001) 355-381.}

{\small [5] \ \ D. Shlyakhtenko, Some Applications of Freeness with
Amalgamation, J. Reine Angew. Math, 500 (1998), 191-212.\strut }

{\small [6] \ \ D. Shlyakhtenko, A-Valued Semicircular Systems, J. of Funct
Anal, 166 (1999), 1-47.\strut }

{\small [7] \ \ D. Voiculescu, Operations on Certain Non-commuting
Operator-Valued Random Variables, Ast\'{e}risque, 232 (1995), 243-275.\strut 
}

{\small [8] \ \ D.Voiculescu, K. Dykemma and A. Nica, Free Random Variables,
CRM Monograph Series Vol 1 (1992).\strut }

{\small [9] \ \ I. Cho, Amalgamated Boxed Convolution and Amalgamated
R-transform Theory (preprint).\strut }

{\small [10] I. Cho, Compressed Amalgamated R-transform Theory,
preprint.\strut }

{\small [11] I. Cho, Perturbed R-transform Theory, preprint.\strut }

{\small [12] I. Cho, Compatibility of a noncommutative probability space and
an amalgamated noncommutative probability space, preprint}

{\small [13] \ R. Speicher, Combinatorics of Free Probability Theory IHP
course note, available at www.mast.queensu.ca/\symbol{126}speicher.\strut }

{\small [14] R. Speicher, Combinatorial Theory of the Free Product with
Amalgamation and Operator-Valued Free Probability Theory, AMS Mem, Vol 132 ,
Num 627 , (1998).}

\end{document}